\pdfoutput=1
\documentclass[journal,onecolumn,draftclsnofoot]{IEEEtran}
\usepackage{amsmath,amssymb}
\usepackage{enumitem}
\usepackage{mathabx}
\usepackage{multirow}
\usepackage{algorithm,algorithmic}
\usepackage{MnSymbol}
\usepackage{graphicx,epstopdf} 
\usepackage[justification=centering,caption=false]{subfig} 
\usepackage[justification=justified,font=footnotesize]{caption}
\usepackage{color}
\usepackage{array}
\usepackage{tikz}
\usepackage[export]{adjustbox}
\usepackage{scrextend}
\usepackage{bbm,float}

\usetikzlibrary{arrows.meta,calc,decorations.markings,math,arrows.meta}
\providecommand{\N}{\ensuremath\mathbb N}
\providecommand{\R}{\ensuremath \mathbb R}
\providecommand{\spt}{\ensuremath \text{spt}}

\newtheorem{theorem}{Theorem}
\newtheorem{corollary}[theorem]{Corollary}
\newtheorem{lem}[theorem]{Lemma}
\newtheorem{assum}[theorem]{Assumption}
\newtheorem{definition}[theorem]{Definition}

\newcommand{\cref}[1]{\ref{#1}}

\newcommand{\ip}[1]{\langle {#1} \rangle}
\newcommand{\Li}{\mathcal{L}_i}

\newcommand{\M}{\mathcal{M}}
\newcommand{\hle}{the Hybrid Liouville Equation }
\newcommand{\ocp}{HOCP}

\newcolumntype{C}[1]{>{\centering\arraybackslash\hspace{0pt}}p{#1}}

\providecommand{\Ram}[1]{\textcolor{blue}{RAM: #1}}

\providecommand{\Later}[1]{}
\providecommand{\Backup}[1]{}

\begin{document}
  \title{Optimal Control for Nonlinear Hybrid Systems via Convex Relaxations}
  \author{Pengcheng Zhao, Shankar Mohan, and Ram Vasudevan
  \thanks{P. Zhao and R. Vasudevan are with the Department of Mechanical Engineering, University of Michigan, Ann Arbor, MI 48109 USA {\tt\scriptsize \{pczhao,ramv\}@umich.edu}}
  \thanks{S. Mohan is with the Department of Electrical Engineering and Computer Science, University of Michigan, Ann Arbor, MI 48109 USA {\tt\scriptsize elemsn@umich.edu}}}
  \maketitle

\begin{abstract}
  This paper considers the optimal control for hybrid systems whose trajectories transition between distinct subsystems when state-dependent constraints are satisfied.
  Though this class of systems is useful while modeling a variety of physical systems undergoing contact, the construction of a numerical method for their optimal control has proven challenging due to the combinatorial nature of the state-dependent switching and the potential discontinuities that arise during switches.
  This paper constructs a convex relaxation-based approach to solve this optimal control problem.
  Our approach begins by formulating the problem in the space of relaxed controls, which gives rise to a linear program whose solution is proven to compute the globally optimal controller. 
  This conceptual program is solved by constructing a sequence of semidefinite programs whose solutions are proven to converge from below to the true solution of the original optimal control problem.
  Finally, a method to synthesize the optimal controller is developed. 
  Using an array of examples, the performance of the proposed method is validated on problems with known solutions and also compared to a commercial solver.
\end{abstract}
  \section{Introduction}

Let ${\mathcal I}$ be a finite set of labels and $U \subset \R^m$ a compact, convex, nonempty set. 
Let $X_i$ be a bounded, compact subset of $\R^{n_i}$ for some $n_i \in \N$.
Let $F_i: [0,T] \times X_i \times U \to \R^{n_i}$ be a vector field on $X_i$ for each $i \in {\mathcal I}$, which defines a \emph{controlled hybrid system}.
That is given an initial condition $x_0 \in X_j$ for some $j \in {\mathcal I}$, $T > 0$, and control input $u:[0,T] \times U$, the flow of the system satisfies the vector field $F_j$ almost everywhere until either the total time of evolution is $T$, the trajectory hits a guard, $G_{(j,j')} \subset X_j$ of the system for some $j' \in {\mathcal I}$, or the trajectory hits the boundary of $X_j$.
If the trajectory hits a guard, then the trajectory is re-initialized according to a reset map, $R_{(j,j')}: G_{(j,j')} \to X_{j'}$, and the flow proceeds from this new point according to the vector field $F_{j'}$ as described earlier
This notion of execution is formalized in Algorithm~\ref{alg:1}.

Let $H: \coprod_{i \in \mathcal I} X_i \rightarrow \R$ and $h:[0,T] \times \coprod_{i \in \mathcal I} X_i \times U \rightarrow \R$ represent the terminal and incremental cost, respectively, which can be distinct for each subsystem.
Let $X_{T_i} \subset X_i$ represent some terminal constraint for each $i \in {\mathcal I}$
Given an initial condition $x_0 \in X_j$, this paper is interested in solving the the following optimal control problem:
\begin{equation}
\label{eq:highlevel_opt}
\begin{aligned}
\inf \big\{ \int_0^T h\left(t,x(t),u(t)\right)\, dt + H\left(x(T)\right) ~ \big\mid ~ &u: [0,T] \to U,~ x: [0,T] \xrightarrow[u,~  x(0) = (x_0, j )]{Algorithm~\ref{alg:1}} \mathcal \coprod_{i \in \mathcal I} X_i, \\ & \text{ and }x(T) \in \coprod_{i \in \mathcal I} X_{T_i} \big\}.
\end{aligned}
\end{equation}
The optimization problem defined in \eqref{eq:highlevel_opt} is numerically challenging to solve since transitions occur between the distinct dynamical systems due to state-dependent constraints, which may give rise to discontinuities in the solution. 
This has rendered the immediate application of derivative-based algorithms to solve this problem impossible.
As a result, this problem has typically been numerically solved by fixing the sequence of transitions between subsystems.
To overcome the limitations of these existing numerical approaches, this paper first presents a conceptual, infinite dimensional linear program over measures, $(P)$, that is proven in Theorem~\ref{thm:P} to compute the global optimum of \eqref{eq:highlevel_opt}.
Subsequently, this paper describes a sequence of implementable, semidefinite programming-based relaxations, $(P_k)$, of the infinite dimensional linear program, which are proven in Theorem \ref{thm:P_k} to converge monotonically from below to the solution to \eqref{eq:highlevel_opt}.

\subsection{Related Work}

Controlled hybrid dynamical systems have been used to describe the dynamics of a variety of physical systems in which the evolution of the system undergoes sudden changes due to the satisfaction of state-dependent constraints such as in bipeds \cite{Westervelt2007}, automotive sub-systems~\cite{Heijden2007}, aircraft control~\cite{Soler2010}, and biological systems~\cite{Elowitz2000}.
Given the practical applications of such systems, the development of algorithms to perform optimal control of hybrid systems has drawn considerable interest amongst both theoreticians and practitioners. 

The theoretical development of both necessary and sufficient conditions for the optimal control of hybrid controlled systems has been considered using extensions of the Pontryagin Maximum Principle~\cite{passenberg2011minimum,shaikh2007hybrid,sussmann1999maximum} and Dynamic Programming~\cite{branicky1998unified,dharmatti2005hybrid,schollig2007hybrid}, respectively.
Recent work has even linked these pair of theoretical approaches for hybrid systems~\cite{pakniyat2014relation}.
Typically these theoretical developments have focused their attention to systems where the sequence of transitions between the systems have been known \textit{a priori}.
As a result, practitioners have typically fixed the sequence of transitions and applied gradient-based methods to perform local optimization over the time spent and control applied within each subsystem~\cite{griffin2015walking,hereid20163d,smit2017energetic,Westervelt2003}. 

More recent work has focused on the development of numerical optimal control techniques for mechanical systems undergoing contact without specifying the ordering of visited subsystems.
For mechanical systems, state-dependent switching arises due to the effect of unilateral constraints.
One approach to address the optimal control problem has focused on the construction of a novel notion of derivative~\cite{pace2016piecewise}. 
Though this method still requires fixing the total number of visited subsystems, assuming \textit{a priori} knowledge of the visited subsystems, and performs optimization only over the initial condition, this gradient-based approach is able to find the locally optimal ordering of subsystems under certain regularity conditions on the nature of contact.
Other approaches have relaxed satisfaction of the unilateral constraint directly and instead focused on treating constraint satisfaction as a continuous decision variable that can be optimized using traditional numerical methods to find local minima~\cite{posa2014direct,westenbroek2015optimal,yunt2005trajectory}.

This paper focuses on developing a numerical approach to find the global optimal to the hybrid optimal control problem. 
Our method relies on treating the optimal control problem in the relaxed sense wherein the original problem is lifted to the space of measures~\cite{Bhatt1996,Kurtz1998}. 
In the instance of classical dynamical systems, this lifting renders the optimal control problem linear in the space of relaxed controls~\cite{berkovitz2013optimal}; however, there were few numerical methods to tackle this relaxed problem directly.

Recent developments in algebraic geometry have made it possible to solve this lifted optimal control problem for classical dynamical systems by relying on moment-based relaxations~\cite{lasserre2008nonlinear}. 
By solving the problem over truncated moment sequences, it is possible to transform the optimal control problem into either a finite-dimensional linear or finite-dimensional semidefinite program. 
Either transformation of the relaxed problem is proven to provide a lower bound on the optimal cost.
In fact, this bound converges to the true optimal cost as the moment sequence extends to infinity under the assumption that the incremental cost is convex in control. 
Recent work has also shown how the optimal control policy can be extracted for systems that are affine in control~\cite{korda2016controller,Zhao2016}. 
Unfortunately this relaxed control formulation for controlled hybrid systems, the subsequent development of a numerically implementable convex relaxation, and optimal control synthesis have remained unaddressed.

Note that the focus of this paper is on the development of an optimal control approach for hybrid systems with state-dependent rather than controlled switching. 
In particular, after state-dependent switching, the state is allowed to change in a discontinuous manner.
This is typically not allowed for systems with just controlled switching.
For this class of systems with controlled switching, there have been a variety of numerical methods proposed to perform optimal control~\cite{bengea2005optimal,claeys2016modal,egerstedt2006transition,johnson2011second,vasudevan2013consistent,vasudevan2013consistent2,wardi2015switched}.

\subsection{Contributions and Organization}

The contributions of this paper are three-fold: first, Section \ref{sec:abs_prob} provides a conceptual infinite dimensional linear programming-based approach that is proven to solve \eqref{eq:highlevel_opt}; second, Section~\ref{sec:implementation} presents a numerically implementable semidefinite programming-based sequence of relaxations to this infinite dimensional linear program that is proven to generate a sequence of convergent lower bounds to the true optimal cost of \eqref{eq:highlevel_opt}; finally Section~\ref{sec:implementation} provides a method to generate a sequence of controllers that converge to the true optimal control of ~\eqref{eq:highlevel_opt}.

The remainder of this paper is organized as follows: Section~\ref{sec:prelim} describes in detail the class of systems under consideration and defines their executions, Section~\ref{sec:abstraction} describes how to lift executions of the hybrid system to the space of measures, Section~\ref{sec:extensions} describes how to extend the optimal control approach to free final time problems, and Section~\ref{sec:examples} illustrates the efficacy of the proposed method on a variety of systems. 

  \section{Preliminaries}

\label{sec:prelim}
This section introduces the notation used throughout the remainder of this paper, define controlled hybrid systems, and formulate the optimal control problem of interest. 
This makes substantial use of measure theory, and the unfamiliar reader may refer to \cite{bogachev2007measure} for an introduction.

\subsection{Notation}
\label{subsec:notation}
Given an element $y \in \R^{n}$, let $[y]_i$ denote the $i$-th component of $y$. 
We use the same convention for elements belonging to any multidimensional vector space. 
Let $\R[y]$ denote the ring of real polynomials in the variable $y$.
Let $\R_k[y]$ denote the space of real valued multivariate polynomials of total degree less than or equal to $k$. 

Let $\{A_i\}_{i\in {\cal I}}$ be a family of non-empty sets indexed by $i$, the \emph{disjoint union} of this family is $\coprod_{i \in {\cal I}}A_i = \bigcup_{i \in {\cal I}} (A_i \times \{i\})$.
Let $\iota_i : A_i \to \coprod_{i \in {\cal I}}A_i$ be the canonical injection defined as $\iota_i(x) = (x,i)$ whose inverse we denote by $\pi_i: \coprod_{i \in {\cal I}}A_i \to A_i \cup \{\emptyset\}$. 
Note that $\pi_i(\iota_j(x)) = \emptyset$ if $i \neq j$.
For the sake of convenience, we denote $x_i := \pi_i(x)$ for the remainder of the paper. 
We similarly define a projection operator onto the indexing set $\lambda: \coprod_{i \in {\cal I}}A_i \to {\cal I}$ such that $\lambda(\iota_i(x)) = i$.


For any set $S$, we denote by $\mathbbm{1}_S$ the indicator function on $S$.
Suppose $Y$ is a metric space, then let $C(Y)$ be the space of continuous functions on $Y$,
let $C_b(Y)$ be the space of bounded continuous functions on $Y$,
let $AC(I)$ be the space of absolutely continuous functions on an interval $I \in \R$,
let $L^1(Y)$ be the space of  $L^1$ functions on $Y$,
let $W^{1,\infty}(Y)$ be the $(1,\infty)$-Sobolev space on $Y$, 
and let $\M(Y)$ be the space of signed Radon measures on $Y$ endowed with the total variation norm (denoted by $\| \cdot \|$), whose positive cone $\M_+(K)$ is the space of unsigned Radon measures on $Y$.
Any measure $\mu \in \M(Y)$ can be viewed as an element of the dual space to $C(Y)$ via the duality pairing
\begin{equation}
    \langle \mu, v \rangle := \int_Y v(z) \mu(z), \quad \forall v \in C(Y).
\end{equation}
For any measure $\mu \in \M(Y)$, let the \emph{support} of $\mu$ be denoted as $\spt(\mu)$.
A probability measure is a non-negative, unsigned measure whose integral is one.
Denote the dual to a vector space $V$ as $V'$.

Suppose $Y$ is a metric space and $Y_1 \subset Y$ is a compact subset with subspace topology, then for any measurable function $f \in L^1(Y_1)$, we define its \emph{zero extension} onto $Y$, denoted as $\hat{f} \in L^1(Y)$, as
\begin{equation}
    \hat{f}(y) =
    \begin{cases}
        f(y), & \text{if } y \in Y_1\\
        0, & \text{if } y \in Y \backslash Y_1
    \end{cases}
\end{equation}
For any measure $\mu \in \M(Y_1)$, we define its \emph{zero extension} onto $Y$, denoted as $\hat{\mu} \in \M(Y)$, as
\begin{equation}
    \hat{\mu}(B) = \mu( B \cap Y_1 )
\end{equation}
for all subsets $B$ in the Borel $\sigma$-algebra of $Y$.
For any measure $\mu \in \M( Y_1 \times Y_2 )$ and variables $(y_1, y_2) \in Y_1 \times Y_2$, we denote by $\mu_{y_1 | y_2} \in \M(Y_1)$ the conditional probability measure of $\mu$ on $Y_1$ given an instance of $y_2 \in Y_2$, and we denote by $\mu_{y_2} \in M(Y_2)$ the marginal of $\mu$ over the space $Y_2$.

For any measure $\mu \in \M_+(\R^n)$ and measurable function $\theta \in L^1(\R^n)$, we define the \emph{convolution} of $\mu$ and $\theta$, denoted as $\mu * \theta$, as
\begin{equation}
    (\mu * \theta)(B) = \int_{\R^n} \int_{\R^n} \mathbbm{1}_B( x + y ) \theta(y) \, dy \, d\mu(x)
\end{equation}
for all subsets $B$ in the Borel $\sigma$-algebra of $\R^n$.
If $Y_1$, $Y_2$ are measurable spaces, $\mu \in \M(Y_1)$, and $f: Y_1 \to Y_2$ is a Borel function, we denote by $f_\# \mu \in \M(Y_2)$ the \emph{pushforward} of $\mu$ through $f$, given by
\begin{equation}
    (f_\# \mu)(B) := \mu( f^{-1}(B) )
\end{equation}
for any $B$ in the Borel $\sigma$-algebra of $Y_2$. Therefore for every $f_\# \mu$-integrable function $v: Y_2 \to \R$, we have:
\begin{equation}
    \int_{Y_2} v \, d(f_\# \mu) = \int_{Y_1}v \circ f \, d\mu
\end{equation}

\subsection{Controlled Hybrid Systems}

Motivated by \cite{burden2015metrization}, we define the class of \emph{controlled hybrid systems} of interest in the remainder of this paper:

\begin{definition}
  \label{def:hybrid_system}
  A \emph{controlled hybrid system} is a tuple ${\cal H} = ({\cal I}, {\cal E}, {\cal D}, U, {\cal F}, {\cal S}, {\cal R})$, where:
  \begin{itemize}
  \item ${\cal I}$ is a finite set indexing the discrete states of ${\cal H}$;
  \item ${\cal E} \subset {\cal I} \times {\cal I}$ is a set of edges, forming a directed graph structure over ${\cal I}$;
  \item ${\cal D} = \coprod_{i\in {\cal I}} X_i$ is a disjoint union of domains, where each $X_i$ is a compact subset of $\R^{n_i}$, and $n_i \in \N$;
  \item $U$ is a compact subset of $\R^m$ that describes the range of control inputs, where $m \in \N$;
  \item ${\cal F} = \{F_i\}_{i \in {\cal I}}$ is the set of vector fields, where each $F_i : \R \times X_i \times U \to \R^{n_i}$ is the vector field defining the dynamics of the system on $X_i$;
  \item ${\cal S} = \coprod_{e \in {\cal E}}S_e$ is a disjoint union of guards, where each $S_{(i,i')} \subset \partial X_i$ is a compact, co-dimension 1 guard defining a state-dependent transition to $X_{i'}$; and,
  \item ${\cal R} = \{R_e\}_{e \in {\cal E}}$ is a set of reset maps, where each map $R_{(i,i')}\colon S_{(i,i')} \to X_{i'}$ defines the transition from guard $S_{(i,i')}$ to $X_{i'}$.
  \end{itemize}
\end{definition}
For convenience, we refer to controlled hybrid systems as just hybrid systems, and we refer to a vertex within the graph structure associated with a controlled hybrid system as a \emph{mode}.

Even though the range space of control inputs are assumed to be the same in each mode, this is not restrictive since we can always concatenate all the control inputs in different modes.
The compactness of each $X_i$ ensures the optimization problem defined in Section 2.3 is well-posed.
To avoid any ambiguity during transitions, we make the following assumptions:
\begin{assum}
\label{assum:g_intersect}
Guards do not intersect themselves, i.e.,
\begin{equation}
S_e \cap S_{e'} = \emptyset, \quad \forall e,e' \in {\cal E},
\end{equation}
and guards do not intersect with the images of reset maps, i.e.,
\begin{equation}
S_e \cap R_{e'}(S_{e'}) = \emptyset, \quad \forall e,e' \in {\cal E}.
\end{equation}
\end{assum}

\begin{assum}
\label{assum:no_scuffing}
The vector field $F_i$ has nonzero normal component on the boundary of $X_i$.
\end{assum}
\noindent Assumption \cref{assum:g_intersect} ensures at most one transition can be executed at a time, and Assumption \cref{assum:no_scuffing} implies transition always happens when a trajectory reaches a guard (i.e., no grazing).

Next, we define an \emph{execution} of a hybrid system up to time $T>0$ via construction in Algorithm \cref{alg:1}. 
Step 1 initializes the execution at a given point $(x_0,i)$ at time $t=0$.
Step 3 defines $\phi$ to be the maximal integral curve of $F_i$ under the control $u$ beginning from the initial point.
Step 4 defines the execution on a finite interval as the curve $\phi$ with associated index $i$. 
As described in Steps 5 - 7, the trajectory terminates when it either reaches the terminal time $T$ or hits $\partial X_i \backslash \bigcup_{(i,i')\in {\cal E}} S_{(i,i')}$ where no transition is defined.
Steps 8 and 9 define a discrete transition to a new domain using a reset map where evolution continues again as a classical dynamical system by returning to Step 3.
We denote the space of such executions as $\cal X$.
Note that for any execution $\gamma$ and any $t \in [0,T]$, we have $\gamma(t) = (\gamma_i(t), \lambda(\gamma(t)))$.

\begin{algorithm}[t]
		  \caption{Execution of Hybrid System ${\cal H}$}
		  \label{alg:1}
		  \begin{algorithmic}[1]
		    \REQUIRE $t = 0$, $T>0$, $i \in {\cal I}$, $(x_0,i) \in {\cal D}$, and $u: \R \to U$ Lebesgue measurable.
		    \STATE Set $\gamma(0) = (x_0,i)$.
		    \LOOP
		    \STATE Let $\phi: I \to X_i$ be an absolutely continuous function such that:
		\begin{enumerate}[label=(\roman*)]
			\item $\dot{\phi}(s)  = F_i(s,\phi(s),u(s))$ for almost every $s \in I$ with respect to the Lebesgue measure on $I \subset [t,T]$ with $(\phi(t),i) = x(t)$ and
			\item for any other $\hat{\phi}\colon \hat{I} \to X_i$ satisfying (i), $\hat{I} \subset I$.
			\end{enumerate}
		    \STATE Let $t' = \sup I$ and $\gamma(s) = (\phi(s),i)$ for each $s \in [t,t')$.
		    \IF{$t' = T$, \textbf{or} $\nexists (i,i') \in {\cal E}$ such that $\phi(t') \in S_{(i,i')}$}
		    \STATE Stop.
		    \ENDIF
		    \STATE Let $(i,i') \in {\cal E}$ be such that $\phi(t') \in S_{(i,i')}$.
		    \STATE Set $\gamma(t') = R_{(i,i')}\left(\phi(t')\right)$, $t = t'$, and $i = i'$. \label{algo:hexec_reset}
		    \ENDLOOP
		  \end{algorithmic}
\end{algorithm}

Hybrid systems can suffer from Zeno executions, i.e. executions that undergo an infinite number of discrete transitions in a finite amount of time.
Since the state of the trajectory after the Zeno occurs may not be well defined \cite{ames2006there}, we do not consider systems with Zeno executions:
\begin{assum}
${\cal H}$ has no Zeno execution.
\end{assum}

\subsection{Problem Formulation}

This paper is interested in finding a pair $(\gamma,u)$ satisfying Algorithm \cref{alg:1} with a given initial condition $x_0$, that reaches a user-specified target set while minimizing a user-specified cost function.
To formulate this problem, we first define \emph{the target set}, $X_T \subset {\cal D}$, as:
\begin{equation}
    X_T = \coprod_{i \in {\cal I}} X_{T_i},
\end{equation}
where $X_{T_i}$ is a compact subset of $X_i$ for each $i \in {\cal I}$.
To avoid any ambiguity, we put the following restriction on the target set:
\begin{assum}
\label{assum:target_set}
    The target set does not intersect with guards, i.e.,
    \begin{equation}
        X_{T_i} \cap S_{(i,i')} = \emptyset, \quad \forall (i,i') \in {\cal E}
    \end{equation}
\end{assum}


Next, we define the system trajectories and control actions of interest.
Given a real number $T>0$ and an initial point $(x_0,j) \in {\cal D}$, a pair of functions $(\gamma,u)$ satisfying Algorithm \cref{alg:1} is called an \emph{admissible pair} if $\gamma(T) \in X_T$. 
The trajectory in this instance is called an \emph{admissible trajectory} and the control input is called an \emph{admissible control}.
The time $T$ at which the admissible trajectory reaches the target set is called the \emph{terminal time}.
For convenience, we denote the space of admissible trajectories and controls by ${\cal X}_T$ and ${\cal U}_T$, respectively. 
The space of admissible pairs is denoted as ${\cal P}_T \subset {\cal X}_T \times {\cal U}_T$.
Without loss of generality, we assume the initial point does not belong to any guard:
\begin{assum}
\label{assum:x0}
The initial condition, $(x_0, i)$, does not belong to any guard, i.e.,
\begin{equation}
    x_0 \not\in S_{(i,i')}, \quad \forall (i,i') \in {\cal E}
\end{equation}
\end{assum}



For any admissible pair $(\gamma,u)$, the associated cost is defined as:
\begin{equation}
\label{eq:J}
J(\gamma,u) = \int_0^T h_{\lambda(\gamma(t))}\left(t,x_{\lambda(\gamma(t))}(t),u(t)\right)\, dt + H_{\lambda(\gamma(T))}\left(x_{\lambda(\gamma(T))}(T)\right)
\end{equation}
where $h_i:[0,T] \times \R^{n_i} \times \R^m \to \R$ and $H_i:\R^{n_i} \to \R$ are measurable functions for each $i \in {\cal I}$.

Our goal is to find an admissible pair that minimizes \eqref{eq:J}, which we refer to as \emph{Hybrid Optimal Control Problem (\ocp)}:
\begin{flalign}
\label{eq:OCP}
  &&\inf_{(\gamma,u) \in {\cal P}_T} &\phantom{4} \int_0^T h_{\lambda(\gamma(t))}\left(t,\gamma_{\lambda(\gamma(t))}(t),u(t)\right)\, dt + H_{\lambda(\gamma(T))}\left(x_{\lambda(\gamma(T))}(T)\right) &&(\ocp) \nonumber\\
  &&\text{s.t.} &\phantom{4} \gamma:[0,T]\to {\cal D} \text{ and } u:[0,T]\to U \text{  defined via Algorithm \cref{alg:1}} \nonumber\\
  &&& \phantom{4} \gamma(T) \in X_T. \nonumber
\end{flalign}
The optimal cost of $(\ocp)$ is defined as:
\begin{equation}
\label{eq:J*}
J^* = \inf_{(\gamma,u) \in {\cal P}_T} J(\gamma,u).
\end{equation}

  \section{The Hybrid Liouville Equation}
\label{sec:abstraction}
It is difficult to directly solve $(\ocp)$ over the space of admissible trajectories for various reasons.
First, the cost function and constraints may be nonlinear and non-convex, respectively, therefore global optimality is not guaranteed for gradient-descent methods.
Second, in the instance of hybrid systems, typically the sequence of transitions between modes must be specified.
To address these limitations, this section constructs measures whose supports model the evolution of families of trajectories, an equivalent form of $J$, and an equivalent form of Algorithm \cref{alg:1} in the space of measures.
These transformations make a convex formulation of $(\ocp)$ feasible.

To begin, consider again the projection $\gamma_i$ of an admissible trajectory $\gamma$ onto $X_i$.
Define the \emph{occupation measure} in mode $i \in {\cal I}$ associated with $\gamma$, denoted as $\mu^i(\cdot \mid \gamma) \in \M_+([0,T] \times X_i)$, as:
\begin{equation}
\label{eq:def:mu_traj}
\mu^i (A \times B \mid \gamma) := \int_0^T \mathbbm{1}_{A\times B}(t,\gamma_i(t))\, dt
\end{equation}
for all subsets $A \times B$ in the Borel $\sigma$-algebra of $[0,T]\times X_i$.
Note that $x_i(t)$ may not be defined for all $t \in [0,T]$, but we use the same notation and let $\mathbbm{1}_{A \times B}(t,\gamma_i(t)) = 0$ whenever $\gamma_i(t)$ is undefined.
The quantity $\mu^i(A \times B \mid \gamma )$ is equal to the amount of time the graph of the trajectory, $(t,\gamma_i(t))$, spends in $A \times B$.

Similarly, define the \emph{initial measure}, $\mu_0^i(\cdot \mid \gamma) \in \M_+(X_i)$, as:
\begin{equation}
\label{eq:def:mu0_traj}
\mu_0^i( B \mid \gamma ) := \mathbbm{1}_B(\gamma_i(0)),
\end{equation}
for all subsets $B$ in the Borel $\sigma$-algebra of $X_i$; 
define the \emph{terminal measure}, $\mu_T^i( \cdot \mid \gamma) \in \M_+(X_{T_i})$, as:
\begin{equation}
\label{eq:def:muT_traj}
\mu_T^i( B \mid \gamma ) := \mathbbm{1}_B(\gamma_i(T)),
\end{equation}
for all subsets $B$ in the Borel $\sigma$-algebra of $X_{T_i}$.
Note in this instance, we have abused notation and the reader should not confuse $\mu^i_0(\cdot \mid \gamma)$ and $\mu^i_T(\cdot \mid \gamma)$ with marginals of $\mu^i(\cdot \mid \gamma)$ evaluated at $t = 0$ and $t=T$. 
Finally, define the \emph{guard measure}, $\mu^{S_{(i,i')}}( \cdot \mid \gamma ) \in \M_+([0,T] \times S_{(i,i')})$, as:
\begin{equation}
    \mu^{S_{(i,i')}}(A \times B \mid \gamma) := \#\{t \in A \mid \lim_{\tau \to t^-}\gamma_i(\tau) \in B\}
\end{equation}
for all subsets $A \times B$ in the Borel $\sigma$-algebra of $[0,T]\times S_{(i,i')}$, given any pair $(i,i') \in {\cal E}$. 
The guard measure counts the number of times a given execution passes through the guard. 


If the admissible control $u$ associated with $\gamma$ according to Algorithm \cref{alg:1} is also given, define the occupation measure in mode $i \in {\cal I}$ associated with the pair $(\gamma,u)$, denoted as $\mu^i(\cdot \mid \gamma,u) \in \M_+([0,T] \times X_i \times U)$, as:
\begin{equation}
\label{eq:def:mu}
\mu^i(A \times B \times C \mid \gamma,u) := \int_0^T \mathbbm{1}_{A \times B \times C} (t,\gamma_i(t),u(t))\, dt
\end{equation}
for all subsets $A \times B \times C$ in the Borel $\sigma$-algebra of $[0,T] \times X_i \times U$.
For notational convenience, it is useful to collect the initial, average, final, and guard occupation measures in each mode into a single distinct object. 
That is, define $\mu_0^{\cal I}( \cdot \mid, \gamma, u ) \in \M_+(\mathcal{D})$ as $\mu^{\cal I}_0( \cdot,i \mid, \gamma, u ) := \mu^i_0( \cdot \mid \gamma, u)$ for each $i \in {\cal I}$. 
For the sake of convenience, we refer to $\mu_0^{\cal I}$ as an initial measure and write $\mu_0^i$ when we refer to the $i$-th slice of $\mu_0^{\cal I}$.
We define $\mu^{\cal I}( \cdot \mid \gamma, u ) \in \M_+([0,T]\times \mathcal{D} \times U)$, $\mu_T^{\cal I}( \cdot \mid \gamma, u ) \in \M_+(X_T)$, and $\mu^{\cal S}( \cdot \mid \gamma, u ) \in \M_+([0,T] \times \mathcal{S})$ in a similar fashion and refer to them in a similar way.

Using the notion of occupation measure and terminal measure, we now rewrite the cost function $J$ as a linear function on measures:
\begin{lem}
Let $\mu^{\cal I}(\cdot \mid \gamma,u)$ and $\mu_T^{\cal I}(\cdot \mid \gamma)$ be the occupation measure and terminal measure associated with the pair $(\gamma,u)$, respectively. 
Then the cost function can be expressed as:
\begin{equation}
J(\gamma,u) = \sum_{i\in {\cal I}}\langle \mu^i(\cdot \mid \gamma,u), h_i \rangle + \sum_{i\in {\cal I}} \langle \mu_T^i(\cdot \mid \gamma), H_i \rangle
\end{equation}
\end{lem}
\begin{IEEEproof}
Notice that $h_i$ and $H_i$ are measurable, and the rest follows directly from \eqref{eq:J}, \eqref{eq:def:muT_traj}, and \eqref{eq:def:mu}.
\end{IEEEproof}

Despite the cost function potentially being a nonlinear function for the admissible pair in the space of functions, the analogous cost function over the space of measures is linear.
A similar analogue holds true for the dynamics of the system.
That is, the occupation measure associated with an admissible pair satisfies a linear equation over measures.
To formulate this linear equation over measures, let $\Li: C^1\left( [0,T]\times X_i \right) \rightarrow C\left( [0,T]\times X_i \times U \right)$ be a linear operator which acts on a test function $v$ as:
\begin{equation}
\label{eq:Li}
(\Li v)(t,x,u) = \frac{\partial v(t,x)}{\partial t} + \sum_{k=1}^{n_i} \frac{\partial v(t,x)}{\partial x_k}[F_i(t,x,u)]_k, \quad \forall i \in {\cal I}
\end{equation}
\Later{$\partial [x]_k$??}
Using the dual relationship between measures and functions, we define $\Li': C([0,T]\times X_i \times U)' \to C^1([0,T]\times X_i)'$ as the adjoint operator of $\Li$, satisfying:
\begin{equation}
\langle \Li' \mu , v \rangle = \langle \mu, \Li v \rangle
\end{equation}
for all $\mu \in \M([0,T] \times X_i \times U)$ and $v \in C^1([0,T]\times X_i)$.

\Later{Need to explain why $\M_+ \subset C'$.}
Each of these adjoint operators can describe the evolution of trajectories of the system within each mode \cite{lasserre2008nonlinear}.
However in the instance of hybrid systems trajectories may not just begin evolution within a mode at $t=0$.
Instead a trajectory can enter mode $i \in {\cal I}$ either by starting from some point in $X_i$ at $t = 0$, or by hitting a guard $S_{(i',i)}$ for some $(i',i) \in {\cal E}$ and subsequently transitioning to a point in $X_i$.
Similarly a trajectory can terminate in mode $i \in {\cal I}$ either by reaching terminal time $T$, or by hitting a guard $S_{(i,i')}$ for some $(i,i') \in {\cal E}$ and transitioning away. 
For notational convenience we modify reset maps to also act on time, namely, define $\tilde{R}_{(i,i')} : [0,T] \times S_{(i,i')} \to [0,T] \times X_{i'}$ by
\begin{equation}
    \tilde{R}_{(i,i')}(t,x) = ( t, R_{(i,i')}(x) )
\end{equation}
for all $(i,i') \in {\cal E}$ and $(t,x) \in [0,T] \times S_{(i,i')}$.
To describe trajectories of a controlled hybrid system using measures we rely on the following result:
\begin{lem}
\label{lem:traj2hle}
Given an admissible pair $(\gamma,u)$, its initial measure, occupation measure, terminal measure, and guard measure, satisfy the following linear equation over measures:
\begin{equation}
\label{eq:hle|x}
\delta_0 \otimes \mu_0^i(\cdot \mid \gamma) + \Li' \mu^i(\cdot \mid \gamma,u) + \sum_{(i',i) \in {\cal E}} \tilde{R}_{(i', i)\#} \mu^{S_{(i',i)}}(\cdot \mid \gamma) = \delta_T \otimes \mu_T^i(\cdot \mid \gamma) + \sum_{(i,i') \in {\cal E}} \mu^{S_{(i,i')}}(\cdot \mid \gamma), \quad \forall i \in {\cal I}
\end{equation}
where the linear operator equation \eqref{eq:hle|x} holds in the sense that:
\begin{equation}
\begin{aligned}
\langle \mu_0^i(\cdot \mid \gamma), v(0,\cdot) \rangle & + \langle \mu^i(\cdot \mid \gamma,u), \Li v \rangle + \\ & + \sum_{(i',i)\in {\cal E}} \langle \mu^{S_{(i',i)}}(\cdot \mid \gamma), v(\cdot, R_{(i',i)}(\cdot)) \rangle
= \langle \mu_T^i(\cdot \mid \gamma), v(T,\cdot) \rangle + \sum_{(i,i') \in {\cal E}} \langle \mu^{S_{(i,i')}}(\cdot \mid \gamma), v \rangle
\end{aligned}
\end{equation}
for all test functions $v \in C^1([0,T] \times X_i)$.
\end{lem}
\begin{IEEEproof}
\Later{Not quite right...}
This lemma is a restatement of Equation (16) in \cite{Shia2014}.
\end{IEEEproof}

Now one can ask whether the converse relationship holds: 
do measures that satisfy \eqref{eq:hle|x} always coincide with trajectories generated by Algorithm \cref{alg:1}? 
More explicitly, do an arbitrary set of measures, $\mu_0^{\cal I} \in \M_+(\mathcal{D})$, $\mu^{\cal I} \in \M_+([0,T]\times \mathcal{D} \times U)$, $\mu^{\cal I}_T \in  \M_+(X_T)$, and $\mu^{\cal S} \in \M_+([0,T] \times \mathcal{S})$, that satisfy \eqref{eq:hle|x} correspond to an initial measure, $\mu_0^{\cal I}( \cdot \mid \gamma )$, occupation measure, $\mu^{\cal I}(\cdot \mid \gamma, u )$, terminal measure, $\mu^{\cal I}_T(\cdot \mid \gamma )$, and guard measure, $\mu^{\cal S} (\cdot \mid \gamma )$?
To answer this question, we first consider a family of admissible trajectories modeled by a non-negative measure $\rho \in \M_+({\cal X}_T)$, and define an \emph{average occupation measure} $\zeta^i \in \M_+([0,T] \times X_i)$ in each mode $i \in {\cal I}$ for the family of trajectories as:
\begin{equation}
\label{eq:def:zeta}
\zeta^i(A \times B) := \int_{{\cal X}_T} \mu^i(A \times B \mid \gamma)\, d\rho(\gamma),
\end{equation}
for any $i \in {\cal I}$; an \emph{average terminal measure}, $\zeta_T^i \in \M_+(X_T)$, by
\begin{equation}
\label{eq:def:zeta_T}
\zeta_T^i( B ) := \int_{{\cal X}_T} \mu_T^i(B \mid \gamma)\, d\rho(\gamma),
\end{equation}
for any $i \in {\cal I}$; and an \emph{average guard measure}, $\zeta^{S_{(i,i')}} \in \M_+([0,T] \times S_{(i,i')})$, by
\begin{equation}
\label{eq:def:zetaS}
\zeta^{S_{(i,i')}}(A \times B) := \int_{{\cal X}_T} \mu^{S_{(i,i')}}(A \times B \mid \gamma)\, d\rho(\gamma)
\end{equation}
for any $(i,i') \in {\cal E}$.

To prove the converse of Lemma \ref{lem:traj2hle}, we define the \emph{Hybrid Liouville Equation} whose solution, as we establish next, can be disintegrated into a set of measures that we eventually prove are related to $\rho$:

\begin{lem}
\label{lem:disintegration}
Let $\mu_0^{\cal I} \in \M_+(\mathcal{D})$, $\mu^{\cal I} \in \M_+([0,T]\times \mathcal{D} \times U)$, $\mu^{\cal I}_T \in  \M_+(X_T)$, and $\mu^{\cal S} \in \M_+([0,T] \times \mathcal{S})$ satisfy \emph{\hle}, which is defined as:
\begin{equation}
\label{eq:hle}
\delta_0 \otimes \mu_0^i + \Li' \mu^i + \sum_{(i',i) \in {\cal E}} \tilde{R}_{(i', i)\#} \mu^{S_{(i',i)}} = \delta_T \otimes \mu_T^i + \sum_{(i,i') \in {\cal E}} \mu^{S_{(i,i')}}
\end{equation}
for each $i \in {\cal I}$. Then each measure $\mu^i$ can be disintegrated as
\begin{equation}
\label{eq:disintegration}
d\mu^i(t,x,u) = d\nu^i_{u \mid t,x}( u ) \, d\mu^i_{t,x}(t,x) = d\nu^i_{u \mid t,x}(u) \, d\tilde{\mu}^i_{x \mid t}(x) \, dt
\end{equation}
where 
$\nu^i_{u \mid t,x}$ is a stochastic kernel on $U$ given $(t,x) \in [0,T] \times X_i$,
$\mu^i_{t,x}$ is the $(t,x)$-marginal of $\mu^i$,
and $\tilde{\mu}^i_{x \mid t}$ is a conditional measure on $X_i$ given $t \in [0,T]$.
\end{lem}

\begin{IEEEproof}
Since each measure $\mu^i$ is defined on a Euclidean space, which is Polish and therefore by definition Souslin, using \cite[Corollary~10.4.13]{bogachev2007measure}, it can be disintegrated as
\begin{equation}
    d\mu^i(t,x,u) = d\nu^i_{u \mid t,x}(u) \, d\mu^i_{t,x}(t,x)
\end{equation}
where $\nu^i_{u \mid t,x}$ is a stochastic kernel on $U$ given $(t,x) \in [0,T] \times X_i$,
and $\mu^i_{t,x}$ is the $(t,x)$-marginal of $\mu^i$. Using the same argument, we can further disintegrate $\mu^i_{t,x}$ into:
\begin{equation}
    d\mu^i_{t,x}(t,x) = d\mu^i_{x \mid t}(x) \, d\mu^i_t(t)
\end{equation}
where $\mu^i_t$ is the $t$-marginal of $\mu^i_{t,x}$.

Next, we show the measure $\mu^i_t$ is absolutely continuous with respect to Lebesgue measure.
Let $\psi \in C^1([0,T])$ be a test function of Equation \eqref{eq:hle}, we have
\begin{align}
         \ip{\delta_T \otimes \mu_T^i + \sum_{(i,i') \in {\cal E}} \mu^{S_{(i,i')}} - \delta_0 \otimes \mu_0^i - \sum_{(i',i) \in {\cal E}} \tilde{R}_{(i',i)\#} \mu^{S_{(i',i)}}, \psi}
        =& \ip{ \mu^i, \Li \psi } \label{eq:lemma7pf.1} \\
        =& \int_{[0,T]} \int_{X_i} \int_{U} \dot{\psi}(t) \, d\mu^i(t,x,u) \label{eq:lemma7pf.2}\\
        =& \int_{[0,T]} \dot{\psi}(t) \int_{X_i} \int_{U} \, d\mu^i(t,x,u) \label{eq:lemma7pf.3}\\
        =& \int_{[0,T]} \dot{\psi}(t) \, d\mu^i_t(t) \label{eq:lemma7pf.4}
\end{align}
where \eqref{eq:lemma7pf.1} and \eqref{eq:lemma7pf.2} follow by definition, \eqref{eq:lemma7pf.3} is from Fubini's theorem, and \eqref{eq:lemma7pf.4} is from the definition of $\mu^i_t$.
Then by applying the results in \cite[Exercise~5.8.78]{bogachev2007measure}, it follows that $\mu^i_t$ is absolutely continuous with respect to the Lebesgue measure.



Since $\mu^i$ is a Radon measure defined over a compact set and therefore $\sigma$-finite, its $t$-marginal measure, $\mu^i_t$, is also $\sigma$-finite. Using the Radon-Nikodym Theorem, there exists a function $l \in L^1([0,T])$ such that
\begin{equation}
    d\mu^i_t(t) = l(t) \, dt
\end{equation}

We now define $d\tilde{\mu}^i_{x \mid t} := l(t) \, d\mu^i_{x \mid t}$ for all $t \in [0,T]$, therefore
\begin{equation}
    d\mu^i_{x \mid t}(x) \, d\mu^i_t(t) = l(t) \, d\mu^i_{x \mid t}( x ) \, dt = d\tilde{\mu}^i_{x \mid t}(x) \, dt,
\end{equation}
and Equation \eqref{eq:disintegration} follows.
\end{IEEEproof}

For notational convenience, in the rest of this paper we abuse notation and denote $\tilde{\mu}^i_{x \mid t}$ as just $\mu^i_{x \mid t}$.

Using the disintegration \eqref{eq:disintegration}, we can rewrite \hle \eqref{eq:hle} as
\begin{equation}
\label{eq:hle2}
\begin{split}
    & \ip{ \delta_T \otimes \mu_T^i + \sum_{(i,i') \in {\cal E}} \mu^{S_{(i,i')}} - \delta_0 \otimes \mu_0^i - \sum_{(i',i) \in {\cal E}} \tilde{R}_{(i',i)\#} \mu^{S_{(i',i)}}, v } \\
    =& \int_{[0,T] \times X_i} \int_U \left( \frac{\partial v(t,x)}{\partial t} + \nabla_x v(t,x) \cdot F_i(t,x,u) \right) \, d\nu^i_{u \mid t,x}(u) \, d\mu^i_{t,x}(t,x) \\
    =& \int_{[0,T] \times X_i} \left( \frac{\partial v(t,x)}{\partial t} + \nabla_x v(t,x) \cdot \left( \int_U F_i(t,x,u) \, d\nu^i(u \mid t,x) \right) \right) \, d\mu^i_{t,x}(t,x) \\
    =& \int_{[0,T] \times X_i} \left( \frac{\partial v(t,x)}{\partial t} + \nabla_x v(t,x) \cdot \bar{F}_i(t,x) \right) \, d\mu^i_{t,x}(t,x)
\end{split}
\end{equation}
where
\begin{equation}
    \bar{F}_i(t,x) := \int_U F_i(t,x,u) \, d\nu^i_{u \mid t,x}(u) \in \text{conv } F_i(t,x,U)
\end{equation}
Here conv denotes the convex hull. 
Therefore we study the trajectories of the \emph{uncontrolled} hybrid system with vector fields ${\cal F} = \{\bar{F}_i\}_{i \in {\cal I}}$, and consider \hle in the form \eqref{eq:hle2}.

To further simplify the notation, we define:
\begin{equation}
\label{eq:sigma&eta}
    \begin{split}
        \sigma^i :=& \delta_0 \otimes \mu_0^i + \sum_{(i',i) \in {\cal E}} \tilde{R}_{(i',i)\#} \mu^{S_{(i',i)}}, \\
        \eta^i :=& \delta_T \otimes \mu_T^i + \sum_{(i,i') \in {\cal E}} \mu^{S_{(i,i')}}
    \end{split}
\end{equation}
and rewrite \hle as a non-homogeneous PDE,
\begin{equation}
    \label{eq:pde}
    \partial_t \mu^i_{t,x} + D_x \cdot (\bar{F}_i \mu^i_{t,x}) = \sigma^i - \eta^i,
\end{equation}
where \eqref{eq:pde} holds in the sense of distributions.
That is when we apply integration by parts, we can write:
\begin{equation}
    \label{eq:pde_explained}
    \int_0^T \int_{X_i} \left( \partial_t v(t,x) + \nabla_x v(t,x) \cdot \bar{F}_i \right)\, d\mu^i_{x \mid t}(x) \, dt + \int_{[0,T] \times X_i} v(t,x) \, d\sigma^i(t,x) = \int_{[0,T] \times X_i} v(t,x) \, d\eta^i(t,x)
\end{equation}
for any test function $v \in C^1([0,T] \times X_i)$.
We later show that $\sigma^i$ and $\eta^i$ capture the trajectories that enter and leave domain $i$, respectively.

\Later{mention that the HLE has at least one POSITIVE solution, citing another theorem yet to be proven.}

\Later{
\Ram{this needs to be moved}
\begin{lem}
For any $\sigma^i$, $\eta^i$, and $\bar{\mu}^i$ that satisfy \hle,
\begin{equation}
    \sigma^i([0,T]\times X_i) = \eta^i([0,T] \times X_i)
\end{equation}
\end{lem}
\begin{IEEEproof}
This follows directly by letting $v = 1$ in Equation \eqref{eq:pde_explained}.
\end{IEEEproof}
}

We establish next that $\mu^i_{x \mid t}$ is related to the solution of the ODE with dynamics $\bar{F}_i$. 
To do this, let $\Phi_i(t,s,x)$ be the solution to the ODE $\bar{F}_i$ at time $t$, starting from $x$ at the initial times $s$, i.e.,
\begin{equation}
\label{eq:Phi}
\frac{d}{dt} \Phi_i(t,s,x) = \bar{F}_i(t,\Phi_i(t,s,x)), \quad \Phi_i(s,s,x) = x.
\end{equation}
Such $\Phi_i(t,s,x)$ is well defined if $0 \leq s \leq t \leq T$, and $\Phi_i(t,s,x) \in X_i$ for all such $t$'s.
Recall that $\Phi(t,\cdot,\cdot)$ are themselves solutions of another homogeneous PDE, formally stated in the following theorem:

\begin{theorem}
Let $\Phi_i(t,s,x)$ be defined as in \eqref{eq:Phi}.
Given any $t \in [0,T]$, $\Phi_i(t,s,x)$ satisfies 
\begin{equation}
    \frac{d}{ds}\Phi_i(t,s,x) + \nabla_x \Phi_i(t,s,x) \cdot \bar{F}_i(s,x) = 0
\end{equation}
whenever $\Phi(t,s,x)$ is well-defined.
\end{theorem}
\begin{IEEEproof}
    The result follows directly by differentiating the semigroup identity
    \begin{equation}
        \Phi_i(t,s,\Phi_i(s,\tau, z)) = \Phi_i(t, \tau, z)
    \end{equation}
    with respect to $s$, and then performing change of variables $x = \Phi_i(s, \tau, z)$.
\end{IEEEproof}

This observation leads to another useful corollary:

\begin{corollary}
\label{cor:bdry}
Let $\sigma^i$ and $\eta^i$ satisfy the non-homogenous PDE \eqref{eq:pde} and let $\Phi_i$ satisfy Equation \eqref{eq:Phi}, then
\begin{equation}
    \int_{[0,T] \times X_i} w(\Phi_i(T,s,x)) \, d\left( \sigma^i(s,x) - \eta^i(s,x) \right) = 0
\end{equation}
for any test function $w \in C^1(X_i)$.
\end{corollary}

\begin{IEEEproof}
    Define a test function $v(s,x) := w(\Phi_i(T,s,x))$. By Equation \eqref{eq:pde_explained},
    \begin{align*}
        & \int_{[0,T] \times X_i} w(\Phi_i(T,s,x)) \, d\left( \sigma^i(s,x) - \eta^i(s,x) \right) \\
        =& \int_0^T \int_{X_i} \left( \partial_s v(s,x) + \nabla_x v(s,x) \cdot \bar{F}_i(s,x) \right) \, d\mu^i_{x \mid s}(x) \, ds\\
        =& \int_0^T \int_{X_i} \left( \frac{\partial w(\Phi_i(T,s,x))}{\partial \Phi_i(T,s,x)} \frac{d}{ds} \Phi_i(T,s,x) + \frac{\partial w(\Phi_i(T,s,x))}{\partial \Phi_i(T,s,x)} \nabla_x \Phi_i(T,s,x) \cdot \bar{F}_i(s,x)\right) \, d\mu^i_{x \mid s}(x) \, ds = 0
    \end{align*}
    where the last step is from chain rule.
\end{IEEEproof}

In the case when the vector field satisfies certain regularity requirements, we can begin to establish a converse to Lemma \ref{lem:traj2hle} by first showing that $\Phi$ is uniquely defined and subsequently showing that it satisfies an important relationship with $\mu^i_{x \mid t}$ almost everywhere. 

\begin{theorem}
\label{thm:pde_sol}
Suppose $\bar{F}_i \in L^1([0,T]; W^{1,\infty}(X_i; \R^{n_i}))$.
Given $\sigma^i, \eta^i \in \M_+([0,T] \times X_i)$, the solution to \eqref{eq:pde_explained} is given by
\begin{equation}
\label{eq:pde_sol}
\mu^i_{x \mid t} = \Phi_i(t,\cdot,\cdot)_\# \, \left(\sigma^i - \eta^i\right)
\end{equation}
for almost every $t \in [0,T]$, where $\Phi_i(t,\cdot,\cdot): [0,t] \times X_i \to X_i$ is defined according to \eqref{eq:Phi}.
\end{theorem}
\begin{IEEEproof}
We prove the result in two steps: first, we show the measure defined by formula \eqref{eq:pde_sol} is a solution of \eqref{eq:pde_explained}, then we show this solution is unique $dt$-almost everywhere.


We first verify \eqref{eq:pde_sol} satisfies Equation \eqref{eq:pde_explained}. 
Notice we need to check the distributional equality only on test functions of the form $\psi(t)w(x)$, i.e.,
\begin{equation}
\label{eq:thm11.1}
\int_0^T \dot{\psi}(t) \, \ip{\mu^i_{x \mid t}, w} \, dt = \int_{[0,T] \times X_i} \psi(t)w(x) \, d\left(\eta^i(t,x) - \sigma^i(t,x)\right) - \int_0^T \psi(t) \, \ip{\mu^i_{x \mid t}, \nabla_x w \cdot \bar{F}_i} \, dt
\end{equation}

We next substitute \eqref{eq:pde_sol} into the left-hand side of \eqref{eq:thm11.1} and show it is equal to the right-hand side of  \eqref{eq:thm11.1}:
\begin{align}
 \int_0^T \dot{\psi}(t) \int_{X_i} &w(x) \, d\mu^i_{x \mid t}(x) \, dt = \int_0^T \dot{\psi}(t) \int_0^t \int_{X_i} w(\Phi_i(t,s,x)) \, d\left(\sigma^i(s,x) - \eta^i(s,x)\right) \, dt \nonumber\\
    =& \int_0^T \int_{X_i} \left( \int_s^T \dot{\psi}(t) w(\Phi_i(t,s,x)) \, dt \right) \, d\left(\sigma^i(s,x) - \eta^i(s,x)\right) \label{eq:arg:1}\\
    \begin{split}
        =&\int_0^T \int_{X_i} \bigg( \psi(T) w(\Phi_i(T,s,x)) - \psi(s) w(\Phi_i(s,s,x) + \\
        & \qquad - \int_s^T \psi(t) \frac{d}{dt} w(\Phi_i(t,s,x))) \, dt \bigg) \, d\left(\sigma^i(s,x) - \eta^i(s,x)\right)
    \end{split} \label{eq:arg:2}\\
    \begin{split}
        =& 0 + \int_{[0,T] \times X_i} \psi(s)w(x) \, d\left(\eta^i(s,x) - \sigma^i(t,x)\right) + \\
        & \qquad - \int_0^T \int_{X_i} \left( \int_s^T \psi(t) \left( \nabla_x w(\Phi_i(t,s,x)) \cdot \bar{F}_i(t, \Phi_i(t,s,x)) \right) \, dt \right) \, d\left(\sigma^i(s,x) - \eta^i(s,x)\right) 
    \end{split} \label{eq:arg:3}\\
    \begin{split}
        =& \int_{[0,T] \times X_i} \psi(s)w(x) \, d\left(\eta^i(s,x) - \sigma^i(t,x)\right) + \\
        & \qquad - \int_0^T \psi(t) \left( \int_0^t \int_{X_i} \left( \nabla_x w(\Phi_i(t,s,x)) \cdot \bar{F}_i(t, \Phi_i(t,s,x)) \right)  d\left(\sigma^i(s,x) - \eta^i(s,x)\right) \right) \, dt
    \end{split} \label{eq:arg:4}\\
    =& \int_{[0,T] \times X_i} \psi(s)w(x) \, d\left(\eta^i(s,x) - \sigma^i(s,x)\right) - \int_0^T \psi(t) \ip{\mu^i_{x \mid t}, \nabla_x w \cdot \bar{F}_i} \, dt \label{eq:arg:5}
\end{align}
where \eqref{eq:arg:1} is deduced by Fubini's theorem; \eqref{eq:arg:2} follows from integration by parts; \eqref{eq:arg:3} is from Corollary \ref{cor:bdry} and \eqref{eq:Phi}; \eqref{eq:arg:4} is from Fubini's theorem; and \eqref{eq:arg:5} follows from \eqref{eq:pde_sol}.
Therefore \eqref{eq:pde_sol} is a solution to the distributional PDE \eqref{eq:pde_explained}.

    We next show the solution is unique $dt$-almost everywhere.
    Suppose there exists measures $\mu^{i}_{x \mid t, 1}$, $\mu^{i}_{x \mid t, 2}\in \M_+(X_i)$ defined for $t \in [0,T]$ that both satisfy Equation \eqref{eq:pde_explained}. Let $\mu^{i}_{x \mid t, 3} := \mu^{i}_{x \mid t, 1} - \mu^{i,}_{x \mid t, 2} \in \M(X_i)$, then $\mu^{i}_{x \mid t, 3}$ satisfies:
    
    \begin{equation}
    \label{eq:homo_pde}
        \int_0^T \int_{X_i} \left( \partial_t v(t,x) + \nabla_x v(t,x) \cdot \bar{F}_i \right) \, d\mu^{i}_{x \mid t, 3} \, dt = 0
    \end{equation}
    According to the proof of \cite[Lemma~3]{henrion2014convex}, such $\mu^{i}_{x \mid t, 3}$ is defined uniquely $dt$-almost everywhere.
    Notice the zero measure $\mu^{i}_{x \mid t, 3}(A) = 0$ for all sets A in the Borel $\sigma$-algebra of $X_i$ and all $t \in [0,T]$ is also a solution to \eqref{eq:homo_pde}, 
    therefore $\mu^{i}_{x \mid t, 3}$ is zero measure for almost every $t \in [0,T]$.
    As a result, $\mu^{i}_{x \mid t, 1}$ and $\mu^{i}_{x \mid t, 2}$ are equal $dt$-almost everywhere, and the solution to \eqref{eq:pde_explained} is unique $dt$-almost everywhere.
\end{IEEEproof}

In practice, the vector field $\bar{F}_i$ may not satisfy the regularity condition required to apply Theorem \ref{thm:pde_sol}; as a result, $\Phi_i$ may not be uniquely defined. 
To deal with solutions to a non-smooth ODE, we construct the notion of evaluation maps that act on the space of all absolutely continuous functions.
Let $AC(\R; \R^{n_i})$ be the space of absolutely continuous functions from $\R$ into $\R^{n_i}$ endowed with the sup norm. 
Define an \emph{evaluation map} $e_t: [0,t]\times [t,T] \times AC(\R;\R^{n_i}) \to \R^{n_i}$:
\begin{equation}
e_t(s, \tau, \gamma) = \gamma(t), \quad s \leq t \leq \tau
\end{equation}
for each $t \in [0,T]$.
As we prove next, this evaluation map establishes a relationship between admissible solutions to vector fields that may not satisfy the regularity conditions described in Theorem \ref{thm:pde_sol} and $\mu_{x \mid t}$:

\begin{theorem}\label{thm:monster}
Let $\mu^i_{x \mid t} \in \M_+(X_i)$ solve the PDE \eqref{eq:pde_explained} and assume that $\bar{F}_i(t,x)$ is pointwise bounded, i.e., $\exists M < +\infty$ such that $\forall (t,x) \in [0,T] \times X_i$, $|\bar{F}_i(t,x)| \leq M$.
Let $\Gamma_i := AC(\R;X_i)$.
Then there exists a measure $\rho^i \in \M_+([0,T] \times [0,T] \times \R^{n_i})$ such that
\begin{enumerate}[label=(\alph*)]
    \item \label{enum:thm12.1}
    $\rho^i$ is concentrated on the triplets $(s,\tau,\gamma)$, where $s \leq \tau$, and $\gamma \in \Gamma_i$ are solutions of the ODE $\dot{\gamma}(t) = \bar{F}_i(t,\gamma(t))$ for almost every $t \in [s,\tau]$. \label{prop:thm11.1}
    \Later{Need to show $s \leq \tau$..}
    \item \label{enum:thm12.2}
    For almost every $t \in [0,T]$, $\mu^i_{x \mid t}$ satisfies the following equality: \label{prop:thm11.2}
    \begin{equation}
    \label{eq:thmmonstereq}
        \mu^i_{x \mid t} = (e_t)_\# \, \rho^i
    \end{equation}
\end{enumerate}
\end{theorem}

\begin{IEEEproof}
This proof consists of several steps: in Step 1, we use a family of mollifiers parameterized by $\epsilon$ to smooth the vector field and all relevant measures and establish a relationship between the smooth measures using the solution to the smooth vector field via Theorem \ref{thm:pde_sol}; 
in Step 2, we prove that all trajectories that satisfy this smooth vector field and enter the domain, eventually leave the domain, and vice versa; 
in Steps 3 and 4, we prove a connection between the time at which each trajectory enters and leaves; 
since Steps 2-4 are all proven for the ``smoothed'' versions of the vector field and measures, in Step 5 we prove that there exists a limiting measure as the parameter controlling smoothness, $\epsilon$, goes to zero and;
in Step 6, we prove that this limit satisfies \eqref{eq:thmmonstereq};
in Step 7, we prove the first part of Theorem \ref{thm:monster} when the vector field is continuous;
in Step 8, we approximate the discontinuous vector field with a sequence of smooth functions and bound the approximation error;
in Step 9, we prove the first part of Theorem \ref{thm:monster} for arbitrary bounded vector fields.

\emph{Step 1 (Regularization).} We first mollify $\mu^i_{x \mid t}$ with respect to the space variable 
using a family of strictly positive mollifiers $\{\theta_\epsilon\} \subset C^\infty(\R^{n_i})$ with unit mass, zero mean, and bounded second moment,
obtaining smooth measures $\mu^i_{x \mid t; \epsilon}$ and smooth vector fields $\bar{F}_i^\epsilon$.

Define
\begin{equation}
\begin{split}
    \mu^i_{x \mid t; \epsilon} :=& \mu^i_{x \mid t} * \theta_\epsilon\\
    \bar{F}_i^\epsilon(t, \cdot) :=& 
        \begin{cases}
            \frac{\bar{F}_i(t,\cdot) \mu^i_{x \mid t} * \theta_\epsilon }{\mu^i_{x \mid t} * \theta_\epsilon}, & \text{if } \|\mu^i_{x \mid t}\| > 0\\
            0, & \text{if } \|\mu^i_{x \mid t}\| = 0
        \end{cases} \\
    \sigma^i_\epsilon :=& \sigma^i * \theta_\epsilon\\
    \eta^i_\epsilon :=& \eta^i * \theta_\epsilon
\end{split}
\end{equation}
Such smooth vector field $\bar{F}_i^\epsilon$ is pointwise bounded, because:
\begin{equation}
\label{eq:Fbar_bd}
    \left|\bar{F}_i^\epsilon(t,\cdot)\right| \leq \frac{\left|\bar{F}_i(t,\cdot)\right| \mu^i_{x \mid t} * \theta_\epsilon}{\mu^i_{x \mid t} * \theta_\epsilon} \leq M \frac{\mu^i_{x \mid t} * \theta_\epsilon}{\mu^i_{x \mid t} * \theta_\epsilon} = M
\end{equation}
The mollified measures are also bounded, because
\begin{equation}
\label{eq:thm12.0.0}
    \left\|\mu^i_{x \mid t; \epsilon} \right\| \leq \left\| \mu^i_{x \mid t} \right\| \left\| \int_{\R^{n_i}} \theta_\epsilon (x) \, dx \right\| = \mu^i_{x \mid t}(X_i)
\end{equation}
and similarly,
\begin{equation}
    \left\|\sigma^i_\epsilon \right\| \leq \sigma^i([0,T]\times X_i), \qquad 
    \left\|\eta^i_\epsilon \right\| \leq \eta^i([0,T] \times X_i)
\end{equation}

We next show that $\mu^i_{x \mid t; \epsilon}$, $\sigma^i_\epsilon$, $\eta^i_\epsilon$ is a solution to \eqref{eq:pde_explained} with respect to $\bar{F}_i^\epsilon$:

\Later{I really don't know how to prove this elegently...}
It is immediate that $\mu^i_{x \mid t; \epsilon}$ is a solution of \eqref{eq:pde_explained} with respect to $\bar{F}_i^\epsilon$, $\sigma^i_\epsilon$, and $\eta^i_\epsilon$. 
Also notice $\bar{F}_i^\epsilon \in L^1([0,T]; W^{1,\infty}(X_i; \R^{n_i}))$, therefore Theorem \cref{thm:pde_sol} can be applied to get
\begin{equation}
\label{eq:mu|t_everywhere}
    \mu^i_{x \mid t; \epsilon} = \Phi_i^\epsilon(t,\cdot,\cdot)_\# \left(\sigma^i_\epsilon - \eta^i_\epsilon\right)
\end{equation}
for almost every $t \in [0,T]$, where $\Phi_i^\epsilon(t,s,x)$ is the solution of
\begin{equation}
\label{eq:phi_eps}
\frac{d}{dt} \Phi_i^\epsilon(t,s,x) = \bar{F}_i^\epsilon(t,\Phi^\epsilon_i(t,s,x)), \quad \Phi_i^\epsilon(s,s,x) = x. 
\end{equation}
for almost every $t \in [s, T]$.
Such function $\Phi^\epsilon_i(\cdot,s,x)$ can be extended to a larger domain $\R$ without any difficulty due to the regularity of $\bar{F}^\epsilon_i$.
If we define a zero extension of $\bar{F}_i^\epsilon$ onto $\R \times \R^{n_i}$, denoted as $\hat{\bar{F}}_i^\epsilon$
(note such $\hat{\bar{F}}_i^\epsilon \in L^1([0,T]; W^{1, \infty}(\R^{n_i}))$), $\Phi^\epsilon_i(\cdot,s,x)$ can be further extended to the domain $\R$,
and we denote the extended version as $\hat{\Phi}^\epsilon_i(\cdot,s,x) \in \Gamma_i$ for any $(s,x) \in [0,T] \times \R^{n_i}$.
The space of all such functions is denoted as $\Gamma_i^\epsilon := \{\hat{\Phi}^\epsilon_i(\cdot,s,x)\}_{(s,x)}$ endowed with the sup norm,
therefore $\Gamma^\epsilon_i$ is a subset of $\Gamma_i$ with the subspace topology.
It follows by the existence and uniqueness theorem for ODE that the evaluation map $e_t(0,T,\cdot): \Gamma_i^\epsilon \to \R^{n_i}$ is an isomorphism for any $t \in [0,T]$.
We now define a map $\Psi^\epsilon: (t,x) \mapsto \hat{\Phi}^\epsilon(\cdot, t,x)$ from $[0,T]\times \R^{n_i}$ to $\Gamma_i^\epsilon$, and also a projection map $\pi^1: (s,x) \mapsto s$ from $[0,T] \times \R^{n_i}$ to $[0,T]$. Let
\begin{equation}
\label{eq:thm12:def_rho_eps}
\begin{split}
    \rho^{i,+}_{\epsilon} &:= \left( \pi^1 \times \Psi^\epsilon \right)_\# \sigma^i_{\epsilon} \in \M_+([0,T] \times \Gamma_i^\epsilon)\\
    \rho^{i,-}_{\epsilon} &:= \left( \pi^1 \times \Psi^\epsilon \right)_\# \eta^i_{\epsilon} \in \M_+([0,T] \times \Gamma_i^\epsilon)
\end{split}
\end{equation}

\emph{Step 2 (Marginals of $\rho^{i,+}_\epsilon$ and $\rho^{i,-}_\epsilon$).} 
This step shows that all trajectories that enter the domain via $\sigma_\epsilon^i$ leave through $\eta_\epsilon^i$ by proving that the $\gamma$-marginals of $\rho^{i,+}_\epsilon$ and $\rho^{i,-}_\epsilon$ are equal.
Since $\Gamma_i^\epsilon$ is isomorphic to $\R^{n_i}$ under the isomorphism $e_T(0,T,\cdot)$, $\R \times \Gamma_i^\epsilon$ is isomorphic to $\R \times \R^{n_i}$ which is Polish.
Therefore by definition $\R \times \Gamma_i^\epsilon$ is Souslin.
Using \cite[Corollary~10.4.13]{bogachev2007measure}, the measures $\rho^{i,+}_\epsilon$ and $\rho^{i,-}_\epsilon$ can be disintegrated as
\begin{equation}
\label{eq:thm12.3:disintegrate_rho}
\begin{split}
    d\rho^{i,+}_\epsilon(s,\gamma) &= d\rho^{i,+}_{s \mid \gamma; \epsilon}(s) \, d\rho^{i,+}_{\gamma; \epsilon}(\gamma)\\
    d\rho^{i,-}_\epsilon(\tau,\gamma) &= d\rho^{i,-}_{\tau \mid \gamma; \epsilon}(\tau) \, d\rho^{i,-}_{\gamma;\epsilon}(\gamma)
\end{split}
\end{equation}
where $\rho^{i,+}_{s \mid \gamma; \epsilon}$ and $\rho^{i,-}_{\tau \mid \gamma; \epsilon}$ are probability measures for all $\gamma \in \spt(\rho^{i,+}_{\gamma; \epsilon})$ and $\gamma \in \spt(\rho^{i,-}_{\gamma;\epsilon})$, respectively. We next show the $\gamma$-marginals are equal.

We first define zero extensions of $\mu^i_\epsilon$, $\sigma^i_\epsilon$, and $\eta^i_\epsilon$ to $[0,2T]\times \R^{n_i}$ as $\hat{\mu}^i_\epsilon$, $\hat{\sigma}^i_\epsilon$, and $\hat{\eta}^i_\epsilon$, respectively. It is immediate that they satisfy the PDE \eqref{eq:pde} on $[0,2T]\times \R^{n_i}$ with respect to $\hat{\bar{F}}_i^\epsilon$.
Using Theorem \cref{thm:pde_sol}, we know
\begin{equation}
    \label{eq:thm12.2.0}
    \hat{\mu}^i_{x \mid t; \epsilon}
    = \hat{\Phi}_i^\epsilon(t,\cdot,\cdot)_\# \left( \hat{\sigma}^i_\epsilon - \hat{\eta}^i_\epsilon \right)
    = \hat{\Phi}_i^\epsilon(t,\cdot,\cdot)_\# \left( \sigma^i_\epsilon - \eta^i_\epsilon \right)
\end{equation}
for almost every $t \in [0,2T]$,
where $\hat{\Phi}_i^\epsilon(t,\cdot,\cdot)$ is defined on $[0,t] \times \R^{n_i}$.
Moreover, since $\hat{\bar{F}}_i^\epsilon(t,x) = 0$ for all $(t,x) \in (T,2T] \times \R^{n_i}$, we know $\hat{\Phi}_i^\epsilon(t,s,x) = \hat{\Phi}_i^\epsilon(T,s,x)$ for all $(t,s,x) \in (T,2T]\times [0,T] \times \R^{n_i}$.
Now suppose there is a set $E_\Gamma \subset \Gamma_i^\epsilon$ such that $\rho^{i,+}_{\epsilon, \Gamma}(E_\Gamma) \neq \rho^{i,-}_{\epsilon, \Gamma}(E_\Gamma)$, 
and define $E := \{e_T(0,T,\gamma)\}_{\gamma \in E_\Gamma} \subset \R^{n_i}$ to be the image of $E_\Gamma$ under the isomorphism $e_T(0,T,\cdot)$. 
Notice
\begin{align}
    \hat{\mu}^i_\epsilon((T,2T] \times E)
    =& \int_{(T,2T]} \int_{\R^{n_i}} \mathbbm{1}_E (x) \, d\hat{\mu}^i_\epsilon(x \mid t) \, dt \label{eq:thm12.2.1}\\
    =& \int_{(T,2T]} \int_{[0,T] \times \R^{n_i}} \mathbbm{1}_E (\hat{\Phi}^\epsilon_i(t,s,x)) \, d \left(\sigma^i_\epsilon(s,x) - \eta^i_\epsilon(s,x) \right)\, dt \label{eq:thm12.2.2}\\
    =& \int_{(T,2T]} \int_{[0,T] \times \R^{n_i}} \mathbbm{1}_E (\hat{\Phi}^\epsilon_i(T,s,x)) \, d \left(\sigma^i_\epsilon(s,x) - \eta^i_\epsilon(s,x) \right)\, dt \label{eq:thm12.2.3}\\
    =& \int_{(T,2T]} \int_{[0,T] \times \R^{n_i}} \mathbbm{1}_E (e_T(0,T,\Psi^\epsilon(s,x))) \, d \left(\sigma^i_\epsilon(s,x) - \eta^i_\epsilon(s,x) \right)\, dt \label{eq:thm12.2.4}\\
    =& \int_{(T,2T]} \int_{[0,T] \times \Gamma_i^\epsilon} \mathbbm{1}_E (e_T(0,T,\gamma)) \, d \left( \rho^{i,+}_{\epsilon}(s,\gamma) - \rho^{i,-}_{\epsilon}(\tau,\gamma) \right) \, dt \label{eq:thm12.2.5}\\
    =& \int_{(T,2T]} \int_{\Gamma_i^\epsilon} \mathbbm{1}_{E_\Gamma} (\gamma) \, d \left( \rho^{i,+}_{\gamma;\epsilon}(\gamma) - \rho^{i,-}_{\gamma;\epsilon}(\gamma) \right) \, dt \label{eq:thm12.2.6}\\
    =& T \left( \rho^{i,+}_{\gamma;\epsilon}(E_\Gamma) - \rho^{i,-}_{\gamma;\epsilon}(E_\Gamma) \right) \label{eq:thm12.2.7}\\
    \neq & 0 \label{eq:thm12.2.8}
\end{align}
where \eqref{eq:thm12.2.1} follows by definition;
\eqref{eq:thm12.2.2} follows by plugging in \eqref{eq:thm12.2.0};
since the evolution of the system stays fixed for $t>T$, \eqref{eq:thm12.2.3} holds;
\eqref{eq:thm12.2.4} follows from the definition of $\Psi^\epsilon$;
\eqref{eq:thm12.2.5} follows from the definition of $\rho^{i,+}_\epsilon$ and $\rho^{i,-}_\epsilon$;
\eqref{eq:thm12.2.6} is because $\gamma \in E_\Gamma$ if and only if $e_T(0,T,\gamma) \in E$, and because the conditional measures $\rho^{i,+}_{s \mid \gamma; \epsilon}$ and $\rho^{i,-}_{\tau, \mid \gamma; \epsilon}$ are probability measures;
\eqref{eq:thm12.2.7} is true by definition.
The result \eqref{eq:thm12.2.8} contradicts $\hat{\mu}^i_\epsilon$ being a zero extension, therefore $\rho^{i,+}_{\gamma;\epsilon} = \rho^{i,-}_{\gamma;\epsilon}$.

\emph{Step 3 (Construct $\rho^i_{\epsilon,\delta}$).} 
We now want to combine $\rho^{i,+}_\epsilon$ and $\rho^{i,-}_\epsilon$ to generate a measure $\rho^i_\epsilon \in \M_+([0,T]\times [0,T] \times \Gamma^\epsilon_i)$ that describes the trajectories that evolve in the domain as well as their entering and exiting time. 
Such a measure can be defined by pushing forward $\rho^{i,+}_\epsilon$ through a map that associates entering time to exiting time. 
However, such a map may not be well defined; for example, two trajectories can enter the domain at the same time but leave at different times. 
To address such issues, we mollify the $t$-component and define a sequence of measures $\rho^i_{\epsilon,\delta}$ first, and then define $\rho^i_\epsilon$ as the limit of this sequence as $\delta \downarrow 0$ which is done in Step 4.

Let $\{\theta_\delta\} \subset C^\infty(\R)$ be a family of smooth mollifiers with first moment equal to 1, and define
\begin{equation}
\begin{split}
    \rho^{i,+}_{s \mid \gamma; \epsilon,\delta} &:= \rho^{i,+}_{s \mid \gamma; \epsilon} * \theta_\delta\\
    \rho^{i,-}_{\tau \mid \gamma; \epsilon,\delta} &:= \rho^{i,-}_{\tau \mid \gamma; \epsilon} * \theta_\delta
\end{split}
\end{equation}
We further define measures $\rho^{i,+}_{\epsilon,\delta}, \rho^{i,-}_{\epsilon,\delta} \in \M_+(\R \times \Gamma_i^\epsilon)$ as
\begin{equation}
\begin{split}
     d\rho^{i,+}_{\epsilon, \delta}(s,\gamma) & = d\rho^{i,+}_{\epsilon, \delta}(s \mid \gamma) \, d\rho^{i,+}_{\epsilon,\Gamma}(\gamma)\\
    d\rho^{i,-}_{\epsilon, \delta}(\tau, \gamma) & = d\rho^{i,-}_{\epsilon, \delta}(\tau \mid \gamma) \, d\rho^{i,-}_{\epsilon,\Gamma}(\gamma)
\end{split}
\end{equation}


For almost every $t \in [0,T]$ and any non-negative test function $w \in L^1(\R^{n_i})$, we know
\begin{align}
    0 \leq & \ip{\mu^i_{x \mid t; \epsilon}, w} \label{eq:thm12.3.21}\\
    = & \int_{[0,t] \times \R^{n_i}} w(\Phi^\epsilon_i(t,s,x)) \, d\left( \sigma^i_\epsilon (s,x) - \eta^i_\epsilon(s,x) \right) \label{eq:thm12.3.22}\\
    = & \int_{[0,t] \times \Gamma_i^\epsilon} w(e_t(0,T,\gamma)) \, d\rho^{i,+}_\epsilon(s,\gamma) - \int_{[0,t] \times \Gamma_i^\epsilon} w(e_t(0,T,\gamma)) \, d\rho^{i,-}_\epsilon(\tau,\gamma) \label{eq:thm12.3.23}\\
    = & \int_{\Gamma_i^\epsilon} w(e_t(0,T,\gamma)) \left( \int_{[0,t]} d\rho^{i,+}_{s \mid \gamma; \epsilon}(s) - \int_{[0,t]} d\rho^{i,-}_{\tau \mid \gamma; \epsilon}(\tau) \right) \, d\rho^i_{\gamma; \epsilon}(\gamma) \label{eq:thm12.3.24}
\end{align}
where
\eqref{eq:thm12.3.21} follows from the fact that $\mu^i_{x \mid t; \epsilon}$ is an unsigned measure;
\eqref{eq:thm12.3.22} follows by substituting in \eqref{eq:mu|t_everywhere};
\eqref{eq:thm12.3.23} follows from the definition \eqref{eq:thm12:def_rho_eps};
\eqref{eq:thm12.3.24} follows from \eqref{eq:thm12.3:disintegrate_rho}.

Equivalently, given any Borel set $E_\Gamma \subset \Gamma^\epsilon_i$,
\begin{equation}
\label{eq:rho_monotone}
\int_{E_\Gamma} \left( \int_{[0,t]} d\rho^{i,+}_{s \mid \gamma; \epsilon}(s) - \int_{[0,t]} d\rho^{i,-}_{\tau \mid \gamma; \epsilon}(\tau) \right) \, d\rho^i_{\gamma;\epsilon}(\gamma) \geq 0
\end{equation}
Since the functions $t \mapsto \int_{[0,t]} d\rho^{i,+}_{s \mid \gamma; \epsilon}(s)$ and $t \mapsto \int_{[0,t]} d\rho^{i,-}_{\tau \mid \gamma; \epsilon}(\tau)$ are absolutely continuous,
Equation \eqref{eq:rho_monotone} is satisfied for all $t \in [0,T]$:
\begin{equation}
    \int_{[0,t]} d\rho^{i,+}_{s \mid \gamma; \epsilon}(s) \geq \int_{[0,t]} d\rho^{i,-}_{\tau \mid \gamma; \epsilon}(\tau)
\end{equation}
for almost every $\gamma \in \spt(\rho^i_{\epsilon,\Gamma})$. 
From the monotonicity of convolution \Later{Prove this using Fubini, change of variables and def of convolution!!}, we know
\begin{equation}
    \int_{-\infty}^t d\rho^{i,+}_{s \mid \gamma; \epsilon,\delta}(s) \geq \int_{-\infty}^t d\rho^{i,-}_{\tau \mid \gamma; \epsilon,\delta}(\tau)
\end{equation}
for almost every $\gamma \in \spt(\rho^i_{\gamma;\epsilon})$.


Because $\rho^{i,+}_{s \mid \gamma; \epsilon,\delta}$ and $\rho^{i,-}_{\tau \mid \gamma; \epsilon,\delta}$ are smooth non-negative measures, the functions $t \mapsto \int_{-\infty}^t \, d\rho^{i,+}_{s \mid \gamma; \epsilon,\delta}(s)$ and $t \mapsto \int_{-\infty}^t \, d\rho^{i,-}_{\tau \mid \gamma; \epsilon,\delta}(\tau)$ are continuous and non-decreasing. Moreover, since
\begin{equation}
0 =　\int_{-\infty}^{-\infty} \, d\rho^{i,+}_{s \mid \gamma; \epsilon,\delta}(s) \leq \int_{-\infty}^t \, d\rho^{i,-}_{\tau \mid \gamma; \epsilon,\delta}(\tau) \leq \int_{-\infty}^t \, d\rho^{i,+}_{s \mid \gamma; \epsilon,\delta}(s) \leq \int_{-\infty}^{\infty} \, d\rho^{i,-}_{\tau \mid \gamma; \epsilon,\delta}(\tau) = 1,
\end{equation}
where the last equality follows because $\rho^{i,-}_{s \mid \gamma ;\epsilon, \delta}$ is a probability measure;
by the Mean Value Theorem, for any $\gamma \in \spt(\rho^i_{\gamma; \epsilon})$ there exists a function $r_\gamma: \R \to \R$ such that
\begin{equation}
r_\gamma(t) \geq t
\end{equation}
and
\begin{equation}
\label{eq:mvt}
\int_{-\infty}^t \, d\rho^{i,+}_{s \mid \gamma; \epsilon,\delta}(s) = \int_{-\infty}^{r_\gamma(t)} \, d\rho^{i,-}_{\tau \mid \gamma; \epsilon,\delta}(\tau)
\end{equation}
for every $\gamma \in \spt(\rho^i_{\gamma; \epsilon})$. Moreover, the function $r_\gamma$ is strictly increasing and therefore invertible, i.e., there exists a function $r_\gamma^{-1}:\R \to \R$ such that $r_\gamma(r_\gamma^{-1}(t)) = r_\gamma^{-1}(r_\gamma(t)) = t$.

Because of the result in Step 2,
Equation \eqref{eq:mvt} can be equivalently written as
\begin{equation}
\label{eq:thm12.3:+vs-}
\int_{\R \times \Gamma^\epsilon_i} \mathbbm{1}_{(-\infty,t]}(s) w(\gamma) \, d\rho^{i,+}_{\epsilon,\delta}(s,\gamma) = \int_{\R \times \Gamma^\epsilon_i} \mathbbm{1}_{(-\infty,r_\gamma(t)]}(\tau) w(\gamma) \, d\rho^{i,-}_{\epsilon,\delta}(\tau,\gamma)
\end{equation}
for any $t \in \R$ and where $w$ is any measurable function on $\Gamma_i$.

We now abuse notation and define a map $r : \R \times \spt(\rho^i_{\gamma; \epsilon}) \to \R$ by letting $r(s,\gamma) := r_\gamma(s)$ for all $\gamma \in \spt(\rho^i_{\gamma; \epsilon})$,
and also projection maps $\pi^1 : (s,\gamma)\in \R \times \Gamma_i^\epsilon \mapsto s \in \R$,
$\pi^2 : (s,\gamma) \in \R \times \Gamma_i^\epsilon \mapsto \gamma \in \Gamma_i^\epsilon$.
We can then define a measure $\rho^i_{\epsilon, \delta} \in \M_+(\R \times \R \times \Gamma_i^\epsilon)$ as
\begin{equation}
\label{eq:thm12.3.0}
\rho^i_{\epsilon,\delta} = \left( \pi^1 \times r \times \pi^2 \right)_\# \rho^{i,+}_{\epsilon, \delta}
\end{equation}

We now establish the relationship between the marginals of $\rho^i_{\epsilon,\delta}$ and the measures $\rho^{i,+}_{\epsilon,\delta}$ and $\rho^{i,-}_{\epsilon,\delta}$.
We use variables $(s,\tau,\gamma) \in \R \times \R \times \Gamma_i^\epsilon$ to denote any point in $\spt(\rho^i_{\epsilon,\delta})$.
Then the $(s,\gamma)$-marginal of $\rho^i_{\epsilon, \delta}$ is equal to $\rho^{i,+}_{\epsilon, \delta}$, because:
\begin{align}
    \rho^i_{\epsilon, \delta}(A \times [0,T] \times B) 
    =& \int_{\R \times \R \times \Gamma^\epsilon_i} \mathbbm{1}_{A\times B}(s, \gamma ) \, d\rho^i_{\epsilon,\delta}(s,\tau, \gamma) \label{eq:thm12.3.1}\\
    =& \int_{\R \times \Gamma_i^\epsilon} \mathbbm{1}_{A \times B}( \pi^1(s,\gamma), \pi^2(s,\gamma) ) \, d\rho^{i,+}_{\epsilon, \delta}(s,\gamma) \label{eq:thm12.3.2}\\
    =& \int_{\R \times \Gamma_i^\epsilon} \mathbbm{1}_{A \times B}( s, \gamma ) \, d\rho^{i,+}_{\epsilon, \delta}(s,\gamma) \label{eq:thm12.3.3}\\
    =& \rho^{i,+}_{\epsilon,\delta}(A \times B) \label{eq:thm12.3.4}
\end{align}
for all Borel sets $A \in [0,T]$ and $B \in \Gamma^\epsilon_i$,
where \eqref{eq:thm12.3.1} follows by definition;
\eqref{eq:thm12.3.2} follows from substituting in \eqref{eq:thm12.3.0};
\eqref{eq:thm12.3.3} follows from substituting the definition of $\pi^1$ and $\pi^2$;
\eqref{eq:thm12.3.4} is true by definition.

To show the $(\tau,\gamma)$-marginal of $\rho^i_{\epsilon,\delta}$ is equal to $\rho^{i,-}_{\epsilon,\delta}$, it is then sufficient to show the following equation
\begin{equation}
\int_{\R \times \R \times \Gamma^\epsilon_i} \mathbbm{1}_{(-\infty,t]}(\tau) w(\gamma) \, d\rho^i_{\epsilon, \delta}(s,\tau,\gamma) =
\int_{\R \times \Gamma^\epsilon_i} \mathbbm{1}_{(-\infty,t]}(\tau) w(\gamma) \, d\rho^{i,-}_{\epsilon, \delta}(\tau, \gamma)
\end{equation}
holds for all indicator functions $\mathbbm{1}_{(\infty,t]}$ and $w \in L^1(\Gamma_i^\epsilon)$. The equation is true because
\begin{align}
    \int_{\R \times \R \times \Gamma^\epsilon_i} \mathbbm{1}_{(-\infty,t]}(\tau) w(\gamma) \, d\rho^i_{\epsilon, \delta}(s,\tau,\gamma) 
    =& \int_{\R \times \Gamma^\epsilon_i} \mathbbm{1}_{(-\infty,t]}(r(s,\gamma)) w(\gamma) \, d\rho^{i,+}_{\epsilon, \delta}(s,\gamma) \label{eq:thm12.3.5}\\
    =& \int_{\R \times \Gamma^\epsilon_i} \mathbbm{1}_{(-\infty,r_\gamma^{-1}(t)]}(s) w(\gamma) \, d\rho^{i,+}_{\epsilon, \delta}(s,\gamma) \label{eq:thm12.3.6}\\
    =& \int_{\R \times \Gamma^\epsilon_i} \mathbbm{1}_{(-\infty, r_\gamma( r_\gamma^{-1}(t) )]}(\tau) w(\gamma) \, d\rho^{i,-}_{\epsilon, \delta}(\tau, \gamma) \label{eq:thm12.3.7}\\
    =& \int_{\R \times \Gamma^\epsilon_i} \mathbbm{1}_{(-\infty, t]}(\tau) w(\gamma) \, d\rho^{i,-}_{\epsilon, \delta}(\tau, \gamma) \label{eq:thm12.3.8}
\end{align}
where
\eqref{eq:thm12.3.5} follows by plugging in \eqref{eq:thm12.3.0};
\eqref{eq:thm12.3.6} is true because $r_\gamma$ is strictly monotonic and therefore $r_\gamma(s) \in (-\infty,t]$ if and only if $s \in (-\infty,r_\gamma^{-1}(t)]$;
\eqref{eq:thm12.3.7} follows by substituting in \eqref{eq:thm12.3:+vs-};
\eqref{eq:thm12.3.8} follows from the fact that $r_\gamma$ is invertible;

\emph{Step 4 (Tightness of the family $\{\rho^i_{\epsilon,\delta}\}_{\delta}$).}
We now show that the limit of $\rho_{\epsilon,\delta}^{i}$ exists as $\delta$ goes to zero and that this limiting measure satisfies $\mu^i_{x \mid t; \epsilon} = (e_t)_\# \rho^i_\epsilon$ for almost every $t$. 
To prove this limiting condition, we use the notion of tightness of measures and apply the following pair of conditions called the \emph{integral condition on the tightness} and \emph{tightness criterion} \cite[pp. 605-606]{maniglia2007probabilistic}:

\noindent {\bf Integral Condition on the Tightness:} Let $X$ be a separable metric space. A family $\mathcal{K} \subset \M_+(X)$ is tight if and only if there exists a function $\Theta: X \to [0,+\infty]$ whose sublevel sets $\{x \in X \mid \Theta(x) \leq c\}$ are compact in $X$ (such functions are called \emph{coercive functions}), such that
\begin{equation}
    \sup_{\mu \in \mathcal{K}} \int_X \Theta(x) \, d\mu(x) < +\infty
\end{equation}

\noindent {\bf Tightness Criterion:} Let $X$, $X_1$, $X_2$ be separable metric spaces and let $r^i: X \to X_i$, $i = 1, 2$ be continuous maps such that the product map $r: r^1 \times r^2: X \to X_1 \times X_2$ is proper. 
Let $\mathcal{K} \subset \M_+(X)$ be such that $\mathcal{K}_i:=r^i_\#(\mathcal{K})$ is tight in $\M_+(X_i)$ for $i = 1, 2$. 
Then also $\mathcal{K}$ is tight in $\M_+(X)$. 
Notice the statement also holds for finitely many maps by induction.

We now choose maps $r^1$, $r^2$ defined in $\R \times \R \times \Gamma^\epsilon_i$ as
\begin{equation}
r^1: (s,\tau,\gamma) \mapsto (s,\gamma) \in \R \times \Gamma^\epsilon_i, \qquad 
r^2: (s,\tau,\gamma) \mapsto \tau \in \R
\end{equation}
It is obvious \Later{explain} that $r = r^1 \times r^2: \R \times \R \times \Gamma^\epsilon_i \to \R \times \Gamma^\epsilon_i \times \R$ is proper. 
The family $\{r^1_\# \rho^i_{\epsilon,\delta}\}_\delta$ is given by $\{\rho^{i,+}_{\epsilon,\delta}\}_\delta$ which are tight by definition \Later{explain},
and the family $\{r^2_\# \rho^i_{\epsilon,\delta}\}_\delta$ is given by the first marginal of $\{\rho^{i,-}_{\epsilon,\delta}\}_\delta$ which are also tight. \Later{explain.}
Applying the tightness criterion, we know the family $\{\rho^i_{\epsilon,\delta}\}_\delta$ is tight, 
and therefore narrowly sequentially relatively compact as $\delta \downarrow 0$ according to Prokhorov compactness theorem.
Let $\rho^i_{\epsilon}$ be any limit of the family $\{\rho^i_{\epsilon,\delta}\}$ as $\delta \downarrow 0$.
Since the $(s,\gamma)$-marginal of $\rho^i_{\epsilon, \delta}$ is equal to $\rho^{i,+}_{\epsilon,\delta}$
and the $(\tau,\gamma)$-marginal of $\rho^i_{\epsilon, \delta}$ is equal to $\rho^{i,-}_{\epsilon,\delta}$,
for arbitrary continuous function $\varphi \in C_b(\R \times \Gamma_i)$ we have
\begin{equation}
\label{eq:tightness_marginals}
\begin{split}
    \int_{\R \times \R \times \Gamma_i} \varphi(s,\gamma) \,            d\rho^i_{\epsilon,\delta}(s,\tau,\gamma)
    =& \int_{\R  \times \Gamma_i} \varphi(s,\gamma) \, d\rho^{i,+}_{\epsilon,\delta}(s,\gamma) \\
    \int_{\R \times \R \times \Gamma_i} \varphi(\tau,\gamma) \,            d\rho^i_{\epsilon,\delta}(s,\tau,\gamma)
    =& \int_{\R  \times \Gamma_i} \varphi(\tau,\gamma) \, d\rho^{i,-}_{\epsilon,\delta}(\tau,\gamma) \\
\end{split}
\end{equation}
We then pass to the limit $\delta \downarrow 0$ in \eqref{eq:tightness_marginals} to obtain
\begin{equation}
\label{eq:tightness_marginals_limit0}
\begin{split}
    \int_{\R \times \R \times \Gamma_i} \varphi(s,\gamma) \,            d\rho^i_{\epsilon}(s,\tau,\gamma)
    =& \int_{\R  \times \Gamma_i} \varphi(s,\gamma) \, d\rho^{i,+}_{\epsilon}(s,\gamma) \\
    \int_{\R \times \R \times \Gamma_i} \varphi(\tau,\gamma) \,            d\rho^i_{\epsilon}(s,\tau,\gamma)
    =& \int_{\R  \times \Gamma_i} \varphi(\tau,\gamma) \, d\rho^{i,-}_{\epsilon}(\tau,\gamma) \\
\end{split}
\end{equation}
Since $\spt(\rho^{i,+}_{\epsilon}) \subset [0,T] \times \Gamma_i^\epsilon$ and $\spt(\rho^{i,-}_{\epsilon}) \subset [0,T] \times \Gamma_i^\epsilon$, we know $\spt(\rho^i_\epsilon) \subset [0,T] \times [0,T] \times \Gamma_i^\epsilon$, and \eqref{eq:tightness_marginals_limit0} can be written as
\begin{equation}
\label{eq:tightness_marginals_limit}
\begin{split}
    \int_{[0,T] \times [0,T] \times \Gamma_i^\epsilon} \varphi(s,\gamma) \,            d\rho^i_{\epsilon}(s,\tau,\gamma)
    =& \int_{[0,T]  \times \Gamma_i^\epsilon} \varphi(s,\gamma) \, d\rho^{i,+}_{\epsilon}(s,\gamma) \\
    \int_{[0,T] \times [0,T] \times \Gamma_i^\epsilon} \varphi(\tau,\gamma) \,            d\rho^i_{\epsilon}(s,\tau,\gamma)
    =& \int_{[0,T]  \times \Gamma_i^\epsilon} \varphi(\tau,\gamma) \, d\rho^{i,-}_{\epsilon}(\tau,\gamma) \\
\end{split}
\end{equation}

In fact, \eqref{eq:tightness_marginals_limit} is also true for arbitrary measurable function $\varphi$ \cite[Theorem 7.14.25]{bogachev2007measure}.
Notice for any triplet $(s,\tau, \gamma) \in \spt(\rho^i_\epsilon)$, we have $\gamma \in \Gamma^\epsilon_i$, therefore $\gamma$ satisfies the ODE:
\begin{equation}
\label{eq:gamma_ode}
    \dot{\gamma} = \bar{F}_i^\epsilon(t,\gamma(t))
\end{equation}
for all $t \in [0,T]$.

For almost every $t \in [0,T]$, and for any $t \in [0,T]$ and measurable test function $w: X_i \to \R$, we have
\begin{align}
    & \int_{X_i} w(x) \, d\mu^i_\epsilon (x \mid t) \nonumber \\
    =& \int_{[0,t] \times X_i} w(\Phi^\epsilon_i(t,s,x)) \, d(\sigma^i_\epsilon(s,x) - \eta^i_\epsilon(s,x)) \label{eq:thm12.3.11}\\
    =& \int_{[0,t] \times X_i} w(e_t(s, T, \Psi^\epsilon(s,x))) \, d\sigma^i_\epsilon(s,x) - \int_{[0,t] \times X_i} w(e_t(\tau,T,\Psi^\epsilon(\tau,x))) \, d\eta^i_\epsilon(\tau,x) \label{eq:thm12.3.12}\\
    =& \int_{[0,t] \times [0,T] \times \Gamma_i^\epsilon} w(e_t(s,T,\gamma)) \, d\rho^i_\epsilon(s,\tau,\gamma) - \int_{[0,t] \times [0,t] \times \Gamma_i^\epsilon} w(e_t(\tau,T,\gamma)) \, d\rho^i_\epsilon(s,\tau,\gamma) \label{eq:thm12.3.13}\\
    =& \int_{[0,t]\times[0,t]\times \Gamma_i^\epsilon} \left( w(e_t(s,T,\gamma)) - w(e_t(\tau,T,\gamma)) \right) \, d\rho^i_\epsilon(s,\tau,\gamma) + \int_{[0,t] \times (t,T] \times \Gamma_i^\epsilon} w(e_t(s,T,\gamma))\, d\rho^i_\epsilon(s,\tau,\gamma) \label{eq:thm12.3.14}\\
    =& \int_{[0,t]\times[0,t]\times \Gamma_i^\epsilon} \left( w(e_t(0,T,\gamma)) - w(e_t(0,T,\gamma)) \right) \, d\rho^i_\epsilon(s,\tau,\gamma) + \int_{[0,t]\times (t,T] \times \Gamma_i^\epsilon} w(e_t(s,\tau,\gamma)) \, d\rho^i_\epsilon(s,\tau,\gamma) \label{eq:thm12.3.15}\\
    =& 0 + \int_{[0,t]\times [t,T] \times \Gamma_i^\epsilon} w(e_t(s,\tau,\gamma)) \, d\rho^i_\epsilon(s,\tau,\gamma) - \int_{[0,t]\times \{t\} \times \Gamma_i^\epsilon} w(e_t(0,T,\gamma)) \, d\rho^i_\epsilon(s,\tau,\gamma) \label{eq:thm12.3.16}
\end{align}
where
\eqref{eq:thm12.3.11} follows by substituting in \eqref{eq:mu|t_everywhere};
\eqref{eq:thm12.3.12} follows from the fact that $\Phi^\epsilon_i(t,s,x)$ is defined only when $t \geq s$ and a change of variable $\tau = s$ for the second half of the equation;
\eqref{eq:thm12.3.13} follows by substiting in \eqref{eq:thm12:def_rho_eps} and \eqref{eq:tightness_marginals_limit};
\eqref{eq:thm12.3.14} follows by splitting the domain of integration;
\eqref{eq:thm12.3.15} follows from the fact that $e_t(t_1,T,\cdot) = e_t(0,T,\cdot)$ and $e_t(t_1,T,\cdot) = e_t(t_1,t_2,\cdot)$ for all $0 \leq t_1 \leq t \leq t_2 \leq T$;
the first term of \eqref{eq:thm12.3.16} is zero because the integrand is a zero function, and the rest
follows by adding and subtracting the set $[0,t] \times \{t\} \times \Gamma_i^\epsilon$ to the domain of integration.

Since $\int_{[0,t]\times \{t\} \times \Gamma_i^\epsilon} w(e_t(0,T,\gamma)) \, d\rho^i_\epsilon(s,\tau,\gamma)$ is non-zero for at most countably many $t$'s (otherwise $\rho^i_\epsilon$ would not be bounded), we can write the previous equation as
\begin{equation}
\label{eq:thm12:mu_eps2rho_eps}
\mu^i_{x \mid t; \epsilon} = (e_t)_\# \rho^i_\epsilon
\end{equation}
for almost every $t \in [0,T]$.

\emph{Step 5 (Tightness of the family $\{\rho^i_\epsilon\}_\epsilon$).}
We now show that the limit of $\rho_{\epsilon}^{i}$ exists as $\epsilon$ goes to zero.
This is done by applying the tightness criterion again. 
To begin, choose maps $r^1$, $r^2$, $r^3$ defined in $[0,T] \times [0,T] \times \Gamma_i$ as
\begin{equation}
r^1:(s,\tau,\gamma) \mapsto s \in [0,T], \qquad
r^2:(s,\tau,\gamma) \mapsto \tau \in [0,T], \qquad
r^3:(s,\tau,\gamma) \mapsto \gamma \in \Gamma_i
\end{equation}
Observe that $r = r^1 \times r^2 \times r^3: \R \times \R \times \Gamma_i \to \R \times \R \times \Gamma_i$ is identity map and therefore proper.
The family $\{r^1_\# \rho^i_\epsilon\}_\epsilon$ and $\{r^2_\# \rho^i_\epsilon\}_\epsilon$ are given by the first marginals of $\sigma^i_\epsilon$ and $\eta^i_\epsilon$, respectively, which are tight (in fact they are independent of $\epsilon$).
About $r^3_\# \rho^i_\epsilon$ we choose a coercive function $\Theta: \Gamma_i \to \R$:
\begin{equation}
    \Theta: \gamma \mapsto |\gamma(T)| + \int_0^T |\dot{\gamma}(t)|^2 \, dt
\end{equation}
where $|\cdot|$ is the Euclidean norm on $\R^{n_i}$.
\Later{explain or citation to why this is coercive? $\Theta$ may not be coercive if $\Gamma_i$ is defined on $\R$. Will come back to this point later.}
\begin{align}
\int_{\Gamma_i} \Theta (\gamma) \, d \left(r^3_\# \rho^i_\epsilon \right)(\gamma)
    =& \int_{[0,T] \times [0,T] \times \Gamma^i_\epsilon} \Theta( \gamma ) \, d\rho^i_\epsilon( s, \tau, \gamma ) \label{eq:thm12.5.1}\\
    =& \int_{[0,T] \times \Gamma^i_\epsilon} \Theta( \gamma ) \, d\rho^{i,+}_\epsilon( s, \gamma ) \label{eq:thm12.5.2} \\
    =& \int_{[0,T] \times \R^{n_i}} \left( \left|\hat{\Phi}^\epsilon_i(T,s,x) \right| + \int_0^T \left|\dot{\hat{\Phi}}^i_\epsilon(t,s,x) \right|^2 \, dt \right) \, d\sigma^i_\epsilon(s,x) \label{eq:thm12.5.3}\\
    \begin{split}
        \leq& \int_{[0,T] \times \R^{n_i}} \left( |\hat{\Phi}^\epsilon_i(s,s,x)| + \int_s^T \left| \bar{F}_i^\epsilon(\hat{\Phi}^i_\epsilon (\tau,s,x)) \right| \, d\tau \right. + \\
        & \qquad + \left. \int_0^T \left| \bar{F}^\epsilon_i (\hat{\Phi}^i_\epsilon(t,s,x)) \right|^2 \, dt \right) \, d\sigma^i_\epsilon(s,x)
    \end{split} \label{eq:thm12.5.4}\\
    \leq& \int_{[0,T] \times \R^{n_i}} |x| \, d\sigma^i_\epsilon(s,x) + MT \sigma^i_\epsilon([0,T] \times \R^{n_i}) + M^2 T \sigma^i_\epsilon([0,T] \times \R^{n_i})  \label{eq:thm12.5.5}\\
    \leq & \int_{[0,T] \times \R^{n_i}} (|x|^2 + 1) \, d\sigma^i_\epsilon(s,x) + (MT + M^2 T) \sigma^i_\epsilon([0,T] \times \R^{n_i}) \label{eq:thm12.5.6}\\
    \leq & \int_{[0,T] \times X_i} \int_{\R^{n_i}} |x+y|^2 \theta_\epsilon(y) \, dy \, d\sigma^i(s,x) + (1 + MT + M^2T) \sigma^i([0,T] \times X_i) \label{eq:thm12.5.7}\\
    \begin{split}
        =& \int_{[0,T] \times X_i} |x|^2 \, d\sigma^i(s,x) 
        + \left( \int_{\R^{n_i}} |y|^2 \theta_\epsilon(y) \, dy \right) \cdot \sigma([0,T] \times X_i) + \\
         & \qquad + \int_{[0,T] \times X_i} \int_{\R^{n_i}} 2x^Ty \cdot \theta_\epsilon(y) \, dy \, d\sigma^i(s,x) 
         + (1 + MT + M^2T) \sigma^i([0,T] \times X_i)
    \end{split} \label{eq:thm12.5.8}\\
    < & +\infty \label{eq:thm12.5.9}
\end{align}
where
\eqref{eq:thm12.5.1} follows from definition of pushforward measure;
\eqref{eq:thm12.5.2} follows from \eqref{eq:tightness_marginals_limit};
\eqref{eq:thm12.5.3} follows from the definition \eqref{eq:thm12:def_rho_eps};
\eqref{eq:thm12.5.4} follows from the \eqref{eq:Phi} and the triangle inequality in $\R^{n_i}$;
\eqref{eq:thm12.5.5} follows from \eqref{eq:Phi} and \eqref{eq:Fbar_bd};
\eqref{eq:thm12.5.6} is because $|x|^2 + 1 \geq |x|$ for all $x \in \R^{n_i}$, and $\sigma^i_\epsilon$ is unsigned measure;
\eqref{eq:thm12.5.7} follows by definition of convolution;
\eqref{eq:thm12.5.8} is because $|x+y|^2 = |x|^2 + |y|^2 + 2x^Ty$ for all $x,y \in \R^{n_i}$;
Since $\sigma^i$ is bounded by assumption and $X_i$ is compact therefore $|x|^2$ is bounded for all $x \in X_i$, the first and last term in \eqref{eq:thm12.5.8} are bounded. Because $\theta_\epsilon$ is assumed to have zero mean and bounded second moment, the second term in \eqref{eq:thm12.5.8} is bounded and the third term in \eqref{eq:thm12.5.8} is zero. Then \eqref{eq:thm12.5.9} follows.
As a result, using the integral condition for the tightness,
$\{r^3_\# \rho^i_\epsilon\}_\epsilon$ is tight, and therefore the family $\{\rho^i_\epsilon\}_\epsilon$ is tight according to the tightness criterion.

\emph{Step 6 (Condition \ref{enum:thm12.2}).} We now prove that the limiting measure of $\rho^i_\epsilon$ as $\epsilon$ goes to zero satisfies \eqref{eq:thmmonstereq}.
Using the Prokhorov Compactness Theorem, the family $\rho^i_\epsilon$ is narrowly sequentially relatively compact. 
We choose a narrowly convergent sequence in $\{\rho^i_\epsilon\}_\epsilon$ and define its limit by $\rho^i \in \M_+([0,T]\times [0,T] \times \Gamma_i)$.
Given any function $w \in C_b(\R^{n_i})$, it follows from \eqref{eq:thm12:mu_eps2rho_eps} that
\begin{equation}
\label{eq:thm12.6.1}
    \int_{\R^{n_i}} w(x) \, d\mu^i_{x \mid t; \epsilon}(x) = \int_{[0,T] \times [0,T] \times \Gamma_i} w(e_t(s,\tau,\gamma)) \, d\rho^i_\epsilon(s,\tau,\gamma)
\end{equation}
for almost every $t \in [0,T]$.
Since $e_t$ is continuous for all $t \in [0,T]$, we know $w \circ e_t \in C_b( [0,T] \times [0,T] \times \Gamma_i)$. 
We then pass to the limit $\epsilon \downarrow 0$ to both sides of \eqref{eq:thm12.6.1} to obtain
\begin{equation}
\label{eq:thm12.6.2}
    \int_{\R^{n_i}} w(x) \, d\mu^i_{x \mid t}(x) = \int_{[0,T] \times [0,T] \times \Gamma_i} w(e_t(s,\tau,\gamma)) \, d\rho^i(s,\tau,\gamma)
\end{equation}
for almost every $t \in [0,T]$.
In fact, \eqref{eq:thm12.6.2} is true for arbitrary measurable function $w \in L^1(\R^{n_i})$ because $C_b(\R^{n_i})$ is dense in $L^1(\R^{n_i})$ \cite[Corollary~4.2.2]{bogachev2007measure}, therefore we have
\begin{equation}
\label{eq:thm12.6.res}
    \mu^i_{x \mid t} = (e_t)_\# \rho^i
\end{equation}
for almost every $t \in [0,T]$.

\emph{Step 7 (Condition \ref{enum:thm12.1} with continuous vector field).}
To prove \ref{enum:thm12.1}, it suffices to show that
\begin{equation}
\label{eq:rho_traj}
    \int_{[0,t] \times [t,T] \times \Gamma_i} \left| \gamma(t) - \gamma(s) - \int_s^t \bar{F}_i(\tau', \gamma(\tau')) \, d\tau' \right| \, d\rho^i(s,\tau,\gamma) = 0
\end{equation}
for any $t \in [0,T]$.
One could attempt to prove this result, by applying the notion of narrow convergence to $\rho_\epsilon^i$; however, the technical difficulty is that the integrand may not be continuous due to the lack of regularity of $\bar{F}_i$.
We may, however, first consider the case where the vector field $\bar{F}_i$ is continuous on $[0,T] \times X_i$ (therefore uniformly continuous since $[0,T] \times X_i$ is compact).
In that instance, we could apply narrow convergence directly to prove our result; however, we instead prove the result using a slightly different technique since the result is useful in subsequent steps.

Let $v \in C_b([0,T] \times X_i; \R^{n_i})$ be a bounded uniformly continuous function,
and let us first prove the estimate:
\begin{equation}
\label{eq:F_estimate}
\int_{[0,t] \times [t,T] \times \Gamma_i} \left| \gamma(t) - \gamma(s) - \int_s^t v(\tau', \gamma(\tau')) \, d\tau' \right| \, d\rho^i(s, \tau, \gamma) \leq \int_{[0,T] \times X_i} \left| \bar{F}_i(\tau, x) - v(\tau, x) \right| \, d\mu^i_{\tau,x}(\tau,x)
\end{equation}
for any $t \in [0,T]$.
Indeed we have
\begin{align}
    & \int_{[0,t] \times [t,T] \times \Gamma_i} \left| \gamma(t) - \gamma(s) - \int_s^t v(\tau',\gamma(\tau')) \, d\tau' \right| \, d \rho^i_\epsilon(s,\tau,\gamma) \nonumber\\
    =& \int_{[0,t] \times [t,T] \times \Gamma_i} \left| \int_s^t \bar{F}^\epsilon_i( \tau', \gamma(\tau') ) \, d\tau' - \int_s^t v(\tau',\gamma(\tau')) \, d\tau' \right| \, d\rho^i_\epsilon(s,\tau,\gamma) \label{eq:thm12.7.1} \\
    =& \int_{[0,t] \times [t,T] \times \Gamma_i} \int_s^t \left| \bar{F}^\epsilon_i( \tau', \gamma(\tau') ) - v(\tau', \gamma(\tau')) \right| \, d\tau' \, d\rho^i_\epsilon(s,\tau,\gamma) \label{eq:thm12.7.2} \\
    =& \int_0^t \int_{[0,\tau'] \times [t,T] \times \Gamma_i} \left| \bar{F}^\epsilon_i( \tau', \gamma(\tau') ) - v(\tau', \gamma(\tau')) \right| \, d\rho^i_\epsilon(s,\tau,\gamma) \, d\tau' \label{eq:thm12.7.3}\\
    =& \int_0^t \int_{[0,\tau'] \times [t,T] \times \Gamma_i} \left| \bar{F}^\epsilon_i( \tau', e_{\tau'}(s,\tau,\gamma) ) - v(\tau', e_{\tau'}(s,\tau,\gamma)) \right| \, d\rho^i_\epsilon(s,\tau,\gamma) \, d\tau' \label{eq:thm12.7.4}\\
    \leq & \int_0^t \int_{[0,\tau'] \times [\tau',T] \times \Gamma_i} \left| \bar{F}^\epsilon_i( \tau', e_{\tau'}(s,\tau,\gamma) ) - v(\tau', e_{\tau'}(s,\tau,\gamma)) \right| \, d\rho^i_\epsilon(s,\tau,\gamma) \, d\tau' \label{eq:thm12.7.5}\\
    =& \int_0^t \int_{\R^{n_i}} \left| \bar{F}^\epsilon_i( \tau', x ) - v(\tau', x) \right| \, d\mu^i_{x \mid \tau' ; \epsilon}(x) \, d\tau' \label{eq:thm12.7.6}\\
    \leq & \int_0^t \int_{\R^{n_i}} \left| \bar{F}_i^\epsilon(\tau,x) - v^\epsilon(\tau,x) \right| \, d\mu^i_{x \mid \tau ; \epsilon} (x) \, d\tau
        + \int_0^t \int_{\R^{n_i}} \left| v^\epsilon(\tau,x) - v(\tau,x) \right| \, d\mu^i_{x \mid \tau; \epsilon}(x) \, d\tau \label{eq:thm12.7.7}\\
    \leq & \int_0^t \int_{\R^{n_i}} \left| \bar{F}_i(\tau,x) - v(\tau,x) \right| \, d\mu^i_{x \mid \tau} (x) \, d\tau
        + \int_0^T \left( \sup_{x \in \R^{n_i}} |v^\epsilon(\tau,x) - v(\tau,x)| \cdot \mu^i_{x \mid \tau; \epsilon}(\R^{n_i}) \right) \, d\tau \label{eq:thm12.7.8}\\
    \leq & \int_0^T \int_{X_i} \left| \bar{F}_i(\tau,x) - v(\tau,x) \right| \, d\mu^i_{\tau,x}(\tau, x)
        + \sup_{\substack{\tau \in [0,T] \\ x \in \R^{n_i}}} |v^\epsilon(\tau,x) - v(\tau,x)| \cdot \mu^i([0,T] \times X_i) \label{eq:thm12.7.9}
\end{align}
where
\eqref{eq:thm12.7.1} follows from \eqref{eq:gamma_ode} and the Fundamental Theorem of Calculus;
\eqref{eq:thm12.7.2} follows from properties of integrals;
\eqref{eq:thm12.7.3} follows from Fubini's Theorem;
\eqref{eq:thm12.7.4} follows from definition of evaluation map $e_t$;
\eqref{eq:thm12.7.5} is because $\rho^i_\epsilon$ is an unsigned measure and $\tau' \leq t$ in the domain of integration;
\eqref{eq:thm12.7.6} follows from \eqref{eq:thm12:mu_eps2rho_eps};
\eqref{eq:thm12.7.7} follows after a change of variables $\tau = \tau'$, adding and subtracting $v^\epsilon(\tau, \cdot) := \frac{(v(\tau, \cdot) \, \mu^i_{x \mid \tau}) * \theta_\epsilon}{ \mu^i_{x \mid \tau; \epsilon}}$ and applying the Triangle Inequality;
\eqref{eq:thm12.7.8} is due to \cite[Lemma~3.9]{maniglia2007probabilistic};
\eqref{eq:thm12.7.9} follows from \eqref{eq:thm12.0.0} and Lemma \ref{lem:disintegration}.
Since the family $\{\rho^i_\epsilon\}_\epsilon$ is tight and the integrand is a continuous and non-negative function on $[0,t] \times [t,T] \times \Gamma_i$, we can take the limit as $\epsilon \downarrow 0$ on both sides of the chain of inequalities. 
Since $v$ is uniformly continuous, $v^\epsilon$ converges to $v$ uniformly as $\epsilon \downarrow 0$, and the second term of \eqref{eq:thm12.7.9} converges to 0, therefore we obtain \eqref{eq:F_estimate}.
Notice if $\bar{F}_i$ is continuous (therefore uniformly continuous on the compact domain $[0,T]\times X_i$), we may let $v := \bar{F}_i$ in \eqref{eq:F_estimate}, and \eqref{eq:rho_traj} follows.

\emph{Step 8 (Error bound of vector field approximation).}
When there is no regularity in $\bar{F}_i$ other than boundedness, we choose a sequence of continuous functions converging to $\bar{F}_i$ in $L^ 1(\mu^i_{t,x}; \R^{n_i})$, and prove an error bound of the approximation:
Let $\{v_k\}_{k \in \N} \subset C([0,T] \times X_i; \R^{n_i})$ be a sequence of continuous functions converging to $\bar{F}_i$ in $L^1(\mu^i_{t,x}; \R^{n_i})$
\cite[Corollary~4.2.2]{bogachev2007measure}.
Given any $t \in [0,T]$, the error between $v_k$ and $\bar{F}_i$ is given by
\begin{align}
    & \int_{[0,t] \times [t,T] \times \Gamma_i} \int_s^t \left| v_k(\tau', \gamma(\tau')) - \bar{F}_i(\tau',\gamma(\tau')) \right| \, d\tau' \, d\rho^i(s,\tau,\gamma) \nonumber\\
    =& \int_0^t \int_{[0,\tau'] \times [t,T] \times \Gamma_i} \left| v_k(\tau', \gamma(\tau')) - \bar{F}_i(\tau',\gamma(\tau')) \right| \, d\rho^i(s,\tau,\gamma) \, d\tau' \label{eq:thm12.8.1}\\
    \leq & \int_0^t \int_{[0,\tau'] \times [\tau',T] \times \Gamma_i} \left| v_k(\tau', \gamma(\tau')) - \bar{F}_i(\tau',\gamma(\tau')) \right| \, d\rho^i(s,\tau,\gamma) \, d\tau' \label{eq:thm12.8.2}\\
    =& \int_0^t \int_{X_i} \left| v_k(\tau',x) - \bar{F}_i(\tau',x) \right| \, d\mu^i_{x \mid \tau'}(x) \, d\tau' \label{eq:thm12.8.3}\\
    \leq & \int_{[0,T] \times X_i} \left| v_k(\tau,x) - \bar{F}_i(\tau,x) \right| \, d\mu^i_{\tau,x}(\tau,x) \label{eq:thm12.8.4}
\end{align}
where
\eqref{eq:thm12.8.1} follows from Fubini's theorem;
\eqref{eq:thm12.8.2} is because $\rho^i$ is unsigned measure and $\tau' \leq t$ in the domain of integration;
\eqref{eq:thm12.8.3} follows from \eqref{eq:thm12.6.res};
in \eqref{eq:thm12.8.4} we performed change of variables $\tau = \tau'$ and the result follows from Lemma \ref{lem:disintegration}.
Observe that as $k \to \infty$ this error goes to zero.

\emph{Step 9 (Condition \ref{enum:thm12.1} with bounded vector field).}
We may now combine Step 7 and Step 8 together and prove condition \ref{enum:thm12.1} in a more general setting.
Using the results in Step 7 and Step 8, we obtain for any $t \in [0,T]$,
\begin{align}
    & \int_{[0,t] \times [t,T] \times \Gamma_i} \left| \gamma(t) - \gamma(s) - \int_s^t \bar{F}_i(\tau', \gamma(\tau')) \, d\tau' \right| \, d\rho^i(s,\tau,\gamma) \nonumber\\
    \begin{split}
        \leq & \int_{[0,t] \times [t,T] \times \Gamma_i} \left| \gamma(t) - \gamma(s) - \int_s^t v_k(\tau', \gamma(\tau')) \, d\tau' \right| \, d\rho^i(s,\tau,\gamma) +\\
        & \qquad + \int_{[0,t] \times [t,T] \times \Gamma_i} \left| \int_s^t v_k(\tau', \gamma(\tau')) \, d\tau' - \int_s^t \bar{F}_i(\tau', \gamma(\tau')) \, d\tau' \right| \, d\rho^i(s,\tau,\gamma)
    \end{split} \label{eq:thm12.9.1}\\
    \begin{split}
        \leq & \int_{[0,T] \times X_i} \left| \bar{F}_i(\tau,x) - v_k(\tau,x) \right| \, d\mu^i_{\tau,x}(\tau,x) +\\
        & \qquad + \int_{[0,t] \times [t,T] \times \Gamma_i} \int_s^t \left| v_k(\tau', \gamma(\tau'))- \bar{F}_i(\tau', \gamma(\tau')) \right| \, d\tau' \, d\rho^i(s,\tau,\gamma)
    \end{split} \label{eq:thm12.9.2}\\
    \leq & 2 \int_{[0,T] \times X_i} \left| \bar{F}_i(\tau,x) - v_k(\tau,x) \right| \, d\mu^i_{\tau,x}(\tau,x) \label{eq:thm12.9.3}
\end{align}
where
\eqref{eq:thm12.9.1} follows from triangle inequality;
\eqref{eq:thm12.9.2} follows from \eqref{eq:F_estimate};
\eqref{eq:thm12.9.3} follows from the error bound proved in Step 8.
When we let $k \to \infty$, \eqref{eq:thm12.9.3} goes to zero, therefore condition \ref{enum:thm12.1} holds.

\end{IEEEproof}

\begin{corollary}
\label{cor:rho_endpoints}
Assume $\bar{F}_i$ is pointwise bounded on $[0,T] \times X_i$.
Let $\mu^i_{x \mid t}$, $\sigma^i$, and $\eta^i$ satisfy the PDE \eqref{eq:pde_explained}, and let $\rho^i$ be defined as in Theorem \cref{thm:monster}.
Define maps $r^1, r^2 \in [0,T] \times [0,T] \times \Gamma_i \to [0,T] \times \R^{n_i}$ by
\begin{equation}
    \begin{split}
        r^1:& (s,\tau,\gamma) \mapsto (s, \gamma(s))\\
        r^2:& (s,\tau,\gamma) \mapsto (\tau,\gamma(\tau))
    \end{split}
\end{equation}
Then
\begin{equation}
\begin{split}
    r^1_\# \rho^i &= \sigma^i\\
    r^2_\# \rho^i &= \eta^i
\end{split}
\end{equation}

\Later{Theorem \cref{thm:monster} establishes a connection between the measure $\mu^i_{x \mid t}$ that solves the PDE \eqref{eq:pde_explained} and trajectories that satisfy the dynamics in mode $i$. We next show those trajectories start in the support of ???}

\begin{IEEEproof}
Recall in the proof of Theorem \cref{thm:monster} we mollified $\sigma^i$ and $\eta^i$ using a family of smooth mollifiers to obtain smooth measures $\sigma^i_\epsilon$ and $\eta^i_\epsilon$. We also defined a \emph{tight} family of measures $\{\rho^i_\epsilon\}_\epsilon \subset \M_+([0,T] \times [0,T] \times \Gamma_i)$ that converges to $\rho^i$ in the narrow sense. The connection between each $\rho^i_\epsilon$ in that family and the mollified measures $\sigma^i_\epsilon$ and $\eta^i_\epsilon$ was established via measures $\rho^{i,+}_\epsilon$ and $\rho^{i,-}_\epsilon$.

For all Borel subsets $A \times B \in [0,T] \times \R^{n_i}$, we have
\begin{align}
    \int_{[0,T] \times \R^{n_i}} \mathbbm{1}_{A \times B}(s,x) \, d(r^1_\# \rho^i_\epsilon)
    =& \int_{[0,T] \times [0,T] \times \Gamma_i} \mathbbm{1}_{A \times B}( s, \gamma(s) ) \, d\rho^i_\epsilon(s,\tau,\gamma) \label{eq:cor12.1}\\
    =& \int_{[0,T] \times \Gamma_i} \mathbbm{1}_{A \times B}( s, \gamma(s) ) \, d\rho^{i,+}_\epsilon(s,\gamma) \label{eq:cor12.2}\\
    =& \int_{[0,T] \times \R^{n_i}} \mathbbm{1}_{A \times B} ( s, \Phi_i^\epsilon(s,s,x) ) \, d\sigma^i_\epsilon(s,x) \label{eq:cor12.3}\\
    =& \int_{[0,T] \times \R^{n_i}} \mathbbm{1}_{A \times B} (s,x) \, \sigma^i_\epsilon(s,x) \label{eq:cor12.4}
\end{align}
where
\eqref{eq:cor12.1} follows by definition of pushforward measure;
\eqref{eq:cor12.2} follows from \eqref{eq:tightness_marginals_limit};
\eqref{eq:cor12.3} follows from \eqref{eq:thm12:def_rho_eps};
\eqref{eq:cor12.4} follows from \eqref{eq:phi_eps}.
Therefore for any continuous function $\varphi \in C_b( [0,T] \times \R^{n_i} )$, we know
\begin{equation}
    \int_{[0,T] \times [0,T] \times \Gamma_i} (\varphi \circ r^1)(s,\tau,\gamma) \, d\rho^i_\epsilon(s,\tau,\gamma) = \int_{[0,T] \times X_{n_i}} \varphi(s,x) \, d\sigma^i_\epsilon( s, x )
\end{equation}
Since the families $\{\sigma^i_\epsilon\}_\epsilon$ and $\{\rho^i_\epsilon\}_\epsilon$ are tight, as was shown in the proof of Theorem \cref{thm:monster},
and therefore narrowly sequentially relatively compact according to Prokhorov Compactness Theorem,
and because $r^1$ is continuous,
we can take the limit (in the narrow sense) as $\epsilon \downarrow 0$ to obtain 
\begin{equation}
\label{eq:cor12.5}
    \int_{[0,T] \times [0,T] \times \Gamma_i} (\varphi \circ r^1)(s,\tau,\gamma) \, d\rho^i(s,\tau,\gamma) = \int_{[0,T] \times X_{n_i}} \varphi(s,x) \, d\sigma^i( s, x )
\end{equation}
In fact \eqref{eq:cor12.5} is true for all measurable functions $\varphi : [0,T] \times \R^{n_i} \to \R$ because $C_b$ is dense \cite[Corollary~4.2.2]{bogachev2007measure}, as a result $r^1_\# \rho^i = \sigma^i$. 
The result for $\eta^i$ can be proved in a similar manner.
\end{IEEEproof}
\end{corollary}

As a consequence, any triplet $(s,\tau,\gamma) \in \spt(\rho^i)$ can be viewed as a trajectory $\gamma$ in mode $i$, well defined on $[s,\tau]$, and satisfying $(s,\gamma(s)) \in \spt(\sigma^i)$, $(\tau,\gamma(\tau)) \in \spt(\eta^i)$.
In fact, such trajectories in different modes are closely related by reset maps, and can be combined together to be admissible trajectories for the entire hybrid system.
To further illustrate this point, we first define an evaluation map that acts on the trajectories for the hybrid system $e^i_t: {\cal X} \to \coprod_{i \in \mathcal I} X_i$ as
\begin{equation}
    e^i_t(\gamma) =
    \begin{cases}
        \gamma_i(t) , & \text{if } \lambda(\gamma(t)) = i\\
        \emptyset, & \text{otherwise}
    \end{cases}
\end{equation}
for each $i \in \mathcal{I}$.
We can then establish a relationship between admissible trajectories and measures that satisfy \eqref{eq:pde_explained}:
\begin{theorem}
\label{thm:hle2traj}
Assume $\bar{F}_i$ is pointwise bounded on $[0,T] \times X_i$,
and let $\mu^i_{x \mid t}$, $\sigma^i$, and $\eta^i$ satisfy the PDE \eqref{eq:pde_explained}.
Then there exists a non-negative measure $\rho \in \M_+(\mathcal{X}_T)$ supported on a family of admissible trajectories, such that
\begin{enumerate}
    \item $\rho$ satisfies
        \begin{equation}
            \mu^i_{x \mid t} = (e_t^i)_\# \rho
        \end{equation}
        for almost every $t \in [0,T]$.
    \item If $\mu_0^i$ satisfies
        \begin{equation}
        \label{eq:thm13.0.1}
            \sum_{i \in {\cal I}} \mu_0^i(X_i) = 1
        \end{equation}
        then $\rho$ is a probability measure.
    \item $\mu^i_{t,x}$ (resp. $\mu^i_T$, $\mu^{S_e}$) is the average occupation measure (resp. average terminal measure, average guard measure) generated by the family of admissible trajectories in the support of $\rho$ for each mode $i \in \mathcal{I}$ and $e \in {\cal E}$.
\end{enumerate}
\end{theorem}
\begin{IEEEproof}
We first show that trajectories defined in support of $\rho^i$ and $\rho^j$ satisfy the reset map for all $(i,j) \in \mathcal{E}$.

\emph{Step 1 (Reset maps are satisfied).}
According to Corollary \ref{cor:rho_endpoints}, it suffices to show
\begin{equation}
\label{eq:rho_reset}
    \sigma^{j} = \delta_0 \otimes \mu_0^i + \sum_{(i,j) \in \mathcal{E}} \tilde{R}_{(i,j)\#} \eta^{i} \quad \forall j\in \mathcal{I}.
\end{equation}
Notice
\begin{equation}
    \eta^i = \delta_T \otimes \mu^i_T + \sum_{(i,i') \in \mathcal{E}} \mu^{S_{(i,i')}}
\end{equation}
where $\spt(\delta_T \otimes \mu^i_T) \cap \left( [0,T] \times S_{(i,j)} \right) = \emptyset$, and $\spt(\mu^{S_{(i,i')}}) \cap \left( [0,T] \times S_{(i,j)} \right) = \emptyset$ for all $i' \neq j$.
Therefore 
\begin{equation}
    \tilde{R}_{(i,j)\#} \eta^i = \tilde{R}_{(i,j)\#} \mu^{S_{(i,j)}}
\end{equation}
and \eqref{eq:rho_reset} follows from \eqref{eq:sigma&eta}.

As a result of Step 1, all trajectories in the support of $\rho^i$ are reinitialized to another trajectory in the support of $\rho^j$ after it reaches the guard $S_{(i,j)}$;
On the other hand, a trajectory can only start in mode $i$ either from the given initial condition $x_0$ at time 0, or by transitioning from another mode $j$ if $(j,i) \in \mathcal{E}$.
To be admissible, we must show that all trajectories start from $x_0$ at time 0, and reach $X_T$ at time $T$.

\emph{Step 2 (Trajectories are defined on $[0,T]$).}
This step shows that trajectories are defined between $[0,T]$. 
To prove this, we first show that for any $i \in \mathcal{I}$ and $(s,\tau,\gamma) \in \spt(\rho^i)$ such that $\tau \neq T$, there is a number $\Delta t > 0$ such that $\tau - s \geq \Delta t$.

Let $(s,\tau,\gamma) \in \spt(\rho^i)$ for some $i \in {\cal I}$, and let $0 \leq s \leq \tau < T$.
According to Theorem \cref{thm:monster} and Corollary \cref{cor:rho_endpoints}, we know 
\begin{equation}
\label{eq:thm13.2.0}
    \gamma(s) \in \{x_0\} \bigcup_{(i',i) \in \mathcal{E}} R_{(i',i)}(S_{(i',i)}), \qquad
    \gamma(\tau) \in \bigcup_{(i,i') \in \mathcal{E}} S_{(i,i')}, \qquad
    \dot{\gamma}(t) = \bar{F}_i(t,\gamma(t)).
\end{equation}
According to Definition \cref{def:hybrid_system}, Assumption \cref{assum:g_intersect}, and Assumption \cref{assum:x0},
$\{x_0\} \bigcup_{(i',i) \in \mathcal{E}} R_{(i',i)}(S_{(i',i)})$ and $\bigcup_{(i,i') \in \mathcal{E}} S_{(i,i')}$ are disjoint compact sets, therefore the distance between those two sets, denoted as $d_i$, is nonzero.
Let $M_i > 0$ be a bound for $\bar{F}_i(t,x)$ over $[0,T] \times X_i$, and define
\begin{equation}
    \Delta t:= \min_{i \in {\cal I}} \frac{d_i}{M_i}
\end{equation}
Therefore for all $(x,\tau,\gamma) \in \spt(\rho^i)$, we have
\begin{align}
    d_i \leq & | \gamma(\tau) - \gamma(s)| \label{eq:thm13.2.1} \\
    =& \left| \int_s^\tau \bar{F}_i( t, \gamma(t) ) \, dt \right| \label{eq:thm13.2.2} \\
    \leq & \int_s^\tau M_i \, dt \label{eq:thm13.2.3} \\
    =& (\tau - s) M_i \label{eq:thm13.2.4}
\end{align}
where 
\eqref{eq:thm13.2.1} and \eqref{eq:thm13.2.2} follows from \eqref{eq:thm13.2.0};
\eqref{eq:thm13.2.3} follows from definition of $M_i$;
\eqref{eq:thm13.2.4} follows from fundamental theorem of calculus.
As a result, $\tau - s \geq \frac{d_i}{M_i} \geq \Delta t$.

To show all trajectories are defined on $[0,T]$, consider the following case:
Let $(s,\tau,\gamma) \in \spt(\rho^i)$.
If $s = 0$ and $\tau = T$, we are done.
Now suppose $\tau < T$. As a result of Step 1 and Corollary \ref{cor:rho_endpoints}, 
$(\tau,\gamma(\tau)) \in S_{(i,i')}$ for some $i' \in {\cal I}$,
and $\gamma$ is reinitialized to another trajectory $\gamma'$ defined on $[\tau,\tau']$, such that $(\tau,\tau',\gamma') \in \rho^{i'}$. 
if $\tau' < T$, the trajectory is defined on $[s,\tau'] \subset [s, \tau + \Delta t]$. 
This process can always be continued until the trajectory is defined up to time $T$.
Similarly, we can show a trajectory defined on $[s,\tau]$ is also defined on $[s-\Delta t, \tau]$ as long as $s > 0$.
Therefore all trajectories are well defined on $[0,T]$.

Notice it follows from the above discussion that for any $i \in {\cal I}$ and $(0, \tau, \gamma) \in \spt(\rho^i)$, $\tau \geq \Delta t$. As a result, $\spt(\mu^{S_e}) \in [\Delta t, T] \times S_e$ for all $e \in {\cal E}$. Then according to the result in Step 1, we know
\begin{equation}
\label{eq:thm15.2.11}
    \spt( \tilde{R}_{e\#} \mu^{S_e} ) \in [\Delta t, T] \times R_e(S_e)
\end{equation}
for all $e \in {\cal E}$.

\emph{Step 3 (Trajectories are admissible).}
We first show all the trajectories end at time $t=T$ in $\spt(\mu_T^i)$ for some $i \in {\cal I}$.
Consider the triplet $(s,T, \gamma) \in \spt(\rho^i)$ for some $i \in {\cal I}$.
According to Corollary \cref{cor:rho_endpoints}, $(T,\gamma(T)) \in \spt(\eta^i) \subset \left( \{T\} \times X_{T_i} \right) \cup \left( [0,T] \times \bigcup_{(i,i') \in {\cal E}}S_{(i,i')} \right)$.
If $\gamma(T) \in X_{T_i}$, we are done.
If $\gamma(T) \in S_{(i,i')}$ for some $(i,i') \in {\cal E}$, then $\gamma$ is reinitialized to another trajectory $\gamma'$ defined on $[T,T]$, such that $(T,T,\gamma') \in \rho^{i'}$.
Applying Corollary \cref{cor:rho_endpoints} again, we know $\gamma'(T) \in X_{T_i} \bigcup_{(i,i') \in {\cal E}}S_{(i,i')}$.
Also, because $\gamma'(T) = R_{(i,i')}(\gamma(T)) \in R_{(i,i')}(S_{(i,i')})$, using Assumption \cref{assum:target_set} we know $\gamma'(T) \in X_{T_i}$.
Using a similar argument, we can also show all the trajectories start from $\spt(\mu^i_0)$, therefore they are admissible by definition.

\emph{Step 4 (Condition 1 and 2).}
As a result of Step 3, there exists a measure $\rho \in \M_+(\mathcal{X}_T)$ such that 
\begin{equation}
\label{eq:rho}
    (e_t^i)_\# \rho(\gamma) = (e_t)_\# \rho^i = \mu^i_{x \mid t}
\end{equation}
for almost every $t \in [0,T]$.

Since all trajectories are defined on $[0,T]$,
to prove $\rho$ is a probability measure, we only need to show
\begin{equation}
    \sum_{i \in {\cal I}} \int_{\{0\} \times [0,T] \times \Gamma_i} \, d\rho^i(s,\tau,\gamma) = 1
\end{equation}

Notice
\begin{align}
    \sum_{i \in {\cal I}} \int_{\{0\} \times [0,T] \times \Gamma_i} \, d\rho^i(s,\tau,\gamma)
    =& \sum_{i \in {\cal I}} \int_{[0,T]\times[0,T] \times \Gamma_i} \mathbbm{1}_{\{0\}}(s) \, d\rho^i(s,\tau,\gamma) \label{eq:thm15.4.1} \\
    =& \sum_{i \in {\cal I}} \int_{[0,T] \times X_i} \mathbbm{1}_{\{0\}}(s) \, d\sigma^i(s,x) \label{eq:thm15.4.2} \\
    =& \sum_{i \in {\cal I}} \sigma^i(\{0\} \times X_i) \label{eq:thm15.4.3} \\
    =& \sum_{i \in {\cal I}} \mu^i_0(X_i) \label{eq:thm15.4.4} \\
    =& 1 \label{eq:thm15.4.5}
\end{align}
where
\eqref{eq:thm15.4.1} and \eqref{eq:thm15.4.3} follows from definition of identity map;
\eqref{eq:thm15.4.2} follows from Corollary \cref{cor:rho_endpoints};
\eqref{eq:thm15.4.4} follows from \eqref{eq:sigma&eta} and \eqref{eq:thm15.2.11}.

\emph{Step 5 (average occupation measure, average terminal measure, and average guard measure).}

Let $A \times B$ be in the Borel $\sigma$-algebra of $[0,T] \times X_i$, then we have
\begin{align}
    \mu^i_{t,x}(A \times B) =& \int_{[0,T] \times X_i} \mathbbm{1}_{A \times B}(t,x) \, d\mu^i_{x \mid t}(x) \, dt \label{eq:thm15.5.1}\\
        =& \int_0^T \int_{\mathcal{X}_T} \mathbbm{1}_A(t) \cdot \mathbbm{1}_B(\gamma_i(t)) \, d\rho(\gamma) \, dt \label{eq:thm15.5.2} \\
        =& \int_{\mathcal{X}_T} \int_0^T \mathbbm{1}_{A \times B}(t, \gamma_i(t)) \, dt \, d\rho(\gamma) \label{eq:thm15.5.3}
\end{align}
where
\eqref{eq:thm15.5.1} follows from Lemma \cref{lem:disintegration};
\eqref{eq:thm15.5.2} follows from \eqref{eq:rho};
\eqref{eq:thm15.5.3} follows from Fubini's theorem.

Also, for all $B$ in the Borel $\sigma$-algebra of $X_{T_i}$, we know
\begin{align}
    \mu_T^i(B) =& (\delta_T \otimes \mu_T^i)(\{T\} \times B) \label{eq:thm15.5.11} \\
        =& \int_{[0,T] \times X_i} \mathbbm{1}_{\{T\} \times B}(t,x) \, d\eta^i(t,x) \label{eq:thm15.5.12} \\
        =& \int_{[0,T] \times [0,T] \times \Gamma_i} \mathbbm{1}_{\{T\}}(\tau) \cdot I_B(\gamma(\tau)) \, d\rho^i(s,\tau,\gamma) \label{eq:thm15.5.13} \\
        =& \int_{[0,T] \times [0,T] \times \Gamma_i} \mathbbm{1}_B(e_T(\gamma)) \, d\rho^i(s,\tau,\gamma) \label{eq:thm15.5.14} \\
        =& \int_{\mathcal{X}_T} \mathbbm{1}_B(\gamma_i(T)) \, d\rho(\gamma) \label{eq:thm15.5.15} 
\end{align}
where
\eqref{eq:thm15.5.11} follows from definition of $\delta_T$;
\eqref{eq:thm15.5.12} follows from Assumption \cref{assum:target_set}, \eqref{eq:sigma&eta}, and the fact that $B \subset X_{T_i}$;
\eqref{eq:thm15.5.13} follows from Corollary \cref{cor:rho_endpoints};
\eqref{eq:thm15.5.14} follows from definition of $e_t$;
\eqref{eq:thm15.5.15} follows from \eqref{eq:rho}.

Finally, for all $(i,i') \in {\cal S}$ and $A \times B$ in the Borel $\sigma$-algebra of $[0,T] \times S_{(i,i')}$, we have
\begin{align}
    \mu^{S_{(i,i')}}(A \times B) =& \int_{[0,T] \times X_i} \mathbbm{1}_{A \times B}(\tau, x) \, d\eta^i(\tau,x) \label{eq:thm15.5.21} \\
    =& \int_{[0,T] \times [0,T] \times \Gamma_i} \mathbbm{1}_{A \times B} (\tau, \gamma(\tau)) \, d\rho^i(s,\tau,\gamma) \label{eq:thm15.5.22} \\
    =& \int_{[0,T] \times [0,T] \times \Gamma_i} \#\{ (\tau, \gamma(\tau)) \in A \times B \} \, d\rho^i(s,\tau,\gamma) \label{eq:thm15.5.23} \\
    =& \int_{{\cal X}_T} \#\{ (\tau, \gamma_i(\tau)) \in A \times B \} \, d\rho(\gamma) \label{eq:thm15.5.24} \\
    =& \int_{{\cal X}_T} \#\{ t \in A \mid \lim_{\tau \to t^-} \gamma_i(\tau) \} \, d\rho(\gamma) \label{eq:thm15.5.25}
\end{align}
where
\eqref{eq:thm15.5.21} follows from Assumption \cref{assum:target_set}, \eqref{eq:sigma&eta}, and the fact that $B \subset S_{(i,i')}$;
\eqref{eq:thm15.5.22} follows from Corollary \cref{cor:rho_endpoints};
\eqref{eq:thm15.5.23} is because of Assumption \cref{assum:no_scuffing};
\eqref{eq:thm15.5.24} follows from \eqref{eq:rho};
\eqref{eq:thm15.5.25} is because all $\gamma_i \in \Gamma_i$ are absolutely continuous.

\end{IEEEproof}

Theorem \cref{thm:hle2traj} illustrates that measures satisfying \hle correspond to trajectories of the convexified inclusion, $\dot{\gamma}_i(t) \in \textrm{ conv } f(t,x(t), U)$, rather than the original specified dynamics within each mode of the system. 
To ensure that the there is no gap between the original specified dynamics and the solutions that correspond to the convexified inclusion, we assume the following condition:
\begin{assum}
\label{assum:F_convex}
The set $F_i(t,x,U)$ is convex for all $t$, $x$, \emph{or} $i \in {\cal I}$ or $F_i$ is control affine.
\end{assum}
\noindent Either condition in the previous assumption is sufficient to ensure that measures satisfying \hle correspond exactly to trajectories described according to Algorithm \cref{alg:1} \cite{vinter1993convex}. 
Finally, notice that Corollary \cref{cor:nu2u} provides a link between the solution measures and the underlying control input, which leads to a method capable of performing control synthesis:
\begin{corollary}
\label{cor:nu2u}
Suppose the dynamics of the hybrid system in each mode is control affine, i.e.,
    \begin{equation}
    F_i(t,x,u) = f_i(t,x) + g_i(t,x)u
    \end{equation}
    for all $t$, $x$, $u$, and $i \in {\cal I}$, where $f_i:\R \times X_i \to \R^{n_i}$ and $g_i:\R \times X_i \to \R^{n_i \times m}$.
Let $\nu^i_{u \mid t,x}$ and $\rho$ be defined as in \eqref{eq:disintegration} and \eqref{eq:rho}, respectively. Then
$t \mapsto ( \theta(t), \int_U u\, d\nu^{\lambda(\theta(t))}_{u \mid t, \theta_{\lambda(\theta(t))}} (u) )$
is an admissible pair for all $\theta \in \spt(\rho)$, where
\begin{equation}
\int_U u\, d\nu^i(u \mid t, x) :=
\begin{bmatrix}
\int_U [u]_1 \, d\nu^i (u \mid t,x)\\
\int_U [u]_2 \, d\nu^i (u \mid t,x)\\
\vdots\\
\int_U [u]_m \, d\nu^i (u \mid t,x)
\end{bmatrix}
\end{equation}
is an $m \times 1$ real vector for each $t$, $x$, and $i \in {\cal I}$.
\end{corollary}

\begin{IEEEproof}
For any $\theta \in \spt(\rho)$, we have:
\begin{equation}
\dot{\theta}_i(t) = f_i(t,\theta_i(t)) + g_i(t,\theta_i(t)) \cdot \int_U u\, d\nu^i_{u \mid t, \theta_i(t)}(u)
\end{equation}
for almost every $t \in [0,T]$.
Since $\nu^i(\cdot \mid t,x)$ is a stochastic kernel and $U$ is convex, we know
\begin{equation}
\int_U u\, d\nu^i_{u \mid t,\theta_i(t)}(u) \in U
\end{equation}
for all $i \in {\cal I}$.
Therefore, 
$t \mapsto ( \theta(t), \int_U u\, d\nu^{\lambda(\theta(t))}_{u \mid t, \theta_{\lambda(\theta(t))}} (u) )$
is an admissible pair.
\end{IEEEproof}

  \section{Infinite Dimensional Linear Program}
\label{sec:abs_prob}
This section reformulates $(OCP)$ as an infinite-dimension al linear program over the space of measures, proves it computes the solution to $(OCP)$, and illustrates how its solution can be used for control synthesis.

We first define $\mu^i_0$ to be dirac measure supported at $x_0$:
\begin{equation}
\label{eq:mu0}
    \mu_0^i = 
    \begin{cases}
        \delta_{x = x_0}, & \text{if } x_0 \in X_i; \\
        0, & \text{otherwise}
    \end{cases}
\end{equation}

Define the optimization problem $(P)$ as:
\begin{flalign} 
&& \inf_\Gamma & \phantom{4} \sum_{i\in {\cal I}}\langle \mu^i, h_i \rangle + \sum_{i\in {\cal I}} \langle \mu_T^i, H_i \rangle && (P) \nonumber \\
&& \text{s.t.} & \phantom{4} \delta_0 \otimes \mu_0^i + \Li' \mu^i + \sum_{(i',i) \in {\cal E}} \tilde{R}_{(i',i)\#} \mu^{S_{(i',i)}} = \delta_T \otimes \mu_T^i + \sum_{(i,i') \in {\cal E}} \mu^{S_{(i,i')}},  &&\forall i \in {\cal I}, \nonumber \\
&&& \phantom{4} \mu^i, \mu_T^i \geq 0 && \forall i \in {\cal I}, \nonumber\\
&&& \phantom{4} \mu^{S_e} \geq 0 && \forall e \in {\cal E} \nonumber
\end{flalign}
where the infimum is taken over a tuple of measures $\Gamma = (\mu^{\cal I}, \mu_T^{\cal I}, \mu^{\cal S}) \in \M_+([0,T]\times \mathcal{D} \times U) \times \M_+(X_T) \times \M_+([0,T] \times \mathcal{S})$ and for each mode $i \in {\cal I}$, $\mu_0^i$ is defined as in \eqref{eq:mu0}.
The dual to problem $(P)$ is given as:
\begin{flalign} 
&& \sup_v & \phantom{4} \sum_{i\in {\cal I}}\langle \mu_0^i(x), v_i(0,x) \rangle && (D) \nonumber \\
&& \text{s.t.} & \phantom{4} \Li v_i(t,x) + h_i(t,x,u) \geq 0,  &&\forall i \in {\cal I},\,(t,x,u) \in [0,T] \times X_i \times U \nonumber \\
&&& \phantom{4} v_i(T,x) \leq H_i(x) &&\forall i \in {\cal I},\, x \in X_{T_i} \nonumber \\
&&& \phantom{4} v_i(t,x) \leq v_{i'}(t,R_{(i,i')}(x)) &&\forall (i,i') \in \mathcal{E},\, (t,x) \in [0,T]\times X_i \nonumber
\end{flalign}
where the supremum is taken over the function $v \in C^1([0,T]\times \mathcal{D})$ and for each mode $i \in {\cal I}$, $\mu_0^i$ is defined as in \eqref{eq:mu0}.
Again, for notational convenience, we denote the $i \in {\cal I}$ slice of $v$ using subscript $i$ (i.e. for every $i \in {\cal I}$ and $(t,x) \in [0,T]\times X_i$, let $v_i(t,x) = v(t,x,i)$).

Next, we have the following result:
\begin{theorem}
There is no duality gap between $(P)$ and $(D)$.
\end{theorem}

\begin{IEEEproof}
The proof follows from \cite[Theorem~3.10]{anderson1987linear}.
\end{IEEEproof}
\Later{Do we need to explain?}

Next, we illustrate the $(P)$ is well-posed by proving the existence of optimal solution:
\begin{lem}
\label{lem:q*}
If $(P)$ is feasible, the minimum to $(P)$, $p^*$, is attained.
\end{lem}
\begin{IEEEproof}
 We prove $(i)$ in a manner similar to that employed in \cite[Theorem 2.3(i)]{lasserre2008nonlinear}. Let $(\mu^{\cal I}, \mu_T^{\cal I}, \mu^{\cal S})$ be a feasible solution to $(P)$. By choosing test functions $v_i=1$ and $v_i = T-t$ for all $i \in {\cal I}$, we may show the tuple of measures $(\mu^{\cal I}, \mu_T^{\cal I}, \mu^{\cal S})$ belongs to the unit ball $B_1$ of $\M([0,T]\times {\cal D} \times U) \times \M(X_T) \times \M([0,T] \times {\cal S})$. By Banach-Alaoglu theorem, $B_1$ is weak-* sequentially compact. Since the operators $\Li'$ and $R_{e\#}$ are continuous, the set of $(\mu^{\cal I}, \mu_T^{\cal I}, \mu^{\cal S})$ satisfying {\hle} is a closed subset of $B_1 \cap \M_+([0,T]\times {\cal D} \times U) \times \M_+(X_T) \times \M_+([0,T] \times {\cal S})$, and therefore is also weak-* sequentially compact. Since the linear functional to be minimized is continuous, $(P)$ is solvable.
\end{IEEEproof}

Now we prove that $(P)$ solves $(OCP)$:

\begin{theorem}
\label{thm:P}
Let $(P)$ be feasible and suppose $h_i(t,x,\cdot)$ is convex for all $i \in {\cal I}$ and $(t,x) \in [0,T] \times X_i$. 
Then $(P)$ solves $(OCP)$, i.e., $p^* = J^*$.
\end{theorem}


\begin{IEEEproof}
 We first prove $p^*$ is a lower bound of $J^*$ by showing that there exists a feasible solution of $(P)$ that achieves $J^*$ as the cost, and then use Theorem \cref{thm:hle2traj} to show $p^*$ cannot be less than $J^*$.

Suppose $(x^*,u^*)$ is an optimal admissible pair to $(OCP)$. 
By Lemma \cref{lem:traj2hle}, its initial measures, occupation measures, terminal measures and guard measures, denoted as $(\xi_0^{\cal I}, \xi^{\cal I}, \xi_T^{\cal I}, \xi^{\cal S})$, are supported on proper domains and satisfy \eqref{eq:hle}. Furthermore, $\xi_0^i = \mu_0^i$ for any $i \in {\cal I}$. 
Therefore, $(\xi^{\cal I}, \xi_T^{\cal I}, \xi^{\cal S})$ is a feasible solution to $(P)$ with cost $J(x^*,u^*) = J^*$, and $p^* \leq J^*$ follows.

We next prove $p^* \geq J^*$. 
Suppose $(\mu^{{\cal I}*}, \mu_T^{{\cal I}*}, \mu^{{\cal S}*})$ are an optimal solution to $(P)$ which exists according to Lemma \cref{lem:q*}.
The optimal tuple satisfies \eqref{eq:hle}.
By Theorem \cref{thm:hle2traj}, there exists a probability measure $\rho \in \M_+({\cal X}_T)$ such that $\\mu^{i*}_{t,x}$ coincides with the occupation measures of a family of admissible trajectories in the support of $\rho$, when restricted to mode $i$.

For the sake of simplicity, we abuse notation in the remainder of this proof and define $[\hat{u}_i(t,x)]_j := \int_U [u]_j\, d\nu^{i*}_{u \mid t,x}(u)$,
for any $i \in {\cal I}$ and $j \in \{1, \cdots, m\}$. 
We have
\begin{align}
q^* = & \sum_{i \in {\cal I}} \left( \int_{[0,T] \times X_i \times U} h_i(t,x,u)\, d\mu^{i*}(t,x,u) + \int_{X_T} H_i(x)\, d\mu_T^{i*}(x) \right) \nonumber \\
=& \sum_{i \in {\cal I}} \left( \int_{[0,T] \times X_i \times U} h_i(t,x,u)\, d\nu^{i*}_{u \mid t,x}(u)\, d\mu^{i*}_{x \mid t}(x)\, dt + \int_{X_T} H_i(x)\, d\mu_T^{i*}(x) \right) \nonumber \\
\geq & \sum_{i \in {\cal I}} \left( \int_{[0,T] \times X_i} h_i(t,x,\hat{u}_i(t,x))\, d\mu^{i*}_{x \mid t}(x)\, dt + \int_{X_T} H_i(x)\, d\mu^{i*}(x \mid T) \right) \label{eq:hconvex}\\
=& \sum_{i \in {\cal I}} \left( \int_{[0,T]} \int_{{\cal X}_T} h_i\left(t,\theta_i(t),\hat{u}_i(t,\theta_i(t))\right)\, d\rho(\theta)\, dt + \int_{{\cal X}_T} H_i(\theta_i(T))\, d\rho(\theta) \right)  \label{eq:subsrho}\\
=& \int_{{\cal X}_T} \sum_{i \in {\cal I}} \left( \int_{[0,T]} h_i\left(t,\theta_i(t),\hat{u}_i\left(t,\theta_i(t)\right)\right)\, dt + H_i(\theta_i(T)) \right)\, d\rho(\theta) \label{eq:fubini}\\
=& \int_{{\cal X}_T} \left( \int_{[0,T]} h_{\lambda(\theta(t))}\left(t,\theta_{\lambda(\theta(t))}(t), \hat{u}_{\lambda(\theta(t))}(t,\theta_{\lambda(\theta(t))}(t))\right)\, dt + H_{\lambda(\theta(T))}\left(\theta_{\lambda(\theta(T))}(T)\right) \right)\, d\rho(\theta) \label{eq:subslambda}\\
=& \int_{{\cal X}_T} J\left(\theta(\cdot), \hat{u}_{\lambda(\theta(\cdot))}(\cdot, \theta_{\lambda(\theta(\cdot))}(\cdot))\right)\, d\rho(\theta) \label{eq:Jtheta} \\
\geq & J^* \label{eq:<=J*}
\end{align}
where \eqref{eq:hconvex} is obtained from the convexity of $h_i(t,x,\cdot)$ and the fact that $\nu^{i*}_{u \mid t,x}$ is a probability measure;
\eqref{eq:subsrho} is from Theorem \cref{thm:hle2traj};
\eqref{eq:fubini} is from Fubini's Theorem;
\eqref{eq:subslambda} is because we let $h_i=0$ where $\theta_i(t)$ is undefined;
\eqref{eq:Jtheta} is because $(\theta(\cdot), \hat{u}_{\lambda_\theta(\cdot)}(\cdot, \theta_{\lambda_\theta(\cdot)}(\cdot)))$ is an admissible pair (according to Corollary \cref{cor:nu2u});
\eqref{eq:<=J*} is because $\rho$ is a probability measure.
\end{IEEEproof}

The previous result provides an extension of the weak formulation in \cite{lasserre2008nonlinear} to hybrid systems, and ensures $(P)$ can be solved to find a solution to $(OCP)$ in a convex manner. 
Next we describe how to perform control synthesis with the solution of $(P)$.

\begin{theorem}
\label{thm:radon-nikodym}
Suppose $(P)$ is feasible, $h_i(t,x,\cdot)$ is convex for all $i \in {\cal I}$ and $(t,x) \in [0,T] \times X_i$, the dynamics of the hybrid system in each mode is control affine, i.e.,
    \begin{equation}
    F_i(t,x,u) = f_i(t,x) + g_i(t,x)u
    \end{equation}
    for all $t$, $x$, $u$, and $i \in {\cal I}$, where $f_i:\R \times X_i \to \R^{n_i}$ and $g_i:\R \times X_i \to \R^{n_i \times m}$,
    and the optimal trajectory $x^*$ is unique $dt$-almost everywhere.
Let $\Gamma^* = (\mu^{{\cal I}*}, \mu_T^{{\cal I}*}, \mu^{{\cal S}*})$ be a vector of measures that achieves the infimum of $(P)$, then
\begin{enumerate}[label=(\alph*)]
    \item One can decompose $\mu^{{\cal I}*}$ in each mode $i \in {\cal I}$ as:
    \begin{equation}
    \label{eq:decomp}
    d\mu^{i*}(t,x,u) = d\nu^{i*}_{u \mid t,x}(u) d\mu^{i*}_{t,x}(t,x) = d\nu^{i*}_{u \mid t,x}(u)\, d\mu^{i*}_{x \mid t}(x)\, dt \qquad \forall i \in {\cal I}
    \end{equation}
    Moreover, $\mu^{i*}_{t,x}(t,x)$ coincides with the occupation measures of $x^*$ in each mode $i \in {\cal I}$ almost everywhere.
    \item  For each $i \in {\cal I}$, $j \in \{1, \cdots, m\}$, define:
    \begin{equation}
    \label{eq:lemmaP4}
    [\hat{u}_i(t,x)]_j := \int_U [u]_j\, d\nu^{i*}_{u \mid t,x}(u),
    \end{equation}
    for every point $(t,x)$ in the support of $\mu^{i*}_{t,x}$, where $d\nu^{i*}$ is as in \eqref{eq:decomp}.
    Then $\hat{u}_i(t,x) \in U$ for each $i \in {\cal I}$ and $(t,x)$ in the support of $\mu^{i*}_{t,x}$, and
    \begin{equation}
    \label{eq:lemmaP42}
    J\left( x^*(\cdot), \hat{u}_{\lambda(x^*(\cdot))}\left(\cdot, x^*_{\lambda(x^*(\cdot))}(\cdot)\right) \right) = J^*
    \end{equation}
    \item There exists a feedback control law, $\tilde{u}_i \in L^1([0,T]\times X_i, U)$ in each mode $i \in {\cal I}$, such that:
\begin{equation}
\label{eq:radon}
 [\tilde{u}_i(t,x)]_j \cdot \int_U\, d\mu^{i*}(t,x,u) = \int_U [u]_j \,d\mu^{i*}(t,x,u) 
\end{equation}
for each $j \in \{1, \cdots, m\}$.

Moreover, if we let $\tilde{u}(t,x,i) := \tilde{u}_i(t,x)$ for all $i \in {\cal I}$ and $(t,x) \in [0,T] \times X_i$, then $\tilde{u}$ is an optimal feedback control law, i.e.,
\begin{equation}
J\left(x^*(\cdot), \tilde{u}(\cdot, x^*(\cdot))\right) = J^*
\end{equation}
\end{enumerate}
\end{theorem}

\Later{Indentation is different here... Can you fix this?}
\begin{IEEEproof}
\hfill 
\begin{enumerate}[label=(\alph*)]
\item First note that the decomposition of $\mu^{{\cal I}*}$ exists as a result of Theorem \cref{thm:hle2traj}. 
Using the notation and result within the proof of Theorem \cref{thm:P}:
\begin{equation}
\label{eq:J=J*}
J \left(\theta(\cdot), \hat{u}_{\lambda(\theta(\cdot))}(\cdot, \theta_{\lambda(\theta(\cdot))}(\cdot))\right) = J^*
\end{equation}
for any $\theta \in \spt(\rho)$. 
Therefore every admissible pair $\left(\theta(\cdot), \hat{u}_{\lambda(\theta(\cdot))}(\cdot, \theta_{\lambda(\theta(\cdot))}(\cdot))\right)$ must be optimal.
Since the optimal trajectory $x^*$ is assumed to be unique $dt$-almost everywhere, we have
\begin{equation}
\label{eq:theta=x*}
\theta(t) = x^*(t) \text{ for almost everywhere } t\in [0,T], \quad \forall \theta(\cdot) \in \spt(\rho)
\end{equation}
According to Theorem \cref{thm:hle2traj}, $\mu^{i*}_{t,x}$ coincides with the occupation measure of the family of admissible trajectories in $\spt(\rho)$. 
Note the similarity between \eqref{eq:def:mu_traj} and \eqref{eq:def:zeta}, therefore  $\mu^{i*}_{t,x}$ coincides with the occupation measure of $x^*$ in each mode $i \in {\cal I}$ almost everywhere.

\item This follows from Corollary \cref{cor:nu2u}, \eqref{eq:J=J*}, and \eqref{eq:theta=x*}.

\item We prove the first result using Radon-Nikodym Theorem, and the second result can be shown using Theorem \cref{thm:P} and the uniqueness of $\tilde{u}_i$.

For notational convenience, we define measures $\mu^{i*}_{t,x}, \eta^{i*}_j \in \M_+([0,T] \times X_i)$ as
\begin{equation}
\begin{split}
    d\mu^{i*}_{t,x} &:= \int_{U_i} d\mu^{i*}(t,x,u)\\
    d\eta^{i*}_j &:= \int_{U_i} [u]_j \, d\mu^{i*}(t,x,u)
\end{split}
\end{equation}
For each mode $i \in {\cal I}$, $\mu^{i*}$ is $\sigma$-finite since it's a Radon measure defined over a compact set, therefore $\mu^{i*}_{t,x}$ and $\eta^{i*}_j$ are also $\sigma$-finite. To apply Radon-Nikodym theorem, we need to show $\eta^{i*}_j$ is absolutely continuous with respect to $\mu^{i*}_{t,x}$:

For any Borel measurable set $A \times B \subseteq [0,T]\times X_i$ such that $\mu^{i*}_{t,x}(A \times B) = 0$, we have
\begin{equation}
\int_{A \times B \times U_i} d\mu^{i*}(t,x,u) = \mu^{i*}_{t,x}(A \times B) = 0
\end{equation}
Note that $\mu^{i*}$ is unsigned measure, therefore $\mu^{i*}$ is zero on any measurable subset of $A \times B \times U_i$. This implies
\begin{equation}
\eta^{i*}_j(A \times B) = \int_{A \times B \times U_i} [u]_j\, d\mu^{i*}(t,x,u) = 0
\end{equation}

As a result of Radon-Nikodym theorem, there exist functions $[\tilde{u_i}]_j \in L^1([0,T] \times X_i, \R_+)$ for each $j \in \{1, \cdots, m_i\}$, such that Equation \eqref{eq:radon} is satisfied. Such functions $[\tilde{u}_i]_j$ are unique $\mu^{i*}_{t,x}$-almost everywhere.

Now we want to show there is a version of $\tilde{u}_i$ whose range is a subset of $U$.
Note that we can disintegrate each $\mu^{i*}$ and define $\hat{u}_i$ as in Theorem \cref{thm:P}, hence Equation \eqref{eq:radon} becomes
\begin{equation}
\begin{split}
[\tilde{u}_i(t,x)]_j \, d\mu^{i*}_{t,x}(t,x) =& \left( \int_U [u]_j \, d\nu^{i*}_{u \mid t,x}(u) \right) \, d\mu^{i*}_{t,x}(t,x)\\
=& [\hat{u}_i(t,x)]_j \, d\mu^{i*}_{t,x}(t,x)
\end{split}
\end{equation}
Since $\hat{u}_i(t,x) \in U$ for all $(t,x)$ in the support of $\mu^{i*}_{t,x}$, we may choose $\tilde{u}_i = \hat{u}_i$ such that Equation \eqref{eq:radon} is satisfied and $\tilde{u}_i \in L^1([0,T]\times X_i, U)$.

Finally, let $x^*$ be the optimal trajectory, and it follows directly from definition that
\begin{equation}
\tilde{u}(t,x^*(t)) = \tilde{u}\left(t,x^*_{\lambda({x^*}(t))}(t),\lambda({x^*}(t))\right) = \tilde{u}_{\lambda({x^*}(t))}\left(t,x^*_{\lambda({x^*}(t))}(t)\right)
\end{equation}
Using Equation \eqref{eq:lemmaP42}, we know
\begin{equation}
J\left( x^*(\cdot), \tilde{u}(\cdot, x^*(\cdot)) \right) = J\left( x^*(\cdot), \hat{u}_{\lambda({x^*}(t))}\left( \cdot, x^*_{\lambda({x^*}(\cdot))}(\cdot) \right) \right) = J^*
\end{equation}
\end{enumerate}
\end{IEEEproof}

Notice that the second result in Theorem \cref{thm:radon-nikodym} requires that we be able to construct the condition measure $\nu^*$ to be able to construct a feedback controller. 
In contrast, the third result within Theorem \cref{thm:radon-nikodym} illustrates how one can construct a feedback controller by computing the Radon-Nikodym derivative using the optimal measures from the solution to $(P)$. 
As we describe in the next section, this latter result can be utilized directly to construct a sequence of controllers that converge to the optimal control. 
Finally notice that in the hypothesis of Theorem \cref{thm:radon-nikodym} we do not assume the uniqueness of the optimal control law, i.e., there may exist different control laws $u_1$ and $u_2$, such that $J^* = J(x^*,u_1) = J(x^*,u_2)$. 
In this instance the admissible pairs $(x^*,u_1)$ and $(x^*,u_2)$ are both optimal, and we are interested in finding either of them.
Instead we only assume that the optimal trajectory is unique almost-everywhere.
As a result, the optimal solution to $(P)$ may not be unique.


  \section{Numerical Implementation}
\label{sec:implementation}
We compute a solution to the infinite-dimensional problem $(P)$ via a sequence of finite-dimensional approximations formulated as semidefinite programs (SDP)s.
These are generated by representing the measures in $(P)$ using a truncated sequence of moments and restricting the functions in $(D)$ to polynomials of finite degree.
As illustrated in this section, the solutions to any of the SDPs in this sequence can be used to synthesize an approximation to the optimal controllers.
A comprehensive introduction to such moment relaxations can be found in \cite{lasserre2009moments}.

To formulate this SDP relaxation, we restrict our interest to polynomial hybrid optimal control problems:
\begin{assum}
The functions $F_i$, $h_i$, and $H_i$ are polynomials, that is, $[F_i]_j \in \R[t,x,u]$, $h_i \in \R[t,x,u]$, and $H_i \in \R[x]$ for all $i \in {\cal I}$ and $j \in \{1, \cdots, n_i\}$.
\end{assum}

Note that in the notation $\R[t,x,u]$, we refer to $x$ as an indeterminate in $X_i$ with dimension $n_i$.
It should not be confused with trajectory of the hybrid system. 
In addition, for notational convenience, the dimension $n_i$ of $x$ is omitted when it is clear in context.

We also make assumptions about the sets $X_i$, $X_{T_i}$, $U$, and $S_e$:
\begin{assum}
$X_i$, $X_{T_i}$, $U$, and $S_e$ are semi-algebraic sets, i.e.,
\begin{align}
X_i &= \left\{x \in \R^{n_i} \mid h_{X_{i_j}} \geq 0, h_{X_{i_j}} \in \R[x], \forall j \in \{1, \cdots, n_{X_i}\}\right\},\\
X_{T_i} &= \left\{ x \in \R^{n_i} \mid h_{T_{i_j}} \geq 0, h_{T_{i_j}} \in \R[x], \forall j \in \{1, \cdots, n_{T_i}\}\right\},\\
U &= \left\{ u \in \R^m \mid h_{U_j} \geq 0, h_{U_j} \in \R[u], \forall j \in \{1, \cdots, n_U\} \right\},\\
S_{(i,i')} &= \left\{ x \in \partial X_i \mid h_{(i,i')_j} \geq 0, h_{(i,i')_j} \in \R[x], \forall j \in \{1, \cdots, n_{(i,i')}\} \right\}
\end{align}
for all $i \in {\cal I}$ and $(i,i') \in {\cal E}$.
\end{assum}
\noindent Since $X_i$ and $X_{T_i}$ are also compact, note that Putinar's condition (see \cite{lasserre2009moments}) is satisfied by adding the redundant constraint $M - \|x\|_2^2$ for some large enough $M$. 

To derive the SDP relaxation, we begin with a few preliminaries.
Any polynomial $p \in \R_k[x]$ can be expressed in the monomial basis as:
\begin{equation}
p(x) = \sum_{|\alpha| \leq k} p_\alpha x^\alpha = \sum_{|\alpha| \leq k} p_\alpha \cdot (x_1^{\alpha_1}\cdots x_n^{\alpha_n} )
\end{equation}
where $\alpha$ ranges over vectors of non-negative integers such that $|\alpha| = \sum_{i=1}^n \alpha_i \leq k$, and we denote $\text{vec}(p) = (p_\alpha)_{|\alpha| \leq k }$ as the vector of coefficients of $p$.
Given a vector of real numbers $y = (y_\alpha)$ indexed by $\alpha$, we define the linear functional $L_y: \R_k[x] \rightarrow \R$ as:
\begin{equation}
L_y(p) := \sum_\alpha p_\alpha y_\alpha
\end{equation}
Note that, when the entries of $y$ are moments of a measure $\mu$:
\begin{equation}
y_\alpha = \int x^\alpha \, d\mu(x),
\end{equation}
then
\begin{equation}
    \langle \mu, p \rangle = \int \left( \sum_\alpha p_\alpha x^\alpha \right)\, d\mu = L_y(p).
\end{equation}
If $|\alpha| \leq 2k$, the \emph{moment matrix}, $M_k(y)$, is defined as:
\begin{equation}
[M_k(y)]_{\alpha \beta} = y_{(\alpha + \beta)}
\end{equation}
Given any polynomial $h \in \R_{l}[x]$ with $l < k$, the \emph{localizing matrix}, $M_k(h,y)$, is defined as:
\begin{equation}
[M_k(h,y)]_{\alpha \beta} = \sum_{|\gamma| \leq l} h_\gamma y_{(\gamma + \alpha + \beta)}.
\end{equation}
Note that the moment and localizing matrices are symmetric and linear in moments $y$.

\subsection{LMI Relaxations and SOS Approximations}

An sequence of SDPs approximating $(P)$ can be obtained by replacing constraints on measures with constraints on moments.
Since $h_i$ and $H_i$ are polynomials, the objective function of $(P)$ can be written using linear functionals as $\sum_{i\in {\cal I}}L_{y_{\mu^i}}(h_i) + \sum_{i\in {\cal I}} L_{y_{\mu_T^i}}(H_i)$, where $y_{\mu^i}$ and $y_{\mu_T^i}$ are the sequence of moments of $\mu^i$ and $\mu_T^i$, respectively.
The equality constraints in $(P)$ can be approximated by an infinite-dimensional linear system, which is obtained by restricting to polynomial test functions: 
$v_i(t,x) \in \R[t,x]$, for any $i \in {\cal I}$.
The positivity constraints in $(P)$ can be replaced with semidefinite constraints on moment and localizing matrices, which guarantees the existence of Borel measures defined on proper domains \cite[Theorem~3.8]{lasserre2009moments}.

A finite-dimensional SDP is then obtained by truncating the degree of moments and polynomial test functions to $2k$.
Let $\Xi_{\cal I} = \coprod_{i \in {\cal I}} \mu^i$,
$\Xi_{\cal E} = \coprod_{e \in {\cal E}} \mu^{S_e}$, $\Xi_T = \coprod_{i \in {\cal I}} \mu_T^i$, and $\Xi = \Xi_{\cal I} \bigcup \Xi_{\cal E} \bigcup \Xi_T$.
Let $(y_{k,\xi})$ be the sequence of moments truncated to degree $2k$ for each $(\xi, i) \in \Xi$, and let $\mathbf{y}_k$ be a vector of all the sequences $(y_{k,\xi})$.
The equality contraints in $(P)$ can then be approximated by a finite-dimensional linear system:
\begin{equation}
A_k(\mathbf{y}_k) = b_k
\end{equation}
Define the $k$-th relaxed SDP representation of $(P)$, denoted $(P_k)$, as
\begin{flalign} 
&& \inf & \phantom{4} \sum_{i\in {\cal I}}L_{y_{k,\mu^i}}(h_i) + \sum_{i\in {\cal I}} L_{y_{k,\mu_T^i}}(H_i) && (P_k) \nonumber \\
&& \text{s.t.} & \phantom{4} A_k(\mathbf{y}_k) = b_k, \nonumber\\
&&& \phantom{4} M_k(y_{k,\xi}) \succeq 0 && \forall (\xi,i) \in \Xi, \nonumber\\
&&& \phantom{4} M_{k_{X_{i_j}}}(h_{X_{i_j}}, y_{k, \mu^i}) \succeq 0 && \forall (j, i) \in \{1,\cdots, n_{X_i}\} \times {\cal I}, \nonumber\\
&&& \phantom{4} M_{k_{U_{i_j}}}(h_{U_j}, y_{k,\mu^i}) \succeq 0 && \forall (j,i) \in \{1,\cdots,n_{U_i}\} \times {\cal I}, \nonumber\\
&&& \phantom{4} M_{k_{S_{e_j}}}(h_{e_j}, y_{k,\xi}) \succeq 0 && \forall (j,\xi,e) \in \{1,\cdots,n_e\}\times \Xi_{\cal E}, \nonumber\\
&&& \phantom{4} M_{k_{T_{i_j}}}(h_{T_{i_j}}, y_{k,\xi}) \succeq 0 && \forall (j,\xi,i) \in \{1,\cdots,n_{T_i}\} \times \Xi_T, \nonumber\\
&&& \phantom{4} M_{k-1}(h_\tau, y_{k,\xi}) \succeq 0 && \forall (\xi,i) \in \Xi_{\cal I} \bigcup \Xi_{\cal E} \nonumber
\end{flalign}
where the infimum is taken over the sequences of moments $(y_{k,\xi})$ for each $(\xi,i) \in \Xi$, $h_\tau = t(T-t)$, $k_{X_{i_j}} = k - \lceil \text{deg}(h_{X_{i_j}})/2 \rceil$, $k_{U_{i_j}} = k - \lceil \text{deg}(h_{U_{i_j}})/2 \rceil$, $k_{S_{e_j}} = k - \lceil \text{deg}(h_{e_j})/2 \rceil$, $k_{T_{i_j}} = k - \lceil \text{deg}(h_{T_{i_j}})/2 \rceil$, and $\succeq$ denotes positive semidefiniteness of matrices.

The dual of $(P_k)$ is a Sums-of-Squares (SOS) program denoted by $(D_k)$ for each $k \in \N$, which is obtained by first restricting the optimization space in $(D)$ to the polynomial functions with degree truncated to $2k$ and by then replacing the non-negativity constraints in $(D)$ with SOS constraints.
For notational convenience, we let $x_i$ be the indeterminate that corresponds to $X_i$.
Define $Q_{2k}(h_{T_{i_1}}, \cdots, h_{T_{i_{n_{T_i}}}}) \subset \R_{2k}[x_i]$ to be the set of polynomials $l \in \R_{2k}[x_i]$ expressible as
\begin{equation}
l = s_0 + \sum_{j=1}^{n_{T_i}} s_j h_{T_{i_j}}
\end{equation}
for some polynomials $\{s_j\}_{i=0}^{n_{T_i}} \subset \R_{2k}[x_i]$ that are sums of squares of other polynomials. 
Every such polynomial is clearly non-negative on $X_{T_i}$.
Similarly, we define $Q_{2k}(h_\tau, h_{X_{i_1}}, \cdots, h_{X_{i_{n_{X_i}}}}, h_{U_1}, \cdots, h_{U_{n_U}}) \subset \R_{2k}[t,x_i,u]$, and $Q_{2k}(h_\tau, h_{(i,i')_1}, \cdots, h_{(i,i')_{n_{(i,i')}}}) \subset \R_{2k}[t,x_i]$ for each $i \in {\cal I}$ and $(i,i') \in {\cal E}$.
Therefore $k$-th relaxed SDP representation of $(D)$, denoted $(D_k)$ is given as
\begin{flalign} 
&& \sup & \phantom{4} \sum_{i\in {\cal I}}\langle \mu_0^i,v_i(0,\cdot) \rangle && (D_k) \nonumber \\
&& \text{s.t.} & \phantom{4} \Li v_i + h_i \in Q_{2k}(h_\tau, h_{X_{i_1}}, \cdots, h_{X_{i_{n_{X_i}}}}, h_{U_1}, \cdots, h_{U_{n_U}}) && \forall i \in {\cal I}, \nonumber\\
&&& \phantom{4} -v_i(T,\cdot) + H_i \in Q_{2k}(h_{T_{i_1}}, \cdots, h_{T_{i_{n_{T_i}}}}) && \forall i \in {\cal I}, \nonumber\\
&&& \phantom{4} v_{i'}\circ \tilde{R}_{(i,i')} - v_i \in Q_{2k}(h_\tau, h_{(i,i')_1}, \cdots, h_{(i,i')_{n_{(i,i')}}}) && \forall (i,i') \in {\cal E}, \nonumber
\end{flalign}
where the supremum is taken over polynomials $v_i \in \R_{2k}[t,x]$ for all $i \in {\cal I}$.
\Later{where is the index $i$ in the definition of polynomials??}

We first prove that these pair of problems are well-posed:
\begin{theorem}
\label{theorem:duality_Pk}
For each $k \in \mathbf{N}$, if $(P_k)$ is feasible, then there is no duality gap between $(P_k)$ and $(D_k)$.
\end{theorem}
\begin{IEEEproof}
This can be proved using Slater's condition (see \cite{boyd2004convex}), which involves noting that $(D_k)$ is bounded below, and then arguing the feasible set has an interior point.
\end{IEEEproof}

Next, we describe how to extract a polynomial control law from the solution of $(P_k)$.
Given moment sequences truncated to $2k$, we want to find an appropriate feedback control law $u_{k,i}^*$ in each mode $i \in {\cal I}$ with components $[u_{k,i}^*]_j \in \mathbb{R}[t,x]$, such that the analogue of \eqref{eq:radon} is satisfied, i.e.,
\begin{equation}
\label{eq:radon_discrete}
\int_{[0,T]\times X_i} t^{\alpha_0}x^{\alpha} \cdot [u_{k,i}^*]_j(t,x) \, \int_U d\mu^{i*}_k(t,x,u) = \int_{[0,T] \times X_i} t^{\alpha_0}x^{\alpha} \cdot \int_U [u]_j \, d\mu^{i*}_k(t,x,u)
\end{equation}
for all $i \in \mathcal{I}$, $j \in \{1,\cdots, m\}$, and $(\alpha_0, \alpha) \in \N \times \N^{n_i}$ satisfying $\sum_{l=0}^n \alpha_l \leq k$, $\alpha_l \geq 0$. \Later{$\alpha_0$ and $\alpha$ might be confusing.}
Here $\mu^{i*}_k$ is any measure whose truncated moments match $y_{\mu_i}^*$.
In fact, when constructing a polynomial control law from the solution of $(P_k)$, these linear equations written with respect to the coefficients of $[u_{k,i}^*]_j$ are expressible in terms of the optimal solution $y_{k,\mu^i}^*$.

To see this, define $(t,x)$-moment matrix of $y_{k,\mu^i}^*$ as:
\begin{equation}
    \left[ M_k^{(t,x)}(y_{k,\mu^i}^*) \right]_{(\alpha_0,\alpha) (\beta_0,\beta)} = L_{y_{k,\mu^i}^*} (t^{\alpha_0+\beta_0}x^{\alpha+\beta}u^\mathbf{0}) = \left[ y_{k,\mu^i}^* \right]_{(\alpha_0,\alpha,\mathbf{0}) (\beta_0,\beta,\mathbf{0})}
\end{equation}
for all $i \in {\cal I}$, and $(\alpha_0, \alpha, \mathbf{0}), (\beta_0, \beta, \mathbf{0}) \in \N \times \N^{n_i} \times \{0\}^m$ satisfying $\sum_{l=0}^n \alpha_l \leq k$, $\alpha_l \geq 0$, $\sum_{l=0}^n \beta_l \leq k$, $\beta_l \geq 0$.
Also define a vector $b_k^{j}$ as
\begin{equation}
\left[ b_k^j(y_{k,\mu^i}^*) \right]_\alpha = L_{y_{k,\mu^i}^*}(t^{\alpha_0}x^{\alpha} \cdot [u]_j)
\end{equation}
for all $j \in \{1, \cdots, m\}$, and $(\alpha_0, \alpha) \in \N \times \N^{n_i}$ satisfying $\sum_{l=0}^n \alpha_l \leq k$, $\alpha_l \geq 0$.
Direct calculation shows Equation \eqref{eq:radon_discrete} is equivalent as the following linear system of equations:
\begin{equation}
\label{eq:ctrl_extract}
M_k^{(t,x)}(y_{k,\mu^i}^*) \, \text{vec}([u_{k,i}^*]_{j}) = b_k^{j}(y_{k,\mu^i}^*)
\end{equation}
To extract the coefficients of the controller, one needs only to compute the generalized inverse of $M_k^{(t,x)}(y_{k,\mu^i}^*)$, which exists since it is positive semidefinite.
\Later{Do we need to prove this? btw what's ``generalized inverse''? Is it pseudo-inverse?}
Note that the degree of the extracted polynomial control law is dependent on the relaxation order $k$. Higher relaxation orders lead to higher degree controllers.

\subsection{Convergence of Relaxed Problems}
Next, we prove the convergence of the pair of approximations:
\begin{theorem}
\label{thm:P_k}
Let $p_k^*$ and $d_k^*$ denote the infimum of $(P_k)$ and supremum of $(D_k)$, respectively. Then $\{p_k^*\}_{k=1}^\infty$ and $\{d_k^*\}_{k=1}^\infty$ converge monotonically from below to the optimal value of $(P)$ and $(D)$.
\end{theorem}

\begin{IEEEproof}
This theorem can be proved using a similar technique adopted in the proof of \cite[Theorem~4.2]{majumdar2014convex}. 
We first establish a lower found of $d_k^*$ by finding a feasible solution to $(D_k)$ for some $k$, and then show that there exists a convergent subsequence of $\{d_k^*\}_{k=1}^\infty$, by arguing the lower bound can be arbitrarily close to $d^*$ for large enough $k$.
Using \cref{theorem:duality_Pk}, we only need to prove $\{d_k^*\}_{k=1}^\infty$ converge monotonically from below to $d^*$.

Note that the higher the relaxation order $k$, the looser the constraint set of the optimization problem $(D_k)$, so $\{d_k^*\}_{k=1}^\infty$ is nond-ecreasing.

Suppose $v \in C^1([0,T] \times {\cal D})$ is feasible in $(D)$. For every $\epsilon>0$ and $i \in {\cal I}$, set
\begin{equation}
\tilde{v}_i(t,x) := v_i(t,x) + \epsilon t - (1+T)\epsilon
\end{equation}
Therefore, $\Li \tilde{v}_i = \Li v_i + \epsilon$, $\tilde{v}_i(T,x) = v_i(T,x) - \epsilon$, and it follows that $\coprod_{i \in {\cal I}} \tilde{v}_i$ is strictly feasible in $(D)$ with a margin at least $\epsilon$. 
Since $[0,T]\times X_i$ and $X_i$ are compact for every $i \in {\cal I}$, and by a generalization of the Stone-Weierstrass theorem that allows for the simultaneous uniform approximation of a function and its derivatives by a polynomial \cite{hirsch2012differential}, we are guaranteed the existence of polynomials $\hat{v}_i$, such that $\|\hat{v}_i - \tilde{v}_i\|_\infty < \epsilon$, and $\|\Li \hat{v}_i - \Li \tilde{v}_i\|_\infty < \epsilon$ for any $i \in {\cal I}$. 
By Putinar's Positivstellensatz \cite[Theorem~2.14]{lasserre2009moments}, those polynomials are strictly feasible for $(D_k)$ for a sufficiently large relaxation order $k$, therefore $d_k^* \geq \sum_{i \in {\cal I}}\hat{v}_i(0,x_0) \geq \sum_{i \in {\cal I}} \tilde{v}_i(0,x_0) - |{\cal I}|\epsilon$, where $|{\cal I}|$ is the number of elements in ${\cal I}$. 
Also, since $\tilde{v}_i(0,x_0) = v_i(0,x_0) - (1+T)\epsilon$, we have $d_k^* > \sum_{i \in {\cal I}} v_i(0,x_0)-(1+T+|{\cal I}|)\epsilon = d^*-(1+T+|{\cal I}|)\epsilon$, where $1+T+|{\cal I}| < \infty$ is a constant. 
Using the fact that $d^*$ is non-decreasing and bounded above by $d$, we know $\{d_k^*\}_{k=1}^\infty$ converges to $d$ from below.

\end{IEEEproof}

Finally we can prove that the sequence of controls extracted as the solution to the linear equation \eqref{eq:ctrl_extract} from the sequence of SDPs converges to the optimal control:

\begin{theorem}
Let $\{y_{k,\xi}^*\}_{(\xi,i) \in \Xi}$ be an optimizer of $(P_k)$, and let $\{\mu_k^{i*}\}_{i \in {\cal I}}$ be a set of measures such that the truncated moments of $\mu_k^{i*}$ match $y_{k,\mu^i}^*$ for each $i \in {\cal I}$.
In addition, for each $k \in \N$, let $u_{k,i}^*$ denote the controller constructed by \eqref{eq:ctrl_extract}, and $\tilde{u}_i$ is the optimal feedback control law in mode $i \in {\cal I}$ defined in \cref{thm:radon-nikodym}.
Then, there exists a subsequence $\{k_l\}_{l \in \N} \subset \N$ such that:
\begin{equation}
\label{eq:mom_conv}
\int_{[0,T] \times X_i} v_i(t,x) [u_{k_l,i}^*]_j(t,x) \, d\mu_{t,x; k_l}^{i*}(t,x) \xrightarrow{l \rightarrow \infty} \int_{[0,T] \times X_i} v_i(t,x) \, [\tilde{u}_i]_j(t,x) \, d\mu^{i*}_{t,x}(t,x)
\end{equation}
for all $i \in {\cal I}$, $v_i \in C^1([0,T] \times X_i)$, and $j \in \{1, \cdots, m\}$.
\end{theorem}

\begin{IEEEproof}
From the proof of \cite[Theorem~4.3]{lasserre2009moments}, if we complete each $y_{k,\xi}^*$ with zeros and make it an infinite vector, then there exists a $y_\xi^* \in l_\infty$ and a subsequence $\{k_l\}_{l \in \N}$ such that for each $(\xi, i) \in \Xi$, $\lim_{l \rightarrow \infty} y_{k_l, \xi}^* = y_\xi^*$ for the weak-* topology $\sigma(l_\infty, l_1)$ of $l_\infty$.
Moreover, as a result of \cite[Theorem~3.8(b)]{lasserre2009moments}, for each $(\xi, i) \in \Xi$, $y_\xi^*$ has a finite Borel representing measure, and this set of represented measures, which we denote by $(\mu^{\cal I*}, \mu_T^{\cal I*}, \mu^{\cal S*})$, is an optimizing vector of measures for $(P)$.

Next, consider any polynomial test function $t^{\alpha_0}x^{\alpha} \in \R[t,x]$, and let $r$ be its degree. 
Then Equation \eqref{eq:radon} and \eqref{eq:radon_discrete} are both true for $k >r$.
Therefore we only need to show
\begin{equation}
\int_{[0,T] \times X_i} t^{\alpha_0}x^{\alpha} \cdot \int_U [u]_j \, d\mu^{i*}_{k_l}(t,x,u) \xrightarrow{l \rightarrow \infty} \int_{[0,T] \times X_i} t^{\alpha_0}x^{\alpha} \cdot \int_U [u]_j \, d\mu^{i*}(t,x,u)
\end{equation}
for all $i \in {\cal I}$.
Define $w_i(t,x,u) := t^{\alpha_0}x^{\alpha} \cdot [u]_j \in \R[t,x,u]$, then
\begin{align}
    &\int_{[0,T] \times X_i \times U} w_i(t,x,u) \, \left( d\mu_{k_l}^{i*}(t,x,u) - d\mu^{i*}(t,x,u) \right) \nonumber \\
    =& \text{vec}(w_i) \left( y_{k_l,\mu^i}^* - y_{\mu^i}^* \right) \xrightarrow{l \rightarrow \infty} 0 \nonumber
\end{align}
where the last argument is true because $\text{vec}(w_i) \in l_1$. Note that the set of polynomials is dense in $C^1$, therefore Equation \eqref{eq:mom_conv} is true for all $C^1$ functions $v_i$.
\end{IEEEproof}

  \section{Extension to Free Final Time Problem}
\label{sec:extensions}
It is useful to sometimes consider the optimal control problem where the system state has to be driven to $X_T$ before a fixed time $T_0$, and not necessarily remain in $X_T$ afterwards (as opposed to reaching $X_T$ exactly at time $T$). 
We refer to this problem as the \emph{free terminal time problem}.
We adapt the notation $T$ to denote the first time a trajectory reaches $X_T$ (the terminal time), and an admissible pair $(x,u)$ can be redefined as follows.
Given a real number $T_0>0$ and a point $(x_0,j) \in {\cal D}$, if there exists a $T$ satisfying $0 < T \leq T_0$, a control $u:[0,T] \to U$, and a trajectory $x:[0,T] \to {\cal D}$ such that $(x,u)$ satisfies Algorithm \cref{alg:1}, and $x(T) \in X_T$ then $x$ is called an \emph{admissible trajectory}, $u$ is called an \emph{admissible control}, and the pair $(x,u)$ is called an \emph{admissible pair}.

In practice this formulation requires, we modify $(OCP)$ by adding in another constraint $0<T \leq T_0$ since $T$ is a free variable now.
The primal LP that solves the free terminal time is obtained by modifying the support of $\mu_T$ in $(P)$ to be $[0,T_0]\times X_T$, and substituting $\delta_T \otimes \mu_T$ with $\mu_T$ in its first constraint. The only modification to $(D)$ is that the second constraint is imposed for all time $t \in [0,T_0]$ instead of just at time $T$. 
All results from the previous sections can be extended to the free-terminal-time case with nearly identical proofs, and the numerical implementation follows in a straightforward manner.
When $h_i \equiv 1$ and $H_i \equiv 0$ for any $i \in {\cal I}$ with free terminal time, the $(OCP)$ can be interpreted as a minimum time problem \Later{should it be ``minimal'' instead of ``minimum''?}, where the optimal control problem must find an admissible pair $(x,u)$ such that the trajectory reaches the target set as quickly as possible.
  \section{Examples}
\label{sec:examples}
\Later{'Degree of relaxation? Relaxation order?'}
This section illustrates the performance of our approach using several examples.
Before proceeding, we begin by describing the numerical implementation of our algorithm.
First, for each of the systems described below, we denote the range space of control inputs in mode $i$ as $U_i$. 
In the implementation, we can always define $U := \prod_{i \in {\cal I}}U_i$ without causing any problems.
Second, for each of the optimal control problems, we implement our algorithm using the MOSEK\cite{mosek} numerical solver in MATLAB and generate a polynomial feedback control law.
Third, the trajectory is obtained by plugging the (saturated) polynomial control law back into the system dynamics in each mode and simulating forward using a standard ODE solver with event detection in MATLAB. 
Once the trajectory hits a guard, another ODE solver is initialized at a new point given by the associated reset map, and the simulation continues in the same way until terminal condition is satisfied.
Next, for the sake of comparison, all the examples are solved either analytically (when possible) or using GPOPS-II\cite{gpops} by iterating through a finite set of of possible transitions.
Notice that in this latter instance we must fix the possible sequences of transitions and provide an initial condition, since existing numerical hybrid optimal control algorithms require this information.
Finally, all of our experiments are performed on an Intel Xeon, 20 core, 2.60 GHz, 128 GB RAM machine.
\Later{do you want to mention that all of the code can be downloaded??}

\subsection{Hybridized Double Integrator}
The double integrator is a two-state, single-input linear system.
\Later{$x_1$ and $x_2$ shouldn't be confused with $\pi_{X_1}x$ and $\pi_{X_2}x$.}
Even though a standard double integrator is a non-hybrid system, we may hybridize it by dividing the domain into two parts, and defining an identity reset map between them  as described in Table \cref{tab:DI:vf} and Table \cref{tab:DI:guard}.

\begin{table}[!h]
\caption{Vector fields and domains of each of the modes of the hybridized double integrator}
\label{tab:DI:vf}
\centering
\begin{tabular}{c|c|c}
\hline
    Mode & $i=1$ & $i=2$ \\
\hline
    Dynamics & $\dot{x}(t) = 
    \begin{bmatrix}
        x_2 \\ 0
    \end{bmatrix} +
    \begin{bmatrix}
        0 \\ 1
    \end{bmatrix} u$ & 
    $\dot{x}(t) = 
    \begin{bmatrix}
        x_2 \\ 0
    \end{bmatrix} +
    \begin{bmatrix}
        0 \\ 1
    \end{bmatrix} u$ \\
\hline
    $X_i$ & $\{(x_1, x_2) \in \R^2 \mid x_1^2 + x_2^2 \leq 0.3\}$ & $\{(x_1, x_2) \in \R^2 \mid x_1^2 + x_2^2 \geq 0.3\}$\\
\hline
    $U_i$ & $[-1,1]$ & $[-1,1]$\\
\hline
\end{tabular}
\end{table}

\begin{table}[t]
\caption{Guards and reset maps of the hybridized double integrator.
The rows are modes in which a transition originates, and the columns are modes to which the transition goes.
}
\label{tab:DI:guard}
\centering
\begin{tabular}{c|c|c}
\hline
    & Mode 1 & Mode 2 \\
\hline
    Mode 1 & N/A & N/A \\
\hline
    Mode 2 &
    $\begin{aligned}
    &S_{(2,1)} = \{(x_1, x_2) \in \R^2 \mid x_1^2 + x_2^2 = 0.3\}\\
    &R_{(2,1)}(x) = x, \quad \forall x \in S_{(2,1)}
    \end{aligned}$
    & N/A \\
\hline
\end{tabular}
\end{table}

Note that $X_2$ is not compact, but we may impose the additional constraint $\|x(t)\|_\infty \leq N$ on $X_2$ for some large $N$ \cite[Section~5.1]{lasserre2008nonlinear}. 
However, this additional constraint is not enforced in the numerical implementations.
We first consider the following minimum time problem: drive the system to the point $(0,0)$ beginning from $x_0 = (0.3,1) \in X_1$ in minimum time.
We assume the minimum time needed is less than 5 and the problem is set up according to Table \cref{tab:DI:minT}.

\begin{table}[!h]
\caption{Optimal control problem setup of the hybridized double integrator}
\label{tab:DI:minT}
\centering
{\begin{tabular}{c|c|c}
\hline
    & $i=1$ & $i=2$ \\
\hline
    $h_i$ & 1 & 1\\
\hline
    $H_i$ & 0 & 0\\
\hline
    $x_0$ & $(0.3,1) \in X_1$ & N/A\\
\hline
    $X_{T_i}$ & $\emptyset$ & $\{(0,0)\} \subset X_2$\\
\hline
    $T_0$ & \multicolumn{2}{c}{5}\\
\hline
\end{tabular}}
\end{table}

\Later{Should put LQR setup after Fig.1}
\begin{table}[!ht]
\caption{Problem setup of the hybridized double integrator LQR problem}
\label{tab:DI:LQR}
\centering
\begin{tabular}{c|c|c}
\hline
    & $i=1$ & $i=2$ \\
\hline
    $h_i$ & $x_1^2 + x_2^2 + 20 \cdot u^2$ & $x_1^2 + x_2^2 + 20 \cdot u^2$\\
\hline
    $H_i$ & 0 & 0\\
\hline
    $x_0$ & $(1,1) \in X_1$ & N/A \\
\hline
    $X_{T_i}$ & $\{(x_1, x_2) \in \R^2 \mid x_1^2 + x_2^2 \leq 0.3\} = X_1$ & $\{(x_1, x_2) \in \R^2 \mid x_1^2 + x_2^2 \geq 0.3\} = X_2$\\
\hline
    $T$ & \multicolumn{2}{c}{5 or 15}\\
\hline
\end{tabular}
\end{table}

\begin{figure}[!h]
\begin{minipage}{0.9\columnwidth}
\centering
	\subfloat[Control action \label{fig:DI_control}]{\includegraphics[trim={1cm 6.5cm 2cm 7cm},clip,width=0.48\textwidth]{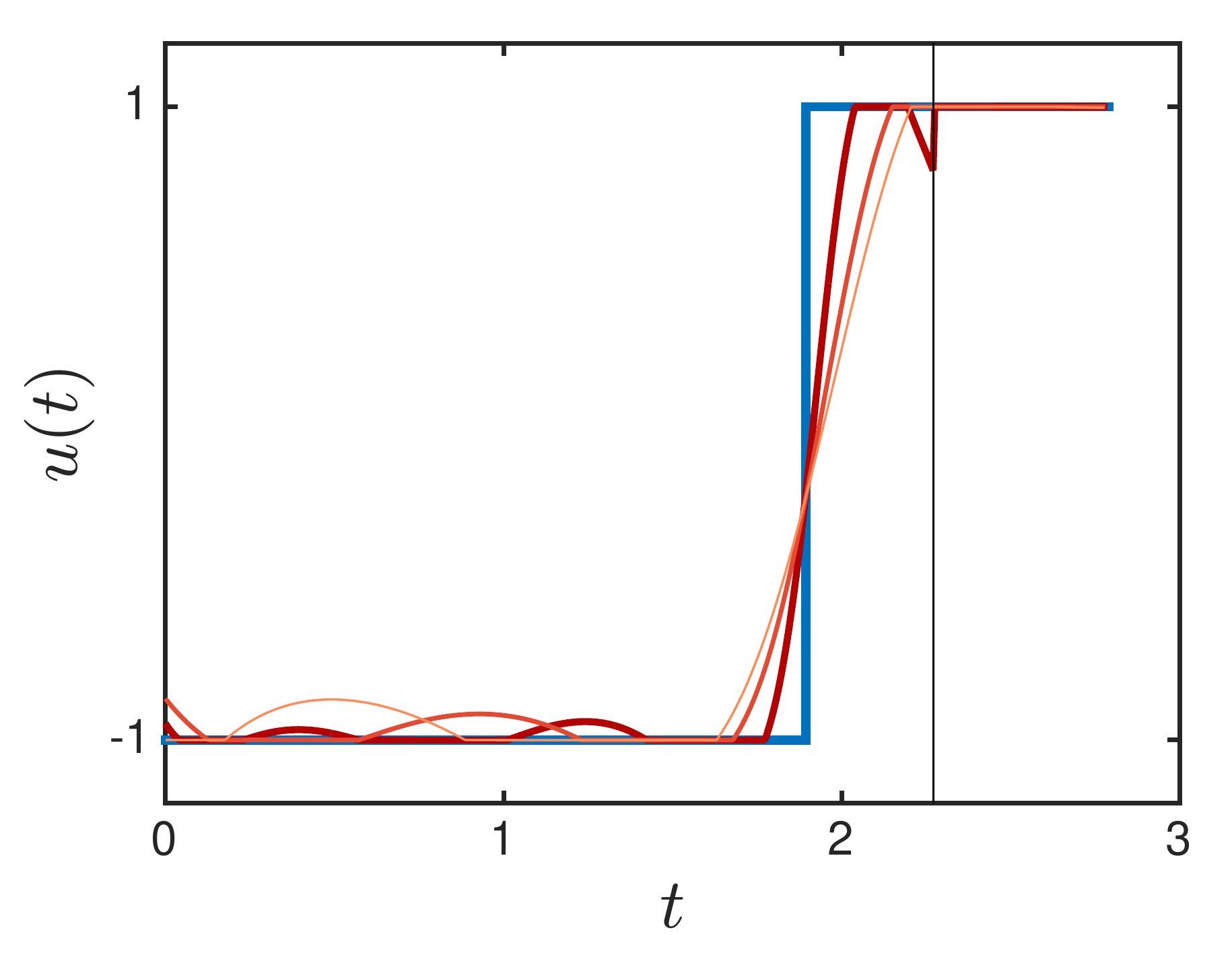}} 
	\hspace{0.4cm}
    \subfloat[Trajectory\label{fig:DI_traj}]{\includegraphics[trim={1cm 6.5cm 2cm 7cm},clip,width=0.48\textwidth]{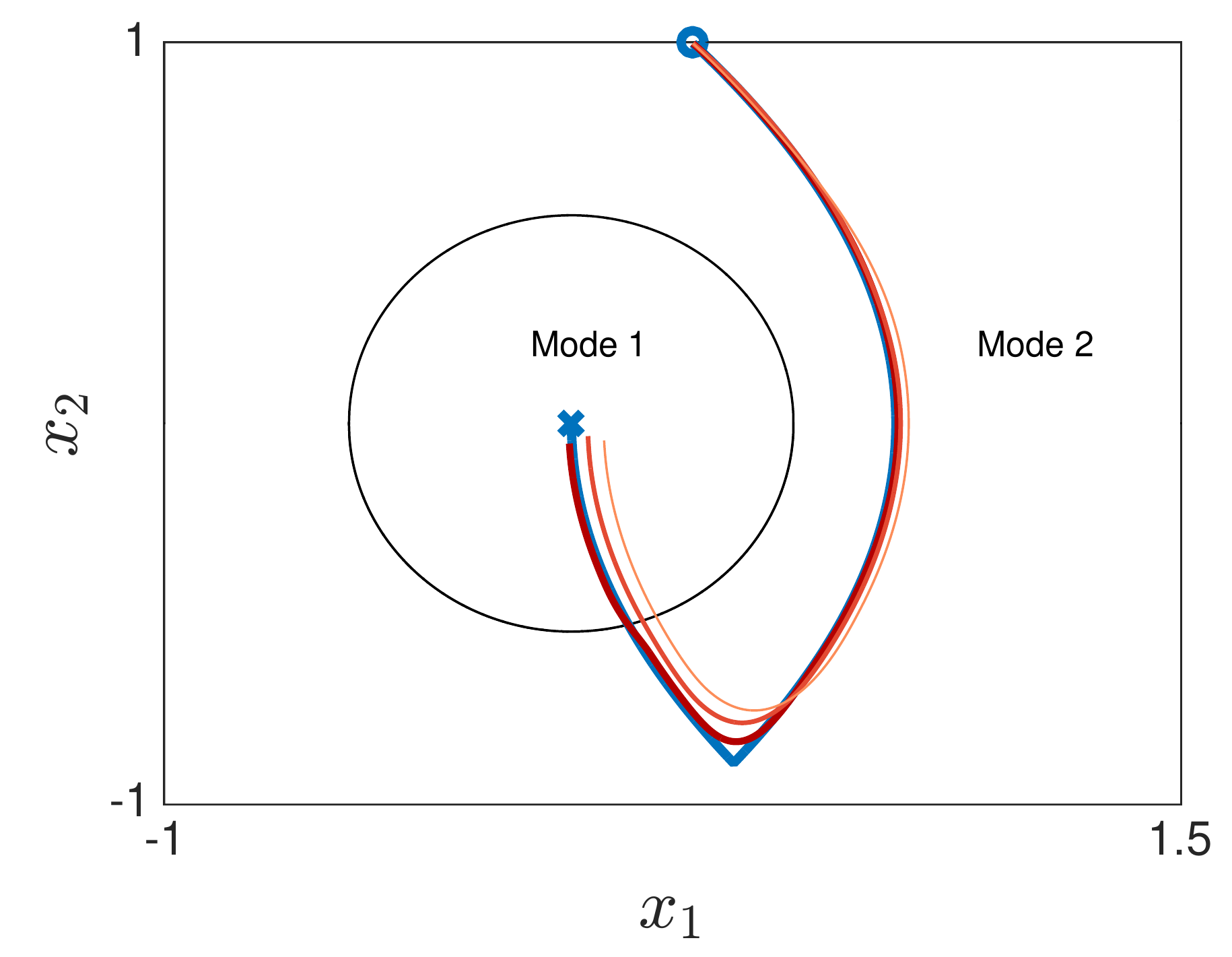}}
	\end{minipage}
\centering
\caption{An illustration of the performance of our algorithm on a free final time version of the hybridized double integrator problem.
The blue circle indicates the given initial point $x_0$, and the blue cross shows the target set.
The blue line is the analytically computed optimal control, while the red lines of various saturation correspond to control actions generated by our method.
When the simulated trajectory does not pass through $(0,0)$ perfectly, the simulation terminates when the closest point is reached.
As the saturation of the color in the illustration increases the corresponding degree of relaxation increases between $2k=6$ to $2k=8$ to $2k=12$.
Figure \cref{fig:DI_control} depicts the control action whereas Figure \cref{fig:DI_traj} illustrates the resultant trajectory when forward simulated through the system. The moment of transition from mode 2 to mode 1 is indicated by a vertical black solid line in Figure \cref{fig:DI_control}.}
\label{fig:DI}
\end{figure}

For this system, the optimal admissible pair is analytically computable, which is used as ground truth and compared to the result of our method with degrees of relaxation $k=6$, $k=8$, and $k = 12$ in Figure \cref{fig:DI}.
The polynomial control law is saturated so that its value is in $U$ for all time.
The cost and computation time are also compared in Table \cref{tab:DI:res}.

\begin{figure}[!h]
\begin{minipage}{0.99\columnwidth}
\centering
    \subfloat[Trajectory, $T=5$\label{fig:DI_H_LQR_traj_T5}]{\includegraphics[trim={1cm 6.5cm 2cm 7cm},clip,width=0.45\textwidth]{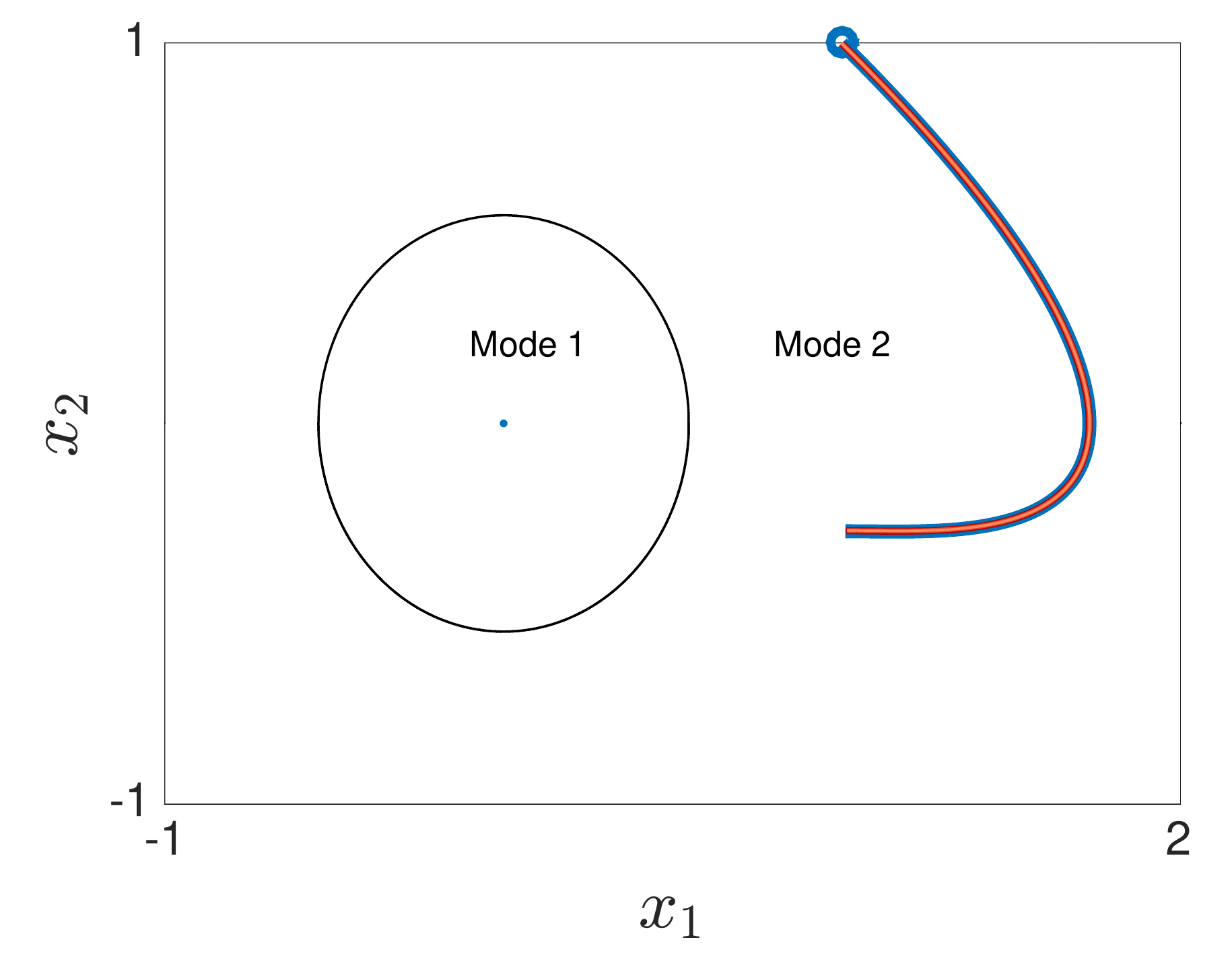}} 
    \hspace{0.4cm}
    \subfloat[Trajectory, $T=15$ \label{fig:DI_H_LQR_traj_T15}]{\includegraphics[trim={1cm 6.5cm 2cm 7cm},clip,width=0.45\textwidth]{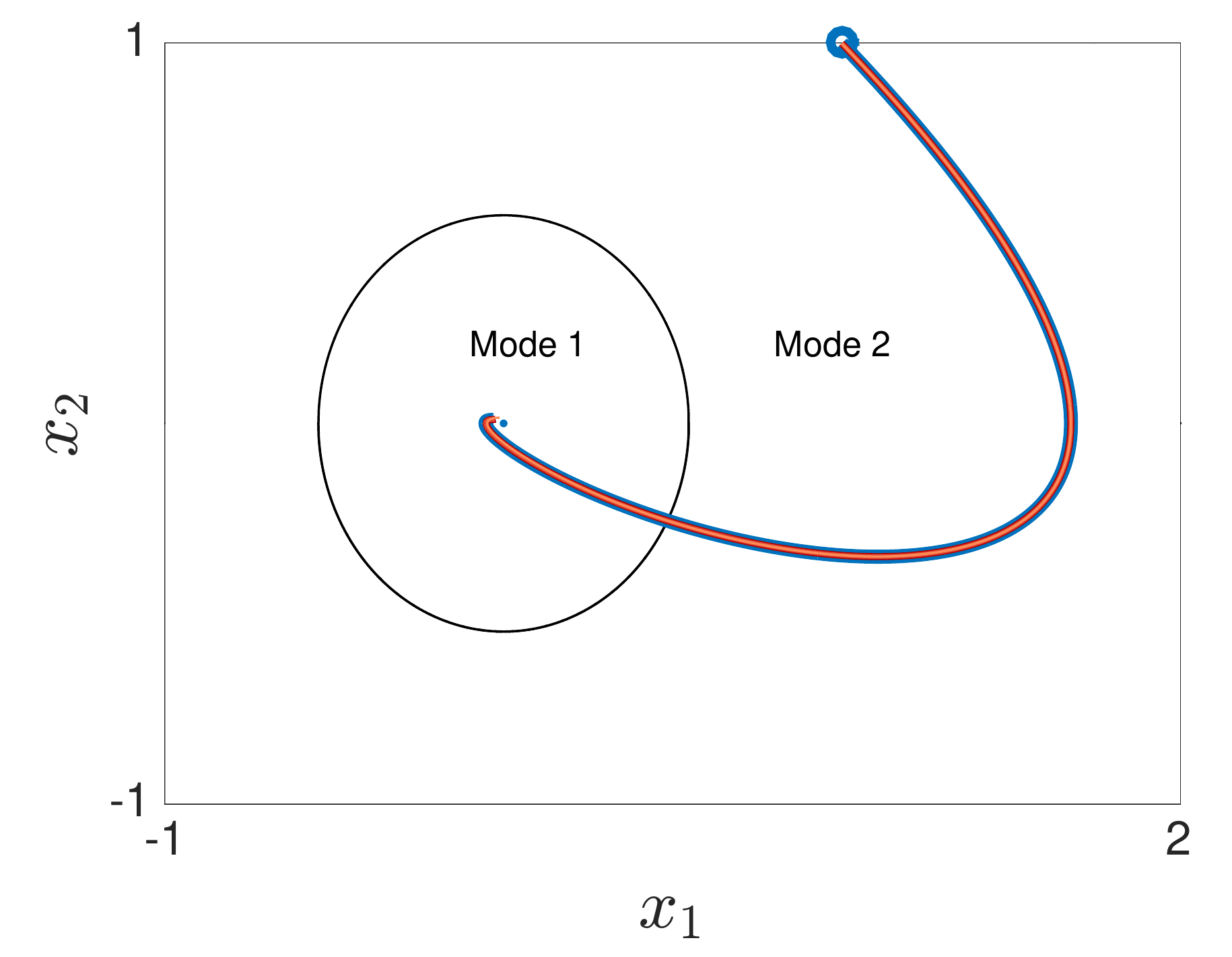}}
	\end{minipage}
\caption{An illustration of the performance of our algorithm on LQR version of the hybridized double integrator problem.
The blue circles indicate the given initial point $x_0$, and the blue dots show the point $(0,0)$.
The blue lines are the analytically computed optimal control, while the red lines of various saturation correspond to control actions generated by our method.
As the saturation of the color increases the corresponding degree of relaxation increases between $2k=6$ to $2k=8$ to $2k=12$.
Figure \cref{fig:DI_H_LQR_traj_T5} shows the trajectories when $T=5$, and Figure \cref{fig:DI_H_LQR_traj_T15} shows the trajectories when $T=15$.
}
\label{fig:DI_LQR}
\end{figure}

Next, we consider an Linear Quadratic Regulator (LQR) problem on the same hybridized double integrator system, where the goal is to drive the system state towards $(0,0)$ while keeping the control action small for all time $t \in [0,T]$.
The problem is set up according to Table \cref{tab:DI:LQR}.
To further illustrate we are able to handle different number of modes visited, two cases where $T=5$ and $T=15$ are considered.
For comparison, the LQR problem is also solved by a standard finite-horizon LQR solver in the non-hybrid case, which we refer to as the ground truth.
The results are compared in Figure \cref{fig:DI_LQR} and Table \cref{tab:DI:res} with degrees of relaxation $2k=6$, $2k=8$, and $2k = 12$.


\begin{table}[htbp]
\caption{Results of the hybridized double integrator examples}
\centering
\label{tab:DI:res}
\begin{minipage}{0.99\textwidth}
\begin{tabular}{c|c|c|C{3cm}|C{3cm}}
\hline
    \multicolumn{2}{c|}{} & Computation time & Cost returned from optimization & Cost returned from simulation \\
\hline
    \multirow{4}{3cm}{Minimum time problem with $T_0=5$}
        & $2k=6$ & 3.1075[s] & 2.7781 & 2.7780 \footnote{\label{footnote1}Trajectory does not reach target set perfectly. The simulation terminates when the closest point is reached} \\ \cline{2-5}
        & $2k=8$ & 10.0187[s] & 2.7847 & 2.7845 \footref{footnote1} \\ \cline{2-5}
        & $2k=12$ & 170.9319[s] & 2.7868 & 2.7865 \footref{footnote1} \\  \cline{2-5}
        & Ground truth & N/A & 2.7889 & N/A\\
\hline
    \multirow{4}{3cm}{LQR problem with $T=5$}
        & $2k=6$ & 2.2299[s] & 24.9496 & 24.9906 \\ \cline{2-5}
        & $2k=8$ & 8.1412[s] & 24.9496 & 24.9906 \\ \cline{2-5}
        & $2k=12$ & 198.2826[s] & 24.9502 & 24.9906 \\  \cline{2-5}
        & Ground truth & N/A & 24.9503 & N/A\\
\hline
    \multirow{4}{3cm}{LQR problem with $T=15$}
        & $2k=6$ & 2.1965[s] & 26.1993 & 26.3428 \\ \cline{2-5}
        & $2k=8$ & 7.7989[s] & 26.1993 & 26.3438 \\ \cline{2-5}
        & $2k=12$ & 168.5383[s] & 26.1996 & 26.3435 \\  \cline{2-5}
        & Ground truth & N/A & 26.2033 & N/A\\
\hline
\end{tabular}
\end{minipage}
\end{table}

\subsection{Dubins Car Model with Shortcut Path}
The next example shows our algorithm can work with different dimensions in each mode, and is capable of choosing the best transition sequence.
Consider a 2-mode hybridized Dubins Car system with identity reset map. 
We now add another 1-dimensional mode to the system, and connect it with the other two modes by defining transitions. 
The vector fields, guards, and reset maps are defined in Table \cref{tab:shortcut:vf} and Table \cref{tab:shortcut:guard}.
In mode 1 and mode 2, the control is $u = (v,\omega)$; In mode 2, the control is $u = v$.
Although the dynamics in mode 1 and mode 2 are not polynomials, they are approximated by 2nd-order Taylor expansion around $x = (0,0,0)$ in the numerical implementation.
We are interested in solving the minimum time problem, where the trajectory starts at $x_0 = (-0.8,0.8,0)$ in mode 1, and ends at $x,y$-position $(0.8,-0.8)$ in mode 2. 
The optimal control problem is defined in Table \cref{tab:shortcut}.


\begin{table}[h]
\caption{Vector fields and domains of the Dubins car model with shortcut path}
\label{tab:shortcut:vf}
\centering
\begin{tabular}{c|c|c|c}
\hline
    Mode & $i=1$ & $i=2$ & $i=3$ \\
\hline
    Dynamics & $\dot{x}(t) = 
    \begin{bmatrix}
        0 \\ 0 \\ 0
    \end{bmatrix} +
    \begin{bmatrix}
        \cos(x_3(t)) & 0 \\ \sin(x_3(t)) & 0 \\ 0 & 1
    \end{bmatrix} u$ & 
    $\dot{x}(t) = 
    \begin{bmatrix}
        0 \\ 0 \\ 0
    \end{bmatrix} +
    \begin{bmatrix}
        \cos(x_3(t)) & 0 \\ \sin(x_3(t)) & 0 \\ 0 & 1
    \end{bmatrix} u$ &
    $\dot{x}(t) = 0 + (-1) \cdot u$\\
\hline
    $X_i$ & $[-1,1] \times [0,1] \times [-\pi,\pi] \subset \R^3$ & $[-1,1]\times[-1,0]\times [-\pi,\pi] \subset \R^3$ & $[-1,1] \subset \R$\\
\hline
    $U_i$ & $[0,1] \times [-3,3]$ & $[0,1] \times [-3,3]$ & $[0,2]$\\
\hline
\end{tabular}
\end{table}

\begin{table}[h]
\caption{Guards and reset maps of the Dubins car model with shortcut path.
The rows are modes in which a transition originates, and the columns are modes to which the transition goes.
}
\label{tab:shortcut:guard}
\begin{tabular}{c|c|c|c}
\hline
    & Mode 1 & Mode 2 & Mode 3 \\
\hline
    Mode 1 &
    $\begin{aligned}
    &S_{(1,2)} = [-1,1] \times \{0\} \times [-\pi, \pi]\\
    &R_{(1,2)}(x) = x, \quad \forall x \in S_{(2,1)}
    \end{aligned}$ & N/A & 
    $\begin{aligned}
    &S_{(1,3)} = [-1,1] \times \{1\} \times [-\pi, \pi]\\
    &R_{(1,3)}(x) = 1, \quad \forall x \in S_{(1,3)}
    \end{aligned}$\\
\hline
    Mode 2 & N/A & N/A & N/A \\
\hline
    Mode 3 & N/A & 
    $\begin{aligned}
    &S_{(3,2)} = \{-1\}\\
    &R_{(3,2)}(x) = (0.6,-0.8,0), \\
    &\forall x \in S_{(3,2)}
    \end{aligned}$ & N/A\\
\hline
\end{tabular}
\end{table}

\begin{table}[!h]
\caption{Problem setup of Dubins car model with shortcut path}
\label{tab:shortcut}
\centering
\begin{tabular}{c|c|c|c}
\hline
    Mode & $i=1$ & $i=2$ & $i=3$ \\
\hline
    $h_i$ & 1 & 1 & 1\\
\hline
    $H_i$ & 0 & 0 & 0\\
\hline
    $x_0$ & $(-0.8,0.8,0) \subset X_1$ & N/A & N/A\\
\hline
    $X_{T_i}$ & N/A & $\{0.8\}\times\{-0.8\}\times[-\pi,\pi] \subset X_2$ & N/A\\
\hline
    $T$ & \multicolumn{3}{c}{3}\\
\hline
\end{tabular}
\end{table}

\begin{table}[!h]
\caption{Results of the Dubins car example with shortcut path}
\centering
\label{tab:shortcut:res}
\begin{tabular}{C{3cm}|c|C{3cm}|C{3cm}}
\hline
     & Computation time & Cost returned from optimization & Cost returned from simulation \\
\hline
        $2k=6$ & 83.0224[s] & 1.5641 & 1.5739 \\
        \hline
        $2k=8$ & $1.2115 \times 10^3$[s] & 1.5647 & 1.5679 \\
        \hline
        $2k=10$ & $1.3206 \times 10^4$[s] & 1.5648 & 1.5703 \\
        \hline
        Ground truth & N/A & 1.5651 & N/A\\
\hline
\end{tabular}
\end{table}

\begin{figure}[htbp]
\begin{minipage}{0.3\textwidth}
    \centering
    \subfloat[Control action: $v$ \label{fig:shortcut_control1}]{\includegraphics[trim={1cm 6.9cm 2cm 7.5cm},clip,width=0.99\textwidth]{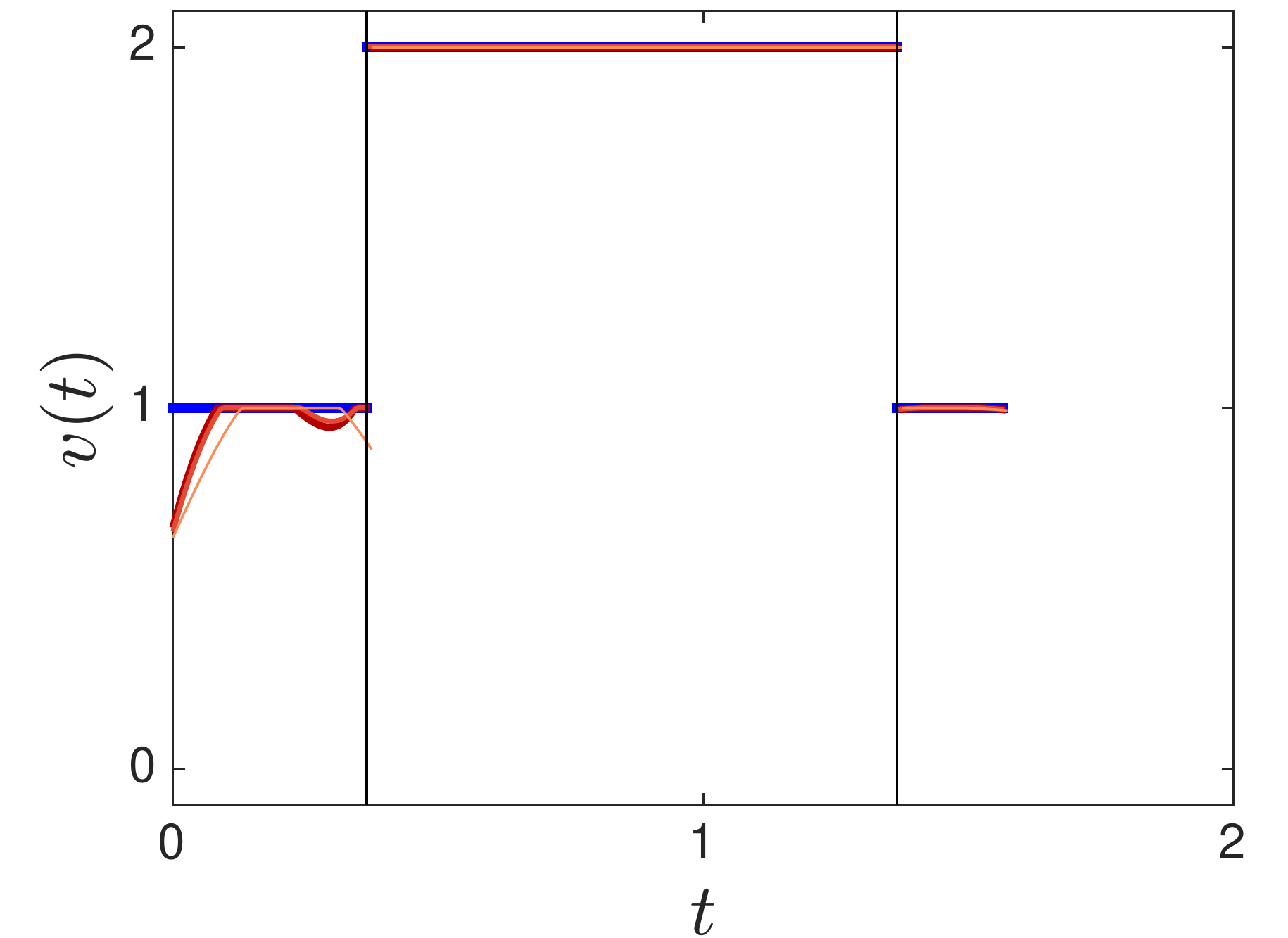}}\\
    \subfloat[Control action: $\omega$ \label{fig:shortcut_control2}]{\includegraphics[trim={1cm 6.9cm 2cm 7.5cm},clip,width=0.99\textwidth]{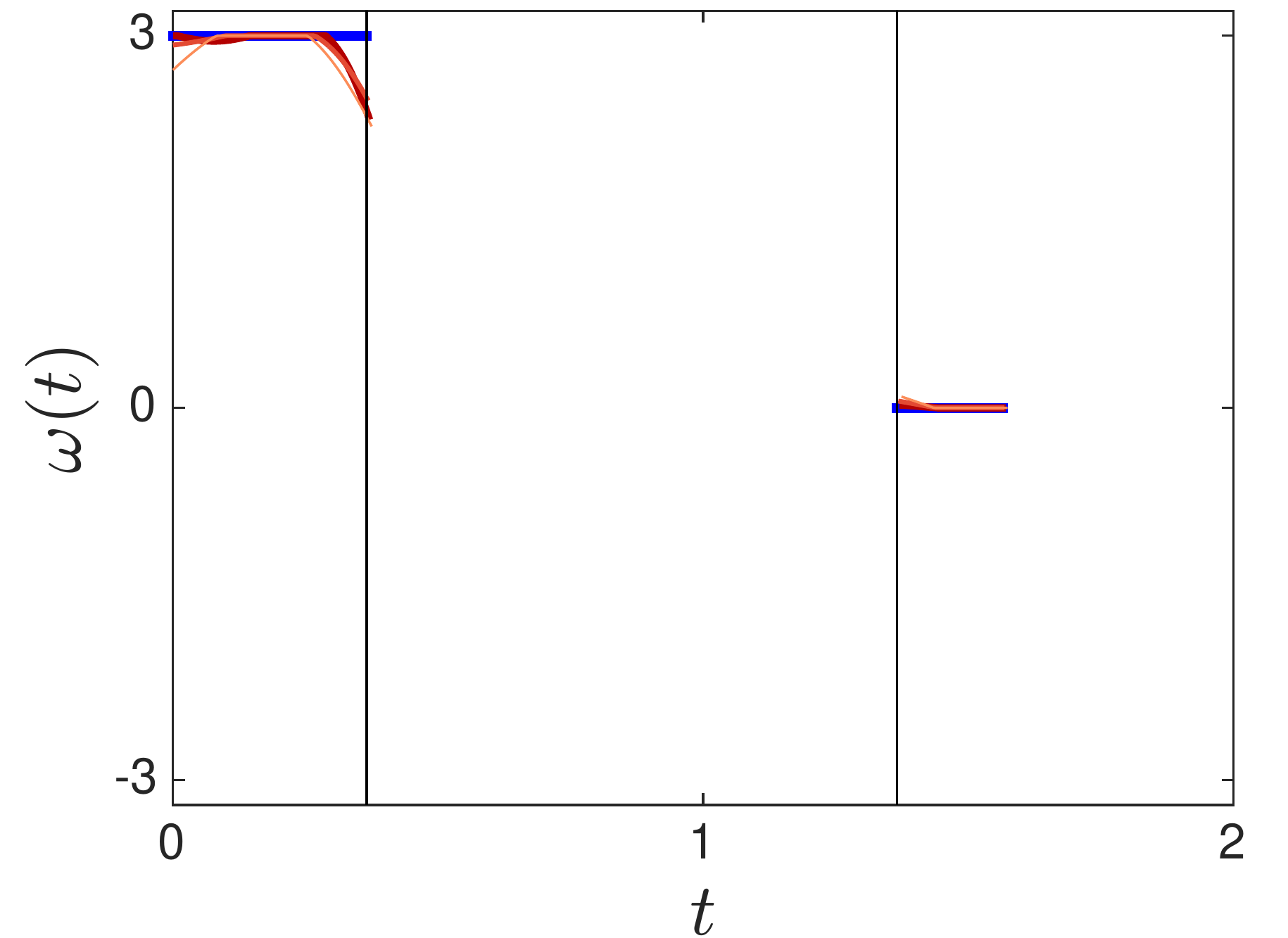}}
\end{minipage}
\begin{minipage}{0.6\textwidth}
\centering
    \subfloat[Trajectory\label{fig:shortcut_traj}]{\includegraphics[trim={1cm 6.5cm 1cm 7cm},clip,width=0.99\textwidth]{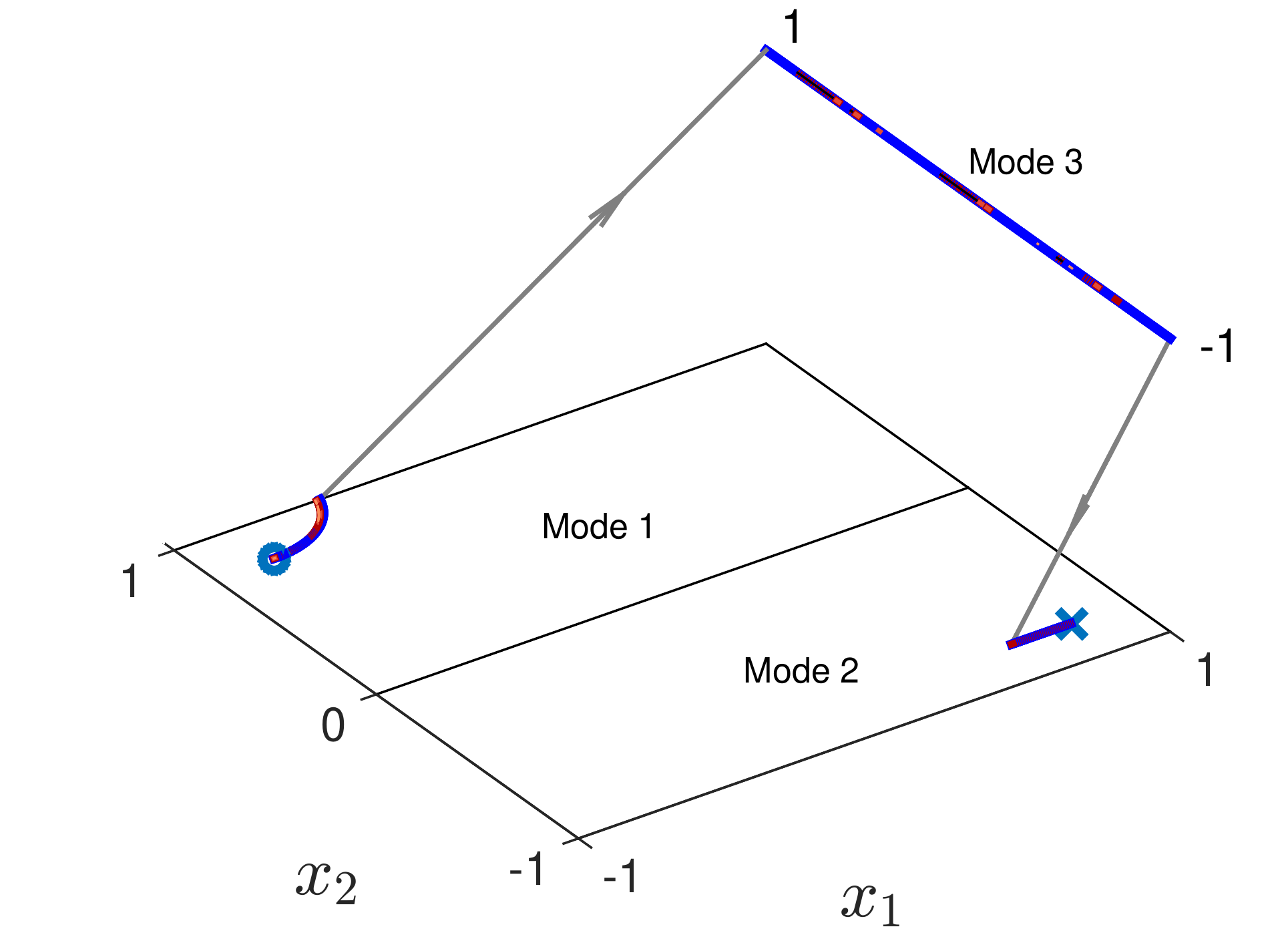}}
	\end{minipage}
\centering
\caption{An illustration of the performance of our algorithm on a minimum time problem of the Dubins car model with shortcut path.
The blue circle indicates the given initial point $x_0$, and the blue cross shows the target set.
The blue solid line is the analytically computed optimal control, and the red lines of various saturations are controls generated by our method.
As the saturation increases the corresponding degree of relaxation increases between $2k=6$ to $2k=8$ to $2k=10$.
Figure \cref{fig:shortcut_control1}, Figure \cref{fig:shortcut_control2} depict the control actions whereas Figure \cref{fig:shortcut_traj} illustrates the corresponding trajectory obtained by forward simulating through the system.
The moment of transitions are indicated by vertical black solid lines in Figure\cref{fig:shortcut_control1} and Figure \cref{fig:shortcut_control2}.}
\label{fig:shortcut}
\end{figure}

Notice the transition sequences ``1-2'' and ``1-3-2'' are both feasible in this instance according to our guard definition, but direct calculation shows that we may arrive at the target point in less time by taking the ``shortcut path'' in mode $3$.
This problem is solved using our algorithm with degrees of relaxation $2k=6$, $2k=8$, and $2k=10$.
As comparison, we treat the analytically computed optimal control as ground truth, and the results are compared in Figure \cref{fig:shortcut} and Table \cref{tab:shortcut:res}.
In this example our algorithm is able to pick the transition sequence ``1-3-2'' and find a tight approximation to the true optimal solution.

\subsection{SLIP Model}
The Spring-Loaded Inverted Pendulum (SLIP) is a classical model that describes the center-of-mass dynamics of running animals and robots, and has been extensively used as locomotion template to perform control law synthesis on legged robots \cite{holmes2006dynamics}. 
Despite its simplicity, an analytical solution to the SLIP dynamics does not exist. 
We may simulate the system numerically, but the optimal control problem is still difficult to solve if the sequence of transition is not known beforehand.

\begin{figure}[!b]
\begin{minipage}{0.9\textwidth}
    \centering
    \hspace{1.5cm}
    \subfloat[SLIP model\label{fig:SLIP:config:1}]{
    \begin{tikzpicture}
    	\tikzstyle{ann} = [fill=white,font=\large,inner sep=1pt]
    	\node[anchor=south west,inner sep=0] (image) at (0,0) {{\includegraphics[trim={4cm 8.5cm 3.5cm 9cm},clip=true,width=0.42\textwidth]{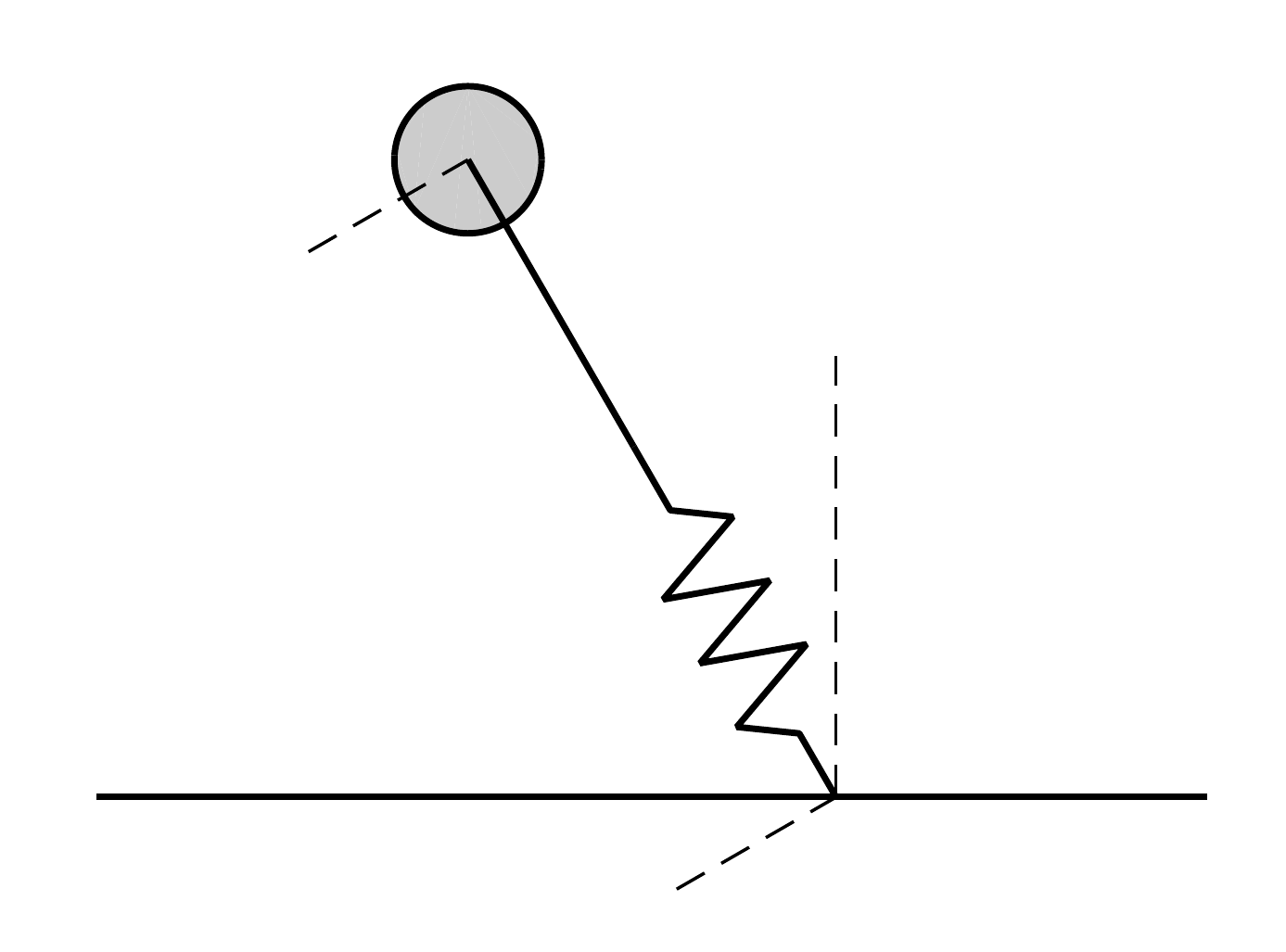}}};
        \begin{scope}[x={(image.south east)},y={(image.north west)}]
        \draw [-{Latex[scale=1.5]}] (0.1,0.165)--(0.26,0.165);
        \draw [-{Latex[scale=1.5]}] (0.1,0.165)--(0.1,0.35);
        \node[ann] at (0.1,0.12) {$o$};
        \node[ann] at (0.13,0.35) {$b$};
        \node[ann] at (0.27,0.12) {$a$};
        \draw [{Latex[scale=1.5]}-{Latex[scale=1.5]}] (0.26,0.75)--(0.55,0.08);
        \node[ann] at (0.405,0.405) {$l$};
        \draw [-{Latex[scale=1.5]}] (0.65,0.55) arc (90:135:0.2);
        \node[ann] at (0.59,0.54) {$\theta$};
        \end{scope}
    \end{tikzpicture}}
    \hspace{0cm}
    \subfloat[Active SLIP model\label{fig:SLIP:config:2}]{
    \begin{tikzpicture}
    	\tikzstyle{ann} = [fill=white,font=\large,inner sep=1pt]
    	\node[anchor=south west,inner sep=0] (image) at (0,0) {{\includegraphics[trim={4cm 8.5cm 3.5cm 9cm},clip=true,width=0.42\textwidth]{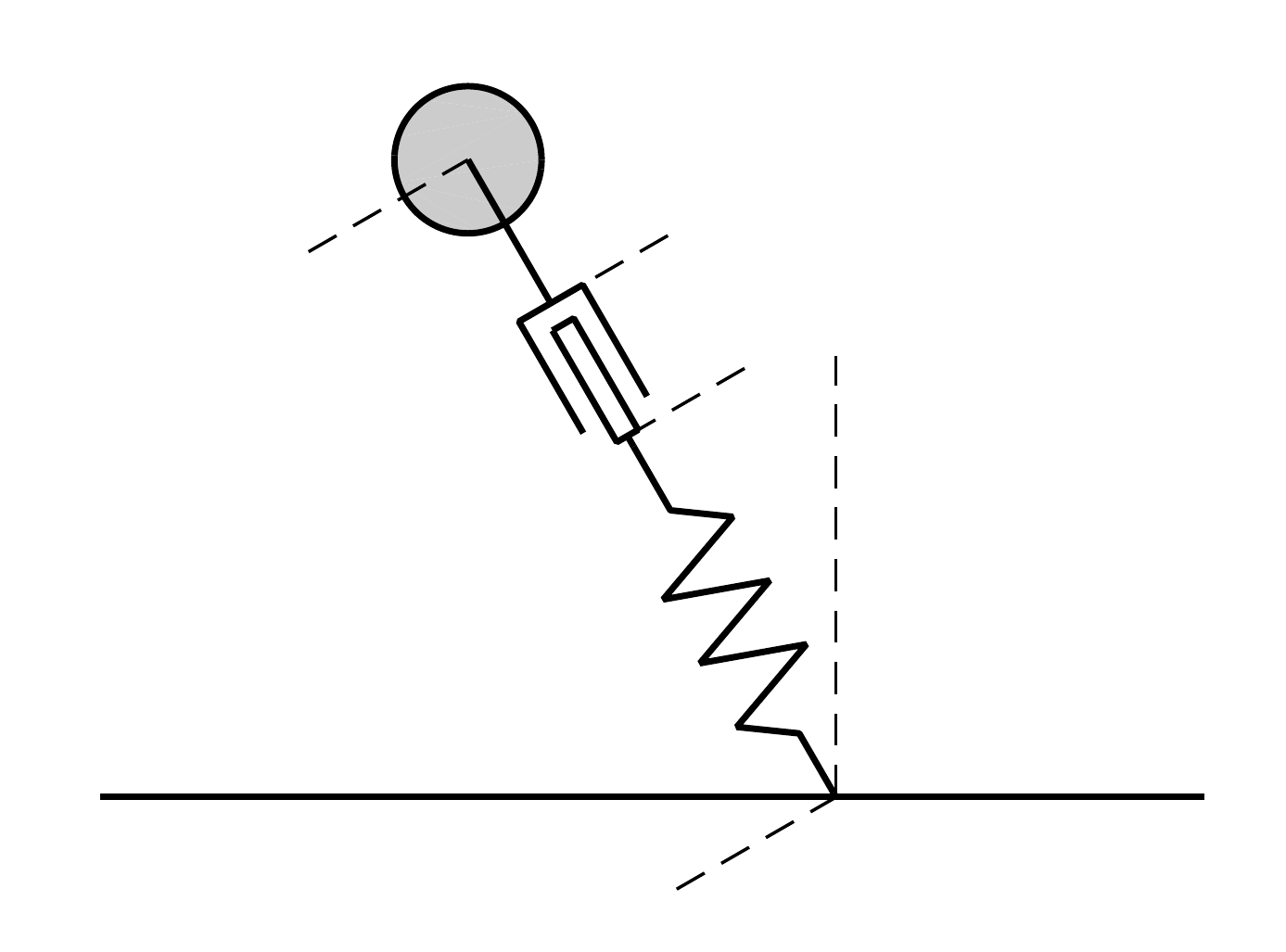}}};
        \begin{scope}[x={(image.south east)},y={(image.north west)}]
        \draw [-{Latex[scale=1.5]}] (0.1,0.165)--(0.26,0.165);
        \draw [-{Latex[scale=1.5]}] (0.1,0.165)--(0.1,0.35);
        \node[ann] at (0.1,0.12) {$o$};
        \node[ann] at (0.13,0.35) {$b$};
        \node[ann] at (0.27,0.12) {$a$};
        \draw [{Latex[scale=1.5]}-{Latex[scale=1.5]}] (0.26,0.75)--(0.55,0.08);
        \node[ann] at (0.405,0.405) {$l$};
        \draw [-{Latex[scale=1.5]}] (0.65,0.55) arc (90:135:0.2);
        \node[ann] at (0.59,0.54) {$\theta$};
        \draw [{Latex[scale=1.5]}-{Latex[scale=1.5]}] (0.56,0.6)--(0.5,0.74);
        \node[ann] at (0.56,0.69) {$u$};
        \end{scope}
    \end{tikzpicture}}
\end{minipage}
\caption{Slip model and system variable definition}
\label{fig:SLIP:config}
\end{figure}

\begin{figure}[!htb]
\centering
\begin{tikzpicture}
	\tikzstyle{ann} = [fill=white,font=\large,inner sep=1pt]
	\tikzstyle{rot} = [fill=white,font=\large,rotate=90]
    \node[anchor=south west,inner sep=0] (image) at (0,0) {{\includegraphics[trim={3cm 10cm 2.5cm 10cm},clip=true,width=0.75\textwidth]{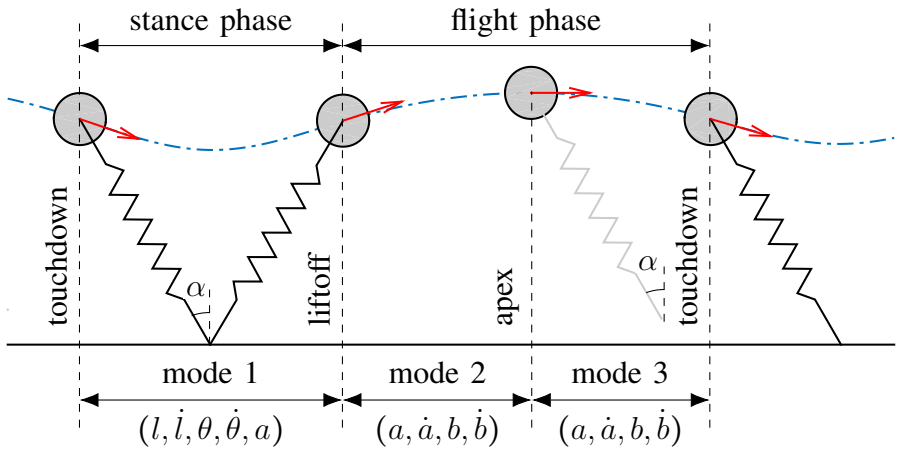}}};
    \begin{scope}[x={(image.south east)},y={(image.north west)}]
    \draw [dashed] (0.108,0.07)--(0.108,0.95);
    \draw [dashed] (0.39,0.07)--(0.39,0.95);
    \draw [dashed] (0.593,0.07)--(0.593,0.8);
    \draw [dashed] (0.784,0.07)--(0.784,0.95);
	\draw[{Latex[scale=1.5]}-{Latex[scale=1.5]}] (0.108,0.13) -- (0.39,0.13);
    \draw[{Latex[scale=1.5]}-{Latex[scale=1.5]}] (0.39,0.13) -- (0.593,0.13);
    \draw[{Latex[scale=1.5]}-{Latex[scale=1.5]}] (0.593,0.13) -- (0.784,0.13);
    \draw[{Latex[scale=1.5]}-{Latex[scale=1.5]}] (0.108,0.9) -- (0.39,0.9);
    \draw[{Latex[scale=1.5]}-{Latex[scale=1.5]}] (0.39,0.9) -- (0.784,0.9);
    \node[ann] at (0.249,0.95) {stance phase};
    \node[ann] at (0.587,0.95) {flight phase};
    
    \node[ann] at (0.249,0.19) {mode 1};
    \node[ann] at (0.249,0.07) {$(l,\dot{l},\theta, \dot{\theta}, a)$};
    \node[ann] at (0.492,0.19) {mode 2};
    \node[ann] at (0.492,0.07) {$(a,\dot{a},b,\dot{b})$};
    \node[ann] at (0.688,0.19) {mode 3};
    \node[ann] at (0.688,0.07) {$(a,\dot{a},b,\dot{b})$};
    \draw [dashed] (0.2475,0.25)--(0.2475,0.37);
    \draw (0.2475,0.32) arc (90:130:0.03);
    \node[ann] at (0.233,0.37) {$\alpha$};
    \draw [dashed] (0.735,0.3)--(0.735,0.42);
    \draw (0.735,0.37) arc (90:130:0.03);
    \node[ann] at (0.718,0.427) {$\alpha$};
    \node[rot,anchor=west] at (0.083,0.27) {touchdown};
    \node[rot,anchor=west] at (0.365,0.27) {liftoff};
    \node[rot,anchor=west] at (0.568,0.27) {apex};
    \node[rot,anchor=west] at (0.759,0.27) {touchdown};
    \end{scope}
\end{tikzpicture}
\caption{SLIP locomotion phases and hybrid system modes}
\label{fig:SLIP:model}
\end{figure}

As is shown in Figure \cref{fig:SLIP:config:1}, the SLIP is a mass-spring physical system, modeled as a point mass, $M$, and a mass-less spring leg with stiffness $k$ and length $l$. 
The dynamics of SLIP consist of two phases: stance phase and flight phase. 
The stance phase starts when the leg comes into contact with the ground with downward velocity, which we call the \emph{touchdown} event, and ends when the leg extends to full length and leaves the ground, which we call the \emph{liftoff} event. 
During the stance phase, the inverted pendulum swings forward around the leg-ground contact point, while the spring contracts due to mass momentum and gravitational force. 
During flight phase, SLIP follows free fall motion where the only external force is the gravity. 
We also assume the leg angle is reset to some fixed value $\alpha$ instantaneously once the SLIP enters flight phase, so that the leg angle at the moment of touchdown stays the same. 
Furthermore, we define the \emph{apex} event to be when the body reaches its maximum height with zero vertical velocity. 
The touchdown, liftoff, and apex events are illustrated in Figure \cref{fig:SLIP:model}.

In the context of this paper, we are interested in the active SLIP model (Figure \cref{fig:SLIP:config:2}), where a mass-less actuator is added to the SLIP leg. 
During stance phase, the actuator may extend from its nominal position within some range, while during flight phase, the actuator has no effect on the system. 
The active SLIP can be modeled as a hybrid system with 3 modes, where the liftoff, apex, and touchdown events define the transitions between them, as shown in Figure \cref{fig:SLIP:model}.

\begin{table}[!h]
\caption{State variables of the active SLIP mode}
\label{tab:SLIP:vars}
\centering
\begin{tabular}{c|l|c|l}
\hline
    $l$ & leg length & $a$ & horizontal displacement \\
    \hline
    $\dot{l}$ & time derivative of $l$ & $\dot{a}$ & time derivative of $a$ \\
    \hline
    $\theta$ & leg angle & $b$ & vertical displacement \\
    \hline
    $\dot{\theta}$ & time derivative of $\theta$ & $\dot{b}$ & time derivative of $b$ \\
\hline
\end{tabular}
\end{table}

\begin{table}[!h]
\caption{Physical parameters of the active SLIP model}
\label{tab:SLIP:params}
\centering
\begin{tabular}{c|l|c}
\hline
    & Explanation & Value\\
    \hline
    $M$ & mass & 1 \\
    \hline
    $k$ & spring constant & 6 \\
    \hline
    $g_0$ & gravitational acceleration & 0.2 \\
    \hline
    $l_0$ & nominal leg length & 0.2 \\
    \hline
    $\alpha$ & reset angle in flight phase & $\pi/6$\\
\hline
\end{tabular}
\end{table}

The behavior of such a system can be fully characterized using 8 variables defined in Table \cref{tab:SLIP:vars}. 
In mode 1, we define the system state to be $x = (l, \dot{l}, \theta, \dot{\theta}, a)$; In mode 2 and mode 3, we define the system state to be $x = (a, \dot{a}, b, \dot{b})$.
The physical parameters, dynamics, and transitions are defined in Table \cref{tab:SLIP:params}, Table \cref{tab:SLIP:vf}, and Table \cref{tab:SLIP:guard}. 
Again, we use 3rd-order Taylor expansion around $(l_0,0,0,0,0)$ to approximate the stance phase dynamics with polynomials.


\begin{table}[htbp]
\scriptsize
\caption{Vector fields and domains of each of the modes of the active SLIP model}
\label{tab:SLIP:vf}
\centering
{\renewcommand{\arraystretch}{0.7}
\begin{tabular}{c|c|C{2.8cm}|C{3.5cm}}
\hline
    Mode & $i=1$ & $i=2$ & $i=3$ \\
\hline
    Dynamics & $\dot{x}(t) = 
    \begin{bmatrix}
        x_2(t) \\
        -\frac{k}{M}(x_1(t)-l_0) - g_0\cos(x_3(t))\\
        x_4(t) \\
        - \frac{2x_2(t)x_4(t)}{x_1(t)} - \frac{g_0\sin(x_3(t))}{x_1(t)}\\
        -x_2(t)\sin(x_3(t)) - x_1(t)x_4(t)\cos(x_3(t))
    \end{bmatrix} +
    \begin{bmatrix}
        0 \\ \frac{k}{M} \\ 0 \\ 0 \\ 0
    \end{bmatrix} u$ & 
    $\dot{x}(t) = 
    \begin{bmatrix}
        x_2(t) \\ 0 \\ x_4(t) \\ -g_0
    \end{bmatrix}$ & 
    $\dot{x}(t) = 
    \begin{bmatrix}
        x_2(t) \\ 0 \\ x_4(t) \\ -g_0
    \end{bmatrix}$ \\
\hline
    $X_i$ & $[0.1,0.2]\times[-0.3,0.3]\times[-1,1]\times[-3,0]\times[-1,1] \subset \R^5$ & $[-1,1]\times[0,0.5]\times[0.15,0.5]\times[0,1] \subset \R^4$ & $[-1,1]\times[0,0.5]\times[l_0\cos(\alpha),0.5]\times[-1,0] \subset \R^4$ \\
\hline
    $U_i$ & $[0,0.1]$ & N/A & N/A\\
\hline
\end{tabular}
}
\end{table}

\begin{table}[htbp]
\scriptsize
\caption{Guards and reset maps of the active SLIP model.
The rows are modes in which a transition originates, and the columns are modes to which the transition goes.
}
\label{tab:SLIP:guard}
\centering
{\renewcommand{\arraystretch}{0.7}
\begin{tabular}{c|c|c|c}
\hline
    & Mode 1 & Mode 2 & Mode 3 \\
\hline
    Mode 1 & N/A &
    $\begin{aligned}
    &S_{(1,2)} = \{x \in X_1 \mid x_1 = l_0, x_2 \geq 0\}\\
    &R_{(1,2)}(x) = 
    \begin{bmatrix}
    x_5 \\ -x_2\sin(x_3)-l_0x_4\cos(x_3) \\ l_0\cos(x_3) \\ x_2\cos(x_3) - l_0x_4\sin(x_3))
    \end{bmatrix}, \\
    & \forall x \in S_{(1,2)}
    \end{aligned}$ & N/A\\
\hline
    Mode 2 & N/A & N/A &
    $\begin{aligned}
    &S_{(2,3)} = \{x \in X_2 \mid x_4 = 0 \}\\
    &R_{(2,3)}(x) = x, \quad \forall x \in S_{(2,3)}
    \end{aligned}$\\
\hline
    Mode 3 & 
    $\begin{aligned}
    &S_{(3,1)} = \{x \in X_3 \mid x_3 = l_0\cos(\alpha) \}\\
    &R_{(3,1)}(x) =
    \begin{bmatrix}
    l_0 \\ -x_2\sin(\alpha)+x_4\cos(\alpha) \\ \alpha \\ -\frac{x_2}{l_0}\cos(\alpha) - \frac{x_4}{l_0}\sin(\alpha))
    \end{bmatrix}, \\
    & \forall x \in S_{(3,1)}
    \end{aligned}$
    & N/A & N/A\\
\hline
\end{tabular}
}
\end{table}

We fix the initial condition, and consider the following two hybrid optimal control problems for the active SLIP:
In the first problem, we maximize the vertical displacement $b$ up to time $T=2.5$. In stance phase, the 1st-order Taylor approximation $b = l \cos(\theta) \approx l$ is used;
In the second problem, we define a constant-speed reference trajectory $a(t)=vt-1$ in the horizontal coordinate, then try to follow this trajectory with active SLIP up to time $T=4$.
The optimal control problems are defined in Table \cref{tab:SLIP}.
Note that these problems are defined such that the optimal transition sequences are different in each instance, and some modes are visited multiple times.

\begin{table}[htbp]
\caption{Problem setup of SLIP}
\label{tab:SLIP}
\centering
\begin{tabular}{c|c|c|c|c}
\hline
    & Mode & $i=1$ & $i=2$ & $i=3$\\
\hline
    \multirow{5}{4cm}{Maximizing vertical displacement}
        & $h_i$ & $-x_1$ & $-x_3$ & $-x_3$ \\ \cline{2-5}
        & $H_i$ & 0 & 0 & 0 \\ \cline{2-5}
        & $x_0$ & N/A & N/A & $(-1,0.3,0.2,0) \in X_3$ \\ \cline{2-5}
        & $X_{T_i}$ & $X_1$ & $X_2$ & $X_3$ \\ \cline{2-5}
        & $T$ & \multicolumn{3}{c}{2.5} \\
\hline
    \multirow{5}{4cm}{Tracking constant-speed trajectory $a(t)=vt-1$ with $v=0.1$}
        & $h_i$ & $(v \cdot t-1-x_5)^2$ & $(v \cdot t-1-x_1)^2$ & $(v \cdot t-1-x_1)^2$ \\ \cline{2-5}
        & $H_i$ & 0 & 0 & 0 \\ \cline{2-5}
        & $x_0$ & N/A & N/A & $(-1,0.3,0.2,0) \in X_3$ \\ \cline{2-5}
        & $X_{T_i}$ & $X_1$ & $X_2$ & $X_3$ \\ \cline{2-5}
        & $T$ & \multicolumn{3}{c}{4} \\
\hline
\end{tabular}
\end{table}

\begin{figure}[H]
\begin{minipage}{0.99\columnwidth}
\centering 
	\subfloat[Maximizing vertical displacement \label{fig:SLIP:b}]{\includegraphics[trim={7cm 7.5cm 7cm 7cm},clip,width=0.35\textwidth]{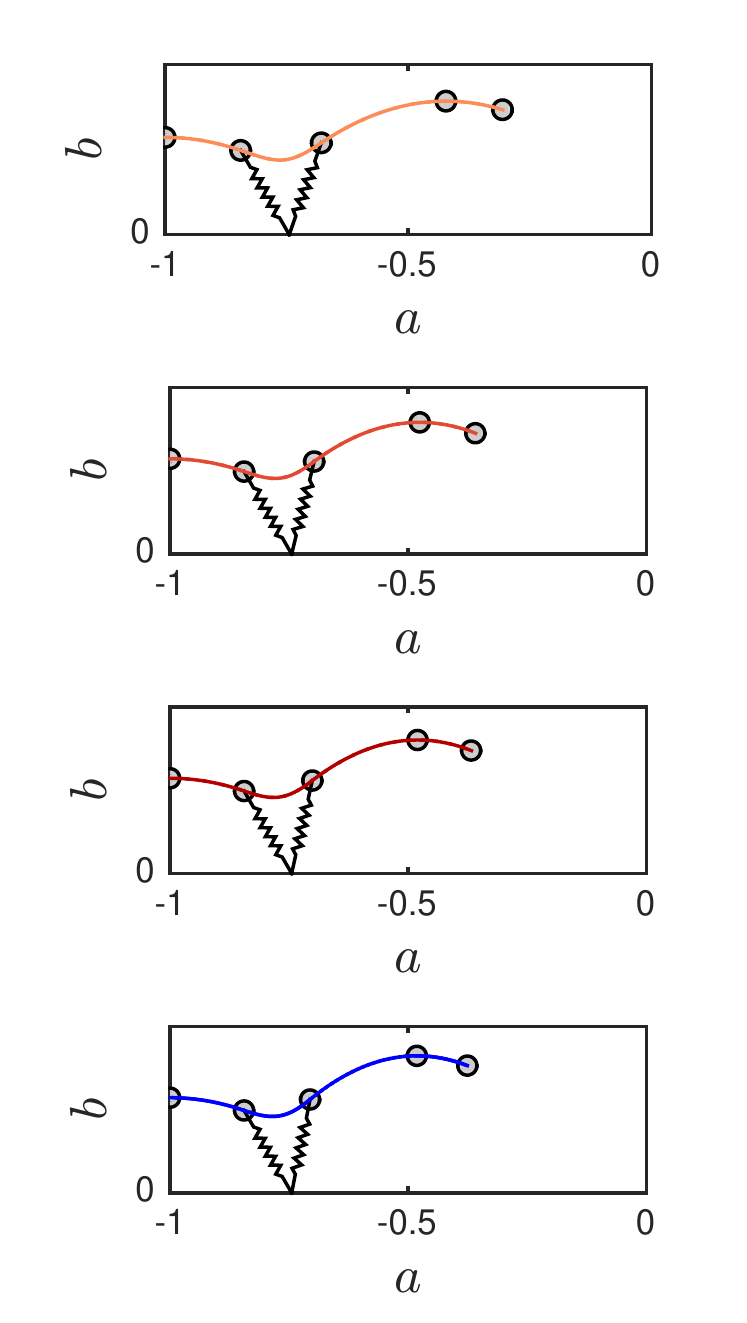}} \hspace{1cm}
	\subfloat[Tracking constant speed $v=0.1$ \label{fig:SLIP:V}]{\includegraphics[trim={7cm 7.5cm 7cm 7cm},clip,width=0.35\textwidth]{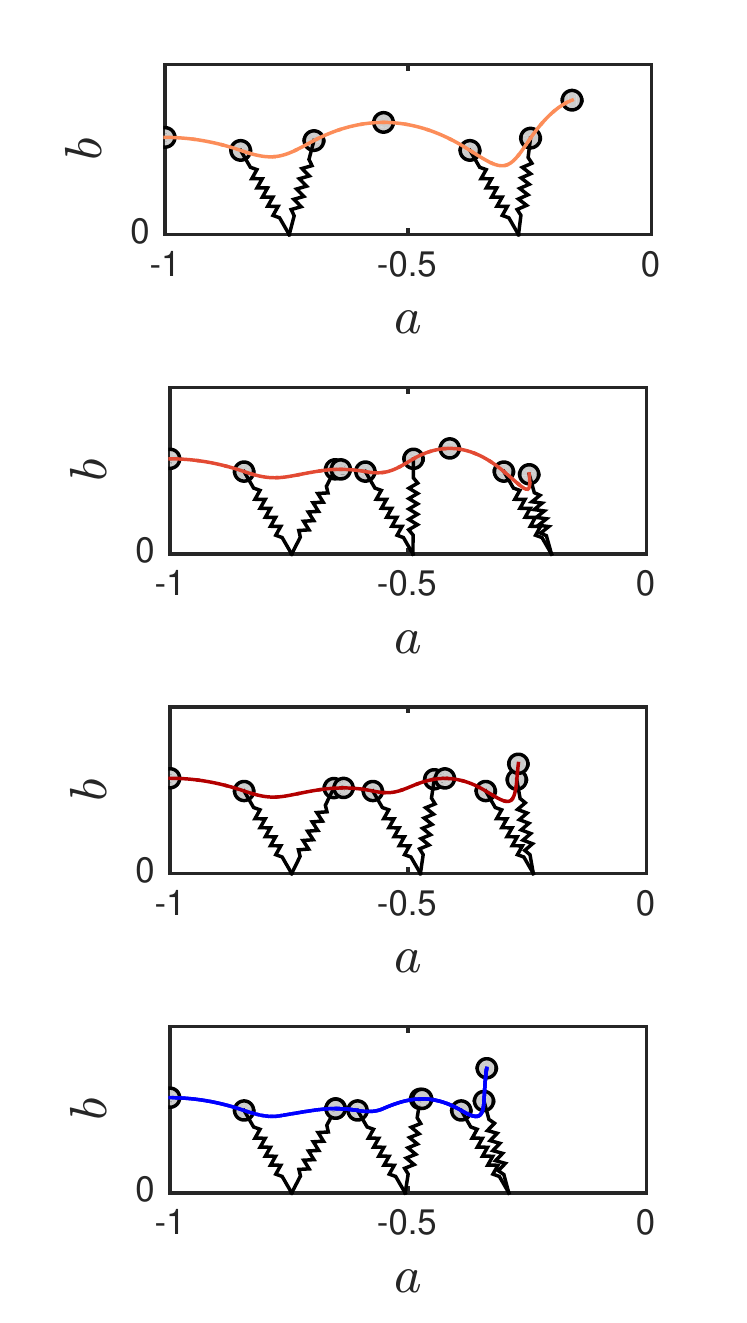}} \hspace{1cm}
	\end{minipage}
\centering
\caption{An illustration of the performance of our algorithm on active SLIP model.
The blue lines are the optimal control computed by GPOPS-II by iterating through all the possible transition sequences, and the red lines of various saturation are controls generated by our method. 
As the saturation increases the corresponding degree of relaxation increases between $2k=4$ to $2k=6$ to $2k=8$. 
Figure \cref{fig:SLIP:b} shows trajectories that maximize vertical displacement, where the optimal solution goes through 3 transitions; 
Figure \cref{fig:SLIP:V} shows trajectories that follow constant speed $v=0.1$, where the optimal solution goes through 8 transitions.
}
\label{fig:SLIP}
\end{figure}

The optimization problems are solved by our algorithm with degrees of relaxation $2k=4$, $2k=6$, and $2k=8$. 
For the sake of comparison, the same problems are also solved using GPOPS-II. 
Since GPOPS-II requires information about the transition sequence, we let the number of transitions to be less than 12, and iterate through all possible transition sequences with GPOPS-II. 
The results are compared in Figure \cref{fig:SLIP} and Table \cref{tab:SLIP:res}.

\begin{table}[htbp]
\caption{Results of the active SLIP}
\centering
\label{tab:SLIP:res}
\begin{tabular}{c|c|c|C{3cm}|C{3cm}}
\hline
    \multicolumn{2}{c|}{} & Computation time & Cost returned from optimization & Cost returned from simulation \\
\hline
    \multirow{4}{3.5cm}{Maximizing vertical displacement}
        & $2k=4$ & 42.1805[s] & -0.7003 & -0.5480 \\ \cline{2-5}
        & $2k=6$ & 722.5955[s] & -0.5773 & -0.5577 \\ \cline{2-5}
        & $2k=8$ & $1.2290 \times 10^{4}$[s] & -0.5754 & -0.5629 \\ \cline{2-5}
        & GPOPS-II & 1453.1083[s] & -0.5735 & N/A \\
\hline
    \multirow{4}{3.5cm}{Tracking constant-speed trajectory $a(t)=vt-1$ with $v=0.1$}
        & $2k=4$ & 50.9472[s] & 0.0422 & 0.38931 \\ \cline{2-5}
        & $2k=6$ & 835.4857[s] & 0.2107 & 0.31507 \\ \cline{2-5}
        & $2k=8$ & $9.1429 \times 10^{3}$[s] & 0.2165 & 0.31142 \\ \cline{2-5}
        & GPOPS-II & 844.3898[s] & 0.2657 & N/A \\
\hline
\end{tabular}
\end{table}

  \bibliography{SICON,refs}
  \bibliographystyle{IEEEtran}
\end{document}